\documentclass[12pt]{amsart}
\usepackage{amssymb,amsmath,amsthm,mathrsfs}
\usepackage{xspace}
\usepackage[dvipsnames]{xcolor}
\usepackage{tikz-cd}
\usepackage{esint}
\usepackage{libertinus}
\usepackage[libertine,vvarbb]{newtxmath}
\usepackage{microtype}
\usepackage[linktocpage=true,backref=page]{hyperref}
\hypersetup{
	colorlinks=true, 	%set true if you want colored links
	linktoc=all,     		%set to all if you want both sections and subsections linked
	linkcolor=RoyalPurple,  	%choose some color if you want links to stand out
	citecolor=Emerald,
}
\calclayout
\newtheorem{theorem}{Theorem}[section]
\numberwithin{theorem}{section} \numberwithin{equation}{section}
\newtheorem{corollary}[theorem]{Corollary}
\newtheorem{proposition}[theorem]{Proposition}
\newtheorem{lemma}[theorem]{Lemma}

\theoremstyle{definition}
\newtheorem{remark}[theorem]{Remark}
\newtheorem{assumption}[theorem]{Assumption}

\newtheorem{definition}[theorem]{Definition}
\newtheorem{conjecture}[theorem]{Conjecture}
\newtheorem{acknowledgment}[theorem]{Acknowledgment}

\makeatletter
\def\@setauthors{%
	\begingroup
	\def\thanks{\protect\thanks@warning}%
	\trivlist
	\centering\footnotesize \@topsep30\p@\relax
	\advance\@topsep by -\baselineskip
	\item\relax
	\author@andify\authors
	\def\\{\protect\linebreak}%
	\MakeUppercase{\authors}%
	\ifx\@empty\contribs
	\else
	\par\@setcontribs
	\@closetoccontribs
	\fi
	\endtrivlist
	\endgroup
}
\providecommand\@wraptoccontribs[2]{}
\newcommand\appendixby[1]{%
	\par\nobreak\vspace{-.5\baselineskip}%
	\begin{center}\footnotesize\MakeUppercase{By #1}\end{center}%
	\nobreak\@afterheading}
\newcommand\printaddresses{\@setaddresses\global\let\addresses\@empty}
\newcommand\appendixaddress[2]{%
	\begingroup\def\addresses{\address{}{#1}\email{}{#2}}\@setaddresses\endgroup}

\def\l@section{\@tocline{1}{0pt}{1pc}{1.5pc}{}}
\makeatother

\title[Reciprocal mirror map derivatives and Quantum Serre duality]{Reciprocal mirror map derivatives, euler-twisted Gromov-Witten invariants, and instanton counts via Quantum Serre duality}

\author{Irit Huq-Kuruvilla}
\email{huqkuruvilla\_irit@berkeley.edu}
\author{Michael T.  Schultz} 
\address{Dept.\!~of Mathematics, Virginia Tech, Blacksburg, VA 24060}
\email{michaelschultz@vt.edu}

\begin{document}
	\begin{abstract}
		We study mirror map derivatives as a means of relating the genus zero Gromov-Witten theories of an anticanonical Calabi-Yau hypersurface in a smooth, almost minimal Fano variety and the local Calabi-Yau formed from the total space of the Fano's canonical bundle with the euler-twisted theory of the ambient Fano. A precise relationship between reciprocal mirror map derivatives is obtained from the non-equivariant limit of quantum Serre duality for the euler-twisted theory of the ambient Fano as formulated by Iritani, Mann, \& Mignon. This allows us to express the generating series of 3-point euler-twisted invariants as the derivative of the local mirror map with respect to the closed one. In low dimensions, this series is shown to be a modular form of Eisenstein type which has arisen before in this area. We prove integrality of the twisted invariants from this identification. The regularized small quantum $D$-module of the Fano is reducible, and we show it is equivalent to an inhomogeneous Picard-Fuchs equation governed by this generating series. The reciprocal mirror map derivative is then shown to be a generating function of local 3-point invariants, which recovers some known results for local Calabi-Yau threefolds and proves a conjecture of the second named author and A. Malmendier in previous work on local fourfolds. Explicit calculations for $\mathbb{P}^2,\mathbb{P}^3$ and $\mathbb{P}^4$ illustrate these relationships, which includes the vanishing of positive-degree ordinary genus-zero invariants for anticanonical elliptic curves and K3 surfaces, the recovery of the quintic threefold’s	instanton numbers, and finally some formulas for genus-zero invariants of local projective spaces. 

	\end{abstract}
	\maketitle
	\tableofcontents
	\section{Introduction}
	\label{s-intro}
	\par The quantum Serre theorem as developed by Coates \& Givental \cite{MR2276766} describes a relationship between certain twisted Gromov-Witten invariants of a smooth projective variety $V$. The data of the twist is specified by the choice $(\mathbf{c},E)$ of an invertible characteristic class and holomorphic vector bundle $E \to V$. In genus zero, it can be expressed as an equality of twisted Lagrangian cones,
	\begin{equation}
		\label{eq-twisted_cones}
		\mathcal{L}^{(\mathbf{c}^*,E^*)} = \mathbf{c}(E)\mathcal{L}^{(\mathbf{c},E)} \, ,
	\end{equation}
	where $(\mathbf{c}^*,E^*)$ is the dual characteristic class and vector bundle, respectively. The former is defined by $\mathbf{c}(B)\mathbf{c}^*(B^*) = 1$ for any vector bundle $B$. The equality of twisted Lagrangian cones implies highly nontrivial relationships between the invariants defining both theories. 
	\par In parallel to ordinary Gromov-Witten (GW) theory of $V$, one may define twisted quantum products, twisted quantum $D$-modules, etc, which encode the same information as the twisted Lagrangian cone. Such twisted invariants were originally used to compute the genus zero invariants of smooth Calabi-Yau (CY) complete intersections $X \subset V$, where $V$ is a toric Fano variety and $X$ is carved out by the vanishing locus of a sufficient number of general sections of $E \to V$ \cite{MR1363062,MR2276766}, or for noncompact local Calabi-Yau varieties $Y = \mathrm{Tot}(E^*)$ \cite{MR1797015}, which are bundles $E^*$ such that the total space of the bundle has trivial canonical class, $K_Y \cong \mathcal{O}_Y$. In this context, $E^*$ is concave with a compatible torus action, and the invertible characteristic class is the equivariant Euler class $\mathbf{c} = \mathbf{e}_\lambda$. The non-equivariant limits of the $(\mathbf{e}_\lambda,E)$-twisted theory yield invariants of the vanishing locus of a generic smooth section of $E$, while the non-equivariant limits of the $(\mathbf{e}_\lambda^*,E^*)$-twisted theory yields invariants of the total space of $E^*$. 
	\par On the other hand, the equality of twisted Lagrangian cones (\ref{eq-twisted_cones}) does not naturally possess a non-equivariant limit, since the ordinary Euler class is not invertible. Iritani, Mann, \& Mignon \cite{MR3519131} provide a powerful extension of the original quantum Serre theorem that does admit a non-equivariant limit. This is done through first casting (\ref{eq-twisted_cones}) as a duality between twisted quantum $D$-modules by means of the Poincar\'e pairing and cup product with $\mathbf{c}(E)$. They show that this duality does admit a non-equivariant limit, and when applied to $\mathbf{e}_\lambda(E)$, yields a tight relationship between the ambient quantum $D$-module of $X^\vee \subset V$, the $(\mathbf{e},E)$-twisted quantum $D$-module of $X$, and the local quantum $D$-module of the total space $Y = \mathrm{Tot}(E^*)$. The maps defining these relationships are determined by the generating function of 2-point Euler-twisted invariants of $X$ and a cohomological shift by the first Chern class $c_1(E)$. Specifics are recalled in \S \ref{ss-QSD}.
	\par When restricted to the small quantum $D$-module of an almost minimal\footnote{Meaning the cohomology is concentrated in even degree and is as small as possible, see \S \ref{sss-prelim}. For $n=2$ we restrict the cohomology of the del Pezzo surface to the anticanonical line, and the Fano index is the degree.} Fano variety $V$ of dimension $n \geq 2$ and index $\mathrm{ind}(V) = \delta$, this form of QSD takes a very simple and elegant form that intertwines the closed and local mirror maps:
	\begin{equation}
		\label{eq-QSD_0}
		Q(q) = (-1)^{\delta}\exp(f(q)) \, ,
	\end{equation}
	where $q = \exp(2\pi i\tau)$, $Q=\exp(2\pi i \tau_{\mathrm{loc}})$ are the exponentiated closed and local mirror maps, and $f$ is a generating function of 2-point euler twisted invariants
	\begin{equation}
		\label{eq-2-point}
		f(q) = \log(q)+\sum_{d=1}^{\infty} q^{d}\left\langle \mathrm{P}^{n-1}, \widetilde{1}\right\rangle_{0,2, d}^{\left(\mathrm{e}, \,  K_V^{-1}\right)} \, 
	\end{equation}
	in terms of the positive primitive generator $\mathrm{P} \in H^2(V,\mathbb{Z})$. Thus QSD has a natural link to both closed and local mirror symmetry. It is our goal in this article to extract meaningful consequences from the equations above by utilizing mirror symmetry. We will see that the non-equivariant limit of QSD places all three GW theories on an equivalent level, and this perspective allows us to obtain integrality results related to (\ref{eq-2-point}) in the euler-twisted theory in low dimensions.
	\par Another such consequence is to recover some known results for instanton counts, or genus zero Gopakumar-Vafa (GV) invariants, of certain local Calabi-Yau (CY) threefolds. In particular, QSD makes it natural to study local CYs that arise as the total space $Y^\vee=\mathrm{Tot}(K_V)$ of the canonical line bundle of $V$, which is the simplest example of the situation described in the preceding paragraphs. A result of Bloch, Kerr, \& Vanhove \cite[Thm. 6.1]{MR3780269} on the local GW potential for $Y^\vee = Y^\vee_{\tau_{\mathrm{loc}}}$, and combined with results of Doran \& Kerr \cite[\S 5]{MR2851155} yields the following statement:
	\begin{theorem}
		\label{thm-intro}
		[Bloch, Kerr, Vanhove; Doran \& Kerr] Let $V$ be an almost minimal toric Fano surface corresponding to an LG model $\mathcal{W}$, and $E_t$ be the corresponding mirror family of elliptic curves to anticanonical elliptic curves $\check{E}_{\tau} \in |-K_V|$. Then the composition $\tau(\tau_{\mathrm{loc}})$ of mirror maps determines the local invariants of $Y^\vee$ as 
		\begin{equation}
			\label{eq-local_threefold_yukawa}
			\frac{\kappa^{\mathrm{loc}}}{2\pi i} \frac{d\tau}{d\tau_{\mathrm{loc}}} = \mathcal{Y}_{3}^{\mathrm{loc}}(\tau_{\mathrm{loc}}) = \kappa^{\mathrm{loc}} + \sum_{d=1}^{\infty} d^3 n_d Q^d \, ,
		\end{equation}
		where $\mathcal{Y}_{3}^{\mathrm{loc}}(\tau_{\mathrm{loc}})$ is the local Yukawa coupling of $Y^\vee$, $\kappa^{\mathrm{loc}} \in \mathbb{Q}$ and $n_d$ are the degree-$d$ genus zero GW invariants of $Y^\vee$.
	\end{theorem}
	\par This result is a rigorous justification of something that was expected by physicists, in particular proving a conjecture of Hosono \cite{MR2282969}. We will see that Theorem \ref{thm-local_3-point} can be obtained naturally from (\ref{eq-QSD_0}). From the physical perspective of local mirror symmetry for local CY threefolds, the third order derivative of the holomorphic prepotential is naturally expressed as the mirror map derivative in the left hand side of (\ref{eq-local_threefold_yukawa}), where $\tau$ is interpreted as the modulus of the mirror elliptic curve. Similar expressions appear elsewhere in the physics literature \cite{MR1486340,MR1468700}.
	\par On the other hand, it has long been known that there is a strong connection to elliptic modular forms for such local CY threefolds. In cases relevant to the context of Theorem \ref{thm-intro}, the closed mirror map $t(\tau)$ is a Hauptmodul and holomorphic period $\omega_0(\tau)$ of the elliptic curve is a modular form of weight-1 for some congruence subgroup of $\Gamma \subseteq \mathrm{SL}(2,\mathbb{Z})$, and the complex structure moduli space of the mirror family is the (genus zero) curve $\Gamma \backslash \mathbb{H}$. Summarized succinctly, one may say that the Picard-Fuchs operator of the mirror family admits a modular parameterization. Moreover, the $Q$-expansion in (\ref{eq-local_threefold_yukawa}) is expressible in quasi-modular forms for $\Gamma$ \cite{MR3273318}. Modularity thus plays an important role in this context, and for example allows one to extend the a priori analytically local calculations in mirror symmetry to a Zariski open set in the complex structure moduli space \cite{MR3251364}. It also has implications for the growth rate of local GW invariants \cite[Cor. 5.3]{MR2851155}. 
	\par One such connection was observed by Stienstra \cite{MR2282958}. He noticed that from computing Mahler measure associated to the modular families of elliptic curves $\mathcal{X} \to \mathbb{P}^1_t$ arising in local mirror symmetry as above, the reciprocal of a certain modular form $\mathcal{R}(\tau)$ of weight-3 of Eisenstein type that factorizes as 
	\begin{equation}
		\label{eq-R_factorization}
		\mathcal{R}(\tau) = \frac{1}{2\pi i} \frac{\omega_0(\tau)}{t(\tau)}\frac{dt}{d\tau}
	\end{equation}
	was proportional by a rational multiple to the instanton expansion in (\ref{eq-local_threefold_yukawa}). This was called an observation as opposed to a theorem, as he at the time had no proof. Notice as well, that since in local mirror symmetry the local mirror map is defined as 
	\begin{equation*}
		\tau_{\mathrm{loc}}(t) = \frac{1}{2\pi i}\int_{t_0}^{t} \frac{\omega_0(s)}{s} ds \, ,
	\end{equation*}
	with respect to a fixed branch of the logarithm defined in a punctured neighborhood of $t=0$, one recovers
	\begin{equation*}
		\mathcal{R}(\tau) =  \frac{d\tau_{\mathrm{loc}}}{d\tau} \, ,
	\end{equation*}
	which is the reciprocal mirror map derivative as in (\ref{eq-local_threefold_yukawa}) in the sense of the inverse function theorem. Theorem \ref{thm-intro} proves Stienstra's observation. 
	\par Equivalently, from the perspective developed in this article, Stienstra's observation is perhaps the first instance of what is now understood as the non-equivariant limit of QSD. What has not been as clearly understood is what significance, if any, an enumerative significance can be attributed to the $q$-series of the modular forms appearing in local mirror symmetry in the local threefold case. The novelty of understanding such results from the unified perspective of modularity and QSD, is that the $q$-series of (\ref{eq-R_factorization}) has a natural enumerative interpretation as the generating series of 3-point euler-twisted invariants obtained from (\ref{eq-2-point}) for the associated del Pezzo; moreover, modularity (in particular, the manifestation of (\ref{eq-R_factorization}) as a series of Eisenstein type) implies such 3-point invariants are integral. 
	\par Doran \& Kerr \cite{MR2851155} elucidated and generalized many of the structures underlying Stienstra's observation. For example, they provided an explanation of the appearance of Eisenstein series using motivic cohomology classes on Kuga varieties. They extended results to the case of anticanonical K3 surfaces in rank-1 (toric) Fano threefolds where similar Eisenstein series of weight-4 appear in the same form of (\ref{eq-R_factorization}). These modular forms played an important role Golyshev \& Zagier's proof of the gamma conjecture for rank-1 Fano threefolds \cite{MR3462676}, and are rooted in Golyshev's analysis of mirror symmetry for rank-1 Fano threefolds \cite{MR2306141}, where modular parameterizations of the associated Picard-Fuchs operators arises naturally. Doran \& Kerr also extended local mirror symmetry results to the local fourfold case by studying modular properties of rank-1 Fano threefolds. 
	\par This led the second named author, in previous work with A. Malmendier \cite{ms2024}, to explore how such reciprocal quantities play a role in local mirror symmetry for CY fourfolds arising as the total space of the canonical bundle of rank-1 Fano threefolds. They were able to recover from an expansion analogous to (\ref{eq-local_threefold_yukawa}) computed purely from the data of the mirror of an anticanonical quartic K3 surface in $\mathbb{P}^3$ the genus zero invariants of local $\mathbb{P}^3$ \cite[\S 4.2.2]{ms2024}. That led them to conjecture \cite[\S 4.3]{ms2024} the local fourfold analogue of Theorem ~\ref{thm-intro}: 
	\begin{conjecture}
		\label{conj-local_fourfolds}
		If $V$ is a rank-1 Fano threefold of index $\delta$, the $\mathsf{Q} := Q^\delta$ expansion of the reciprocal mirror map derivative
		\begin{equation*}
			 \frac{d\tau}{d\tau_{\mathrm{loc}}} = \frac{1}{\mathcal{R}(\tau_{\mathrm{loc}})}
		\end{equation*}
		is, up to rational normalization, a generating function for the genus zero invariants of the local CY fourfold $Y^\vee = \mathrm{Tot}(K_V)$. Moreover, this expression can be expressed in terms of quasi-modular forms.
	\end{conjecture}
	\par The precise statement is provided in \S \ref{ss-reciprocal_mirror_map_derivative}. The goal of this article began with the goal of proving this conjecture, and trying to deeply understand the significance of modularity in the low dimensional cases. We found that the key of understanding the results described here, and also generalizing to higher dimensional CY varieties, is the non-equivariant limit of the euler-twisted QSD as captured by the elegant Equation (\ref{eq-QSD_0}) above. The appearance of Eisenstein series in the rank-1 Fano threefold case also implies that the associated 3-point euler-twisted invariants are integral, exactly analogous to the case described above. The reader should notice that after taking the logarithm of this equation and differentiating with respect to $\tau$, one recovers exactly the mirror map derivative $\frac{d\tau_{\mathrm{loc}}}{d\tau}$. This observation ties the whole story together. Our goal for the article is to appropriately describe the setup in which (\ref{eq-QSD_0}) is natural, contextualize it in the context of mirror symmetry, and then use the machinery of GW theory to rigorously prove a generalization of the Conjecture \ref{conj-local_fourfolds}, from which it will follow as a corollary. The power of QSD is that it is formulated in a much greater generality, and in particular applies to non-toric Fano varieties, and to higher dimensional such varieties. We were thus able to also obtain results for local CY fivefolds, where a modular picture is in general impossible \cite[\S 10.6]{MR2851155}. We found that a good amount of exposition was useful for aligning all of these ingredients, which we do before the aforementioned proofs. 
	\par The layout of this paper is as follows. In \S \ref{ss-QSD} we recall relevant results on genus zero Gromov-Witten theory and Quantum Serre duality, following the presentation of Iritani, Mann, \& Mignon \cite{MR3519131}. In \S \ref{s-mirror} we then recall basic notions of closed and local mirror symmetry for CY varieties determined by an almost minimal Fano variety $V$. \S \ref{s-non_eq} is functionally the A-side of mirror symmetry and \S \ref{s-mirror} is functionally the B-side, from the perspective of the LG-CY correspondence. Although mirror symmetry for elliptic curves and K3 surfaces is in general simpler than the CY threefold case, in particular because the genus zero theory is more or less vacuous (for different reasons), QSD necessitates that we understand the situation carefully as (\ref{eq-QSD_0}) is a bridge between the closed and local theories through euler-twisted GW invariants. 
	\par \S \ref{s-twisted_invariants} contains our main results. In Theorem \ref{thm-3-point_euler-twisted}, we carefully spell out the observation above on differentiating the logarithm of (\ref{eq-QSD_0}), then show that the resulting generating function of 3-point euler-twisted invariants arises naturally as the inhomogeneity of an inhomogeneous Picard-Fuchs equation. It follows from analysis in \S \ref{ss-fano_ms} that this inhomogeneous Picard-Fuchs equation is essentially the full second structure connection of the ambient Fano variety. Combining with the modularity results in low dimensions described above, in Theorem \ref{thm-euler-twisted_modular} we conclude that the coefficients in the $q$-series of the modular form $\mathcal{R}(\tau)$ are 3-point euler-twisted invariants of the ambient Fano; on the other hand, modularity allows us to conclude that these invariants are integral. Our main result is Theorem \ref{thm-local_3-point}, which concludes that the mirror map derivative $\frac{d\tau}{d\tau_{\mathrm{loc}}}$ is always a generating function of \emph{particular} local 3-point invariants; for local threefolds and fourfolds they are the invariants described above. Corollary \ref{cor-local_fourfold_instantons} is our proof of Conjecture \ref{conj-local_fourfolds}.
	\par Then in \S \ref{s-projective}, we provide computational examples of quantum Serre duality. In particular, we compute low degree 3-point euler-twisted invariants of $\mathbb{P}^n$ for $n=2,3,4$, and use them to relate genus zero invariants of anticanonical and local CYs determined by $\mathbb{P}^n$ for low values of $n=2,3,4$. We show how the classical instanton counts (= genus-zero GV invariants) $N_d \in \mathbb{Z}$ of the quintic threefold in $\mathbb{P}^4$ can be expressed in terms of 3-point Euler-twisted invariants of $\mathbb{P}^4$, and in turn how these 3-point invariants on local $\mathbb{P}^4$ can be expressed in terms of the $N_d$. This latter computation extends one that appears in \cite[\S 6.5]{MR3519131}, which was the inspiration for us to carry out the computations in this section. Interestingly, such invariants also appear to be integral, though it is not obvious from the definition of their generating function that this should be true.
	\par Lastly, in Appendix \ref{app-quantum_serre}, for convenience we give some more background information on GW theory and on quantum Serre duality. 
	\begin{remark}
		\label{rmk-generalities}
		We operate from the framework of Iritani, Mann, \& Mignon \cite{MR3519131}, whose work extends to a more general class of varieties. In particular, the non-equivariant limits of the euler-twisted and local small quantum $D$-module described in \S \ref{ss-SQDM} depends only on $V$ smooth with nef anticanonical class. For simplicity, we take $V$ to be a smooth, almost minimal Fano variety, for which there is a well understood notion of mirror symmetry. 
	\end{remark}
	\begin{acknowledgment}
		The authors acknowledge the use AI in the preparation of this article in the following way: the authors used AI tools (mostly Claude's Opus models) to perform literature searches and general queries at certain stages of the project, to troubleshoot and overcome various \LaTeX ~ formatting issues that arose in the process of writing, and for copy editing purposes to find notational inconsistencies and grammar / spelling errors in the pre-final draft. All of the mathematical content and arguments in the article are original, as is the entire writing and editing of the article.
	\end{acknowledgment}
	\subsubsection{Preliminaries - almost minimal Fano varieties}
	\label{sss-prelim}
		\par  Let $V$ be a smooth variety of dimension $n \geq 2$. Unless otherwise stated, in this article we will assume that $V$ is \emph{almost minimal}. Recall that $V$ is said to be \emph{minimal} if the cohomology is as small as possible, that is $H^{2k}(V,\mathbb{Z}) \cong \mathbb{Z}$ and $H^{2k+1}(V,\mathbb{Z}) = 0$ for $k=0,\dots,n$; one says that $V$ is \emph{almost minimal} if is satisfies the previous condition except in the middle dimensional cohomology group, which for example allows for $V$ to admit complex structure deformations. We will generally take $V$ to be a smooth Fano $n$-fold, that is, a smooth, complex projective variety of dimension-$n$ with ample anticanonical class. For such a Fano variety of dimension $n \geq 3$, an important numerical invariant is the so-called index $\delta= \mathrm{ind}(V) \in \mathbb{N}_{\geq 1}$, which is defined as $\delta = \left[H^2(V,\mathbb{Z}) \, : \, \mathbb{Z}\mathrm{P}\right]$, where $\mathrm{P} \in H^2(V,\mathbb{Z})$ is a primitive positive generator. We use the notation $H=c_1(-K_V) = c_1(V) \in H^2(V,\mathbb{Z})$ to be the anticanonical class, and thus we have ~$H = \delta\mathrm{P}$. 
		\par For example, every projective space $\mathbb{P}^n$ is a minimal Fano $n$-fold of index $\delta=n+1$. Our interest in Fano varieties is in large part that they are often a vehicle to carry Calabi-Yau (CY) divisors in their complete anticanonical linear system $|-K_V|$. In the context of this article, a CY variety of dimension $m > 1$ is a (smooth) variety with trivial canonical class and $h^{i,0}=0$ for $i = 1,\dots,m-1$. Such a variety may be closed (resp. open), and hence projective (resp. quasiprojective). A 2-dimensional CY variety is a K3 surface; a 1-dimensional CY variety is an elliptic curve. A quasiprojective CY variety is called a \emph{local} CY.
		\par We will mainly, but not exclusively, be concerned with low-dimensional Fano varieties, of dimension $\dim(V)=2,3,4$. A del Pezzo surface is a Fano variety of dimension two; they are (almost) classified by a single numerical invariant, the degree $(K_V,K_V) = m =1,\dots,9$. Up to isomorphism, a del Pezzo surface is either the ruled surface $\mathbb{P}^1\times \mathbb{P}^1$ of degree $8$, or is the surface $S_m$ of degree $m \geq 1$ realized as the blowup of $\mathbb{P}^2$ in $9-m$ points in general position. The only (almost) minimal del Pezzo surface is $S_9 = \mathbb{P}^2$. Because of this situation in dimension 2, we will take the following convention that modifies the previous notion: an \emph{almost minimal} del Pezzo surface $S$ is the pair $(S,\langle -K_S \rangle)$ of the surface and the degree-2 cohomology restricted to the anticanonical line; here 
		\begin{equation*}
			\langle -K_S \rangle = \mathrm{span}_{\mathbb{C}}\{-K_S\}\subseteq H^2(S,\mathbb{C}) \, .
		\end{equation*} 
		\par Increasing the dimension by one, basic properties of all Fanos (see e.g., \cite[Thm. 1.3]{MR4922010}) and Hodge theory imply that any Fano threefold of Picard rank-1 is either minimal or almost minimal. Recall that each deformation class of rank-1 Fano threefolds is determined by two positive integer invariants, the index $\delta = \mathrm{ind}(V)$, and the level $N =  \frac{1}{2\delta^2} \left\langle H^3,[V] \right\rangle$. There are precisely seventeen deformation classes of smooth rank-1 Fano threefolds as determined by Iskovskikh \cite{MR0463151,MR0503430}, which correspond to the following possible combinations of invariants $(N, \delta)$:
	\begin{align}
		\label{eq-degree_level_invariants}
		(1,1), \; \; (2,1), \; \; (3,1), \; \; (4,1), \; \; & (5,1), \; \; (6,1), \; \;  (7,1), \; \; (8,1), \; \; (9,1), \; \; (11,1), \nonumber \\
		(1,2), \; \; (2,2), \; \; (3,2), \; \; & (4,2), \; \; (5,2), \; \; (3,3), \; \; (2,4)\, . 		
	\end{align}
	 More generally, there are 105 deformation classes of smooth Fano threefolds, classified and constructed by Mori \& Mukai \cite{MR641971,MR1102273,MR1969009}. A relatively uniform description of these varieties is given for example by Coates et al \cite{MR3470714}. More generally, for any Fano $V$, every class in $H^2(V,\mathbb{Z})$ is algebraic, so in particular one has $\mathrm{Pic}(V) \cong H^2(V,\mathbb{Z})$ and thus any $V$ almost minimal has Picard rank $\rho(V)=1$. For simplicity, we will often assume that $\delta = 1$. 
	\par We will study divisors in the complete anticanonical linear system of $V$. For $V = \mathbb{P}^n$, these are global sections of the line bundle $-K_{\mathbb{P}^n} \cong \mathcal{O}(n+1)$ which have of course been widely studied, as they are a natural example of CY varieties. For low dimensions $\dim(V) = 2,3$, it is well-known that the anticanonical linear system is non-empty. For example, we have the following result which plays an important role in the context of this article.
	\begin{proposition}
		\label{prop-fano_threefold_divisor}
		[Prokhorov, \cite[Prop. 1.10]{MR4922010}] Let $V$ be a Fano threefold, and $X^\vee \subset V$ be a reduced, irreducible component of a divisor $D \in |-K_V|$. Then only one of the following possibilities holds:
		\begin{enumerate}
			\item[(i)] $X^\vee$ is a birationally ruled surface;
			\item[(ii)] $X^\vee=D$ and $X^\vee$ is a K3 surface, possibly with du Val singularities.
		\end{enumerate}
	\end{proposition}
	Property (ii) follows in the smooth case from adjunction, which shows the $K_{X^\vee} \cong \mathcal{O}_{X^\vee}$. There are other results characterizing effective anticanonical divisors of Fano fourfolds and fivefolds. It was shown by H\"oring \& Voison \cite{MR2918165} that any smooth Fano fourfold has an effective anticanonical divisor with isolated canonical singularities; this was sharpened to terminal singularities by Heuberger \cite{heuberger2015}. Later, results of H\"oring \& \'Smiech showed that any smooth Fano fivefold has an effective canonical divisor, and if it is  reduced, has at worst canonical singularities \cite{MR4060366}. More generally, results on the non-emptiness of $|-K_V|$ depend heavily on the index $d$:
	\begin{proposition}
		\label{prop-nonempty}
		[Prokhorov, \cite[Thm. 1.15]{MR4922010}]
		Suppose that $V$ is a smooth Fano $n$-fold of index $\delta \geq n-2$. Then $|-K_V| \neq \emptyset$.
	\end{proposition} 
	At present, the general non-emptiness of $|-K_V|$ for higher dimensional Fano varieties of low index is not known to the best of our knowledge. This leads us to the following critical assumption that generalizes Proposition \ref{prop-fano_threefold_divisor} that we will utilize throughout the remainder of the article. 
	\begin{assumption}
		\label{ass-nonempty}
		Let $V$ be a smooth, almost minimal Fano $n$-fold. Then we assume that $|-K_V| \neq \emptyset$, and a very general, reduced, irreducible divisor $X^\vee \in |-K_V|$ is a CY variety, possibly with canonical (Gorenstein) singularities. In the singular case, we freely assume that we may work with a minimal crepant resolution, which by an abuse of notation we also write as $X^\vee$.
	\end{assumption}
	\par Given our assumptions, associated to a Fano $n$-fold $V$ are two natural CY varieties of respective dimension $n-1, \, n+1$:  $X^\vee \in |-K_V|$, a very general smooth anticanonical divisor, and $Y^\vee = \mathrm{Tot}(K_V)$, the total space of the canonical line bundle. The quasiprojective CY variety $Y^\vee$ is a fundamental nontrivial example of a local CY, and will be a principal object of study in the remainder of the article. By Noether-Lefschetz type results, $\rho(X^\vee)=1$, and by the Hard Lefschetz theorem, $X^\vee$ may only acquire nontrivial cohomology in middle dimension. The compactly supported cohomology groups of $Y^\vee$ satisfy the same conditions, by more basic considerations for cohomology of the total space of vector bundles.
	\begin{remark}
		\label{rmk-Q_Fano}
		There are singular Fano varieties with $|-K_V| = \emptyset$, see e.g. \cite[Ch. 16]{MR1798978}. Another explicit example due to Johnson \& Koll\'ar \cite{MR1821068} shows that the $\mathbb{Q}$-del Pezzo surface $V$ of degree $8k+4$ in the weighted projective space $\mathbb{P}(2,2k+1,2k+1,4k+1)$ has $|-K_V| = \emptyset$. However, a celebrated result of Birkar \cite{MR3997127} shows that for $V$ a klt $\mathbb{Q}$-Fano $n$-fold, there is an integer $m \in \mathbb{N}_{\geq 1}$ depending only on $n = \dim(V)$ such that the anti-pluricanonical system $|-mK_V| \neq \emptyset$.
	\end{remark}
	\subsubsection{Notation}
	\label{sss-notation}
	\par We use the following notation consistently throughout the article. \\
	\vspace{5mm}
	\noindent
	\begin{tabular}{@{}l p{0.78\textwidth}@{}}
		$V$ & a smooth almost minimal Fano variety of dimension $n \geq 2$ over $\mathbb{C}$ \\
		$\mathrm{P}$ & the positive, primitive generator $H^2(V,\mathbb{Z})\cong \mathbb{Z}$ \\
		$H$ & the anticanonical class $H = -K_V = \delta \mathrm{P}$, where $\delta$ is the Fano index of $V$ \\
		$X^\vee = X^\vee_\tau$ & a smooth anticanonical CY variety in $V$, with choice of complexified symplectic moduli $\tau$ \\
		$\mathcal{X}^\vee$ & the family of $X^\vee_\tau$ parameterized by $\tau$ \\
		$X = X_t$ & a mirror variety of $X^\vee_\tau$ with choice of complex structure $t$ \\
		$\mathcal{X}$ & the family of $X_t$ parameterized by $t$ \\
		$\tau(t)$ & the closed mirror map for $\mathcal{X}^\vee$ and $\mathcal{X}$ \\
		$q(t)$ & the exponentiated closed mirror map $q(t) = \exp(2\pi i\tau(t))$ \\
		$Y^\vee = Y^\vee_{\tau_{\mathrm{loc}}}$ & the local CY variety $\mathrm{Tot}(K_V)$ with choice of complexified symplectic moduli $\tau_{\mathrm{loc}}$ \\
		$Y = Y_t$ & the mirror of $Y^\vee_{\tau_{\mathrm{loc}}}$ \\ 
		$\tau_{\mathrm{loc}}(t)$ & the local mirror map for $Y^\vee$ and $Y$ \\
		$Q(t)$ & the exponentiated local mirror map $Q(t) = \exp(2\pi i\tau_{\mathrm{loc}}(t))$ \\
		$\mathrm{Eff}(V)$ & the set of effective curves classes in $H_2(V,\mathbb{Z})$ \\
		$\langle \psi ,d \rangle$ & the pairing $\int_d \psi$ between $\gamma \in H^2(V)$ and $d \in H_2(V)$. \\
	\end{tabular}
	\section{Non-equivariant limit of quantum Serre duality}
	\label{s-non_eq}
	\par In this section we briefly review aspects of the non-equivariant limit of quantum Serre duality (QSD) that are pertinent to our analysis. More specifics for $V$ an almost minimal Fano variety are given in \S \ref{ss-fano_ms}. For convenience, a deeper look into the technical structure of QSD that is beyond the scope of this article is contained in Appendix \ref{app-quantum_serre}. Gromov-Witten invariants and quantum $D$-module are reviewed in Appendix \ref{ss-GW_invariants} and twisted invariants are reviewed in Appendix \ref{ss-twisted_invariants}. The formulation of QSD in terms of twisted Lagrangian cones is reviewed in Appendix \ref{ss-lagrangian_cone}, and the non-equivariant limit of the euler-twisted theory in Appendix \ref{ss-euler_twisted}.
	\subsection{Small quantum $D$-module}
	\label{ss-SQDM}
	\par Let $V$ be a smooth projective variety of dimension $n$ over $\mathbb{C}$. Choose a homogeneous basis $T_0=1,\ldots,T_N \in H^\bullet(V,\mathbb{C})$ and its
	Poincar\'e dual basis $T^0,\ldots,T^N$, so that
	$(T_\alpha,T^\beta)=\delta_\alpha^\beta$. Write  
	$$ \tau=\sum_{\alpha=0}^N t^\alpha T_\alpha,\qquad
	\partial_\alpha=\frac{\partial}{\partial t^\alpha}.$$ 
	Quantum cohomology of $V$, denoted by $QH(V) := (H^\bullet(V,\mathbb{C}), \star_\tau)$ is a deformation of the usual cup product on the cohomology ring, providing a family of multiplications parameterized by $\tau \in H^\bullet(V,\mathbb{C})$. It is also called the big quantum cohomology ring, and is well-known to have the structure of a Frobenius algebra \cite{MR1397274,MR1702284}. Here $\star_\tau$ is the quantum product (\ref{eq-quantum_product}), which encodes genus zero Gromov-Witten (GW) invariants of $V$. GW invariants of $V$ are equivalently encoded in a triple 
	\begin{equation}
		\label{eq-QDM}
		\operatorname{QDM}(V):=(\mathcal F,\nabla,S)
	\end{equation} 
	called the A-side quantum $D$-module. Here $\mathcal F := H^\bullet(V,\mathbb{C}) \otimes_{\mathbb C}\mathbb C[z][[q,\tau]]$, which is thought of as a trivial holomorphic vector bundle with fibre $H^\bullet(V,\mathbb{C})$ over $H^\bullet(V,\mathbb{C}) \times \mathbb{C}_z$, and $\nabla$ is a flat connection on $\mathcal{F}$ called the quantum connection. The connection $\nabla$ is defined in terms of $\star_\tau$ as 
	\begin{equation}
		\label{eq-quantum_connection}
		\nabla_\alpha
		:=
		\partial_\alpha+\frac{1}{z}(T_\alpha\bullet_\tau),
	\end{equation}
	and is flat, moreover, it preserves the extension of the Poincar\'e pairing $S(u,v):=\bigl(u(-z),v(z)\bigr)$ to $\mathcal{F}$. $\nabla$ may be extended to a flat connection in the $z$-direction by incorporating quantum multiplication with the so-called euler vector field; it may also be thought of as a pencil of flat connections $\nabla^z$ by projecting onto the first factor $\mathcal{F} \to H^\bullet(V,\mathbb{C})$. 
	\par The quantum connection $\nabla^z$ has a fundamental solution  $L(\tau,z) = \mathrm{id} + O(z^{-1})$ which is defined in terms of descendant invariants (\ref{eq-fundamental_solution}) and acts isometrically on the Poincar\'e pairing. One then defines the $J$-function as
	\begin{equation}
		\label{eq-J-function}
		J(\tau,z):=zL(\tau,z)^{-1}1
		=z1+\tau+O(z^{-1}) \, .
	\end{equation}
	It is the primary generating function of all (descendant) genus zero GW invariants. 
	\par The small quantum $D$-module $\mathrm{SQDM}(V)$ is the restriction of the data above to the divisorial directions in quantum multiplication $\tau \in H^2(V,\mathbb{C})$, and the resulting algebraic structure is the small quantum cohomology ring $\mathrm{QH}_{\mathrm{sm}}(V)$. In the case that $V$ is an almost minimal variety, the small quantum connection reduces to a linear ODE after setting $z=1$. It is irregular, with precisely one regular singularity at $t=0$ and one irregular singular point at $t = \infty$ (we follow the conventions of Golyshev \& Zagier \cite[\S 1.2]{MR3462676}). The fundamental solution at the irregular singular point is the small $L$-function.
	\subsection{Non-equivariant limit of Euler-twisted theory}
	\label{ss-QSD}
	\par The key insight of Givental \cite{MR1653024}, Givental \& Coates \cite{MR2276766} is that in the case of a subvariety $X^\vee \subset V$ determined by the vanishing of some generic sections of a holomorphic vector bundle $E \to V$, that not only should the $\mathrm{QDM}(X^\vee)$ be naturally related to $\mathrm{QDM}(V)$, but also there is an induced quantum $D$-module $\mathrm{QDM}(Y^\vee)$, where $Y^\vee = \mathrm{Tot}(E^*)$ is the total space of the dual bundle. It follows that there is an associated local quantum connection $\nabla^{\mathrm{loc}}$ \cite[\S 3.5]{MR3519131}. QSD makes this relationship precise. It is naturally formulated in the language of twisted GW invariants and twisted Lagrangian cones, the latter of which naturally encodes the (twisted) $J$-function. It requires the data of an invertible characteristic class $\mathbf{c}$, and reads in genus zero as an equality
		\begin{equation}
		\label{eq-twisted_cones0}
		\mathcal{L}^{(\mathbf{c}^*,E^*)} = \mathbf{c}(E)\mathcal{L}^{(\mathbf{c},E)} \, ,
	\end{equation}
	where $(E^*,\mathbf{c}^*)$ is the vector bundle and characteristic class dual of $(E,\mathbf{c})$. 
	\par In the geometric case of a subvariety in $V$, the characteristic class is taken as the equivariant euler class $\mathbf{c} = \mathbf{e}_\lambda$. Here the action is taken to be the natural $\mathbb{C}^*$-action scaling the fibres of $E$, extended to all of $Y^\vee$ by acting trivially on $V$ (or, compatibly in the case that $V$ also admits a torus action). Then if $E$ is convex, the non-equivariant limits of the $(\mathbf{e}_\lambda,E)$-twisted theory yield GW invariants of $X^\vee$, while the non-equivariant limits of the $(\mathbf{e}_\lambda^*,E^*)$-theory yield GW invariants of $Y^\vee$. However, (\ref{eq-twisted_cones0}) does not admit a non-equivariant limit since $\mathbf{e}_\lambda$ is not invertible. The work of Iritani, Mann, \& Mignon \cite{MR3519131} remedies this by constructing non-equivariant limits of the associated twisted quantum $D$-modules and showing that they are related in a precise manner as described shortly.
	\par We are particularly interested in the context of this article in the geometry of anticanonical Calabi-Yau (CY) hypersurfaces $X^\vee$ in an almost minimal Fano variety $V$ of dimension $n \geq 2$. Particulars of these varieties are recalled in \S \ref{sss-prelim}. Then in the context of the discussion above, we take $E = K_V^{-1} = \wedge^n TV$ as the anticanonical line bundle. Then $K_V^{-1}$ is convex, and the desired duality should relate to the so-called local CY variety $Y^\vee = \mathrm{Tot}(K_V)$. Such a situation is familiar to both geometers and physicists since the beginning of the subject of mirror symmetry and the Givental's fundamental work on the algebrogeometric incarnation of GW theory \cite{MR1653024,MR1408320} in studying the quintic threefold in $V=\mathbb{P}^4$, whence $K_V^{-1} = \mathcal{O}(5)$. In the context of mirror symmetry, it is natural to relate the quantum $D$-modules of $X^\vee$, $V$, and $Y^\vee$ together. We will use the following machinery to do so in \S \ref{s-twisted_invariants}.
	\begin{theorem}
		\label{thm-IMM_1} \cite[Theorem 1.1]{MR3519131}
		There exists a typically nonlinear map $f_\circ : H^\bullet(V) \to ~H^\bullet(V)$ such that:
		\begin{itemize}
			\item[(1)] The twisted quantum $D$-modules $\mathrm{QDM}_{(K_V^{-1},\mathbf{e}_\lambda)}$ and $f_\circ^*\mathrm{QDM}_{(K_V,\mathbf{e}^*_\lambda)}$ are dual with respect to the Poincar\'e pairing, and
			\item[(2)] The map $\mathrm{QDM}_{(K_V^{-1},\mathbf{e}_\lambda)} \to f_\circ^*\mathrm{QDM}_{(K_V,\mathbf{e}^*_\lambda)}$ defined by $\alpha \mapsto \mathbf{e}_\lambda(K_V) \cup \alpha$ is a morphism of twisted quantum $D$-modules.
		\end{itemize}
	\end{theorem}
	\par Here the map $f_\circ$ is given in (\ref{eq-equivariant_map}), and we have $\mathbf{e}_\lambda(-K_V) = c_1(V) + \lambda$. Crucially, the map $f_\circ$ admits a non-equivariant limit $f := \lim_{\lambda \to 0} f_{\circ}$ \cite[Lem. 3.2]{MR3519131}; it is given explicitly in terms of 2-point euler-twisted invariants as 
	\begin{equation}
		\label{eq-twisted_gf0}
		f(q) = \log(q)+\sum_{d=1}^{\infty} q^{d}\left\langle \mathrm{P}^{n-1}, \widetilde{1}\right\rangle_{0,2, d}^{\left(\mathrm{e}, \,  K_V^{-1}\right)} \, 
	\end{equation}
	where $q$ is the Novikov variable on $QH(V)$. Then we have the following. 
	\begin{corollary}
		\label{cor-IMM_2} \cite[Cor. 1.2]{MR3519131}
		Let $V$ be an almost minimal Fano variety and $\mathbf{e}$ the non-equivariant euler class. Let $h : H^\bullet(V,\mathbb{C}) \to H^\bullet(V,\mathbb{C})$ be the map given by $h(\tau) = \tau + ~\pi \sqrt{-1}c_1(V)$ and $f : H^\bullet(V,\mathbb{C}) \to H^\bullet(V,\mathbb{C})$ be the map (\ref{eq-twisted_gf0}). Then:
		\begin{itemize}
			\item[(1)] The twisted quantum $D$-modules $\mathrm{QDM}_{(\mathbf{e},K_V^{-1})}(V)$ and $(h \circ f)^*\mathrm{QDM}(Y^\vee)$ are dual to each other, and
			\item[(2)] If $X^\vee \in |-K_V|$ and $\iota : X^\vee \to V$ is the inclusion map such that the Poincar\'e pairing on $H^\bullet_{\mathrm{amb}}(X^\vee) := \mathrm{Im}(\iota^*H^\bullet(V) \to H^\bullet(X^\vee))$ is nondegenerate, then the morphism $\mathbf{e}(-K_V) : \mathrm{QDM}_{(\mathbf{e},K_V^{-1})}(V) \to (h \circ f)^*\mathrm{QDM}(Y^\vee)$ factors through the ambient part of the quantum $D$-module $\mathrm{QDM}_{\mathrm{amb}}(X^\vee)$ as
			\begin{equation}
			\label{eq-quantum_serre_factorization}
				\begin{tikzcd}[column sep=large]
					\mathrm{QDM}_{(\mathbf{e},K_V^{-1})}(V)
					\arrow[r, "\mathbf{e}(-K_V)\cup"]
					\arrow[d, two heads, "\iota^*"']
					& (h\circ f)^*\,\mathrm{QDM}(Y^\vee) \\
					(\iota^*)^*\,\mathrm{QDM}_{\mathrm{amb}}(X^\vee)
					\arrow[ur, hook, "\iota_*"]
					&
				\end{tikzcd}
			\end{equation}
		\end{itemize}
	\end{corollary}
	\par The diagram (\ref{eq-quantum_serre_factorization}) descends to the analogous symbols with small quantum $D$-modules \cite[\S 5.5]{MR3519131}. Given our assumption that $V$ is almost minimal, $\mathrm{SQDM}_{\mathrm{amb}}(V) \cong \mathrm{SQDM}(V)$. Here one has 
	\begin{equation}
		\label{eq-SQDM_twisted}
		\begin{aligned}
			\operatorname{SQDM}_{\left(\mathbf{e}, K_{\mathrm{V}}^{-1}\right)}(V) &:=\left(H^{\bullet}(V) \otimes \mathcal{O}_{U \times \mathbb{C}_z}, \nabla^{\mathrm{eu}}, S_{\mathrm{eu}}\right) \, , \\
			\operatorname{SQDM}\left(Y^\vee\right) &:=\left(H^{\bullet}(V) \otimes \mathcal{O}_{U \times \mathbb{C}_z}, \nabla^{\mathrm{loc}}\right) \, , \\
			\operatorname{SQDM}_{\mathrm{amb}}(X^\vee) &:=\left(H_{\mathrm{amb}}^{\bullet}(X^\vee) \otimes \mathcal{O}_{U \times \mathbb{C}_z}, \nabla^Z, S_Z\right) \, ,
		\end{aligned}
	\end{equation}
	and $U \subseteq H^2(X^\vee,\mathbb{C})$ defined  to be some analytic region where $q \to 0$ in a meaningful way; such a point is called the \emph{large volume limit}. We will discuss this in much more detail in \S \ref{s-mirror}. The relevant data for the euler-twisted connection in (\ref{eq-SQDM_twisted}) is given in Appendix \ref{ss-twisted_invariants}.
	\par In spite of the technical machinery above, the expression of the non-equivariant QSD on the small quantum $D$-module is extremely clean and elegant. We delay the expression until Lemma \ref{lem-exp_mirror_maps} after we have had a chance to fully discuss mirror symmetry in \S \ref{s-mirror}. First, we will recall the A-side machinery relevant to the B-side of mirror symmetry: $I$-functions.
	\subsection{$I$-functions}
	\label{ss-I_functions}
	\par Given a variety $W$, the $I$-function, when it exists, is a regularized version of the $J$-function that lies on the same Lagrangian cone. The (small) $I$-function is usually determined by some extrinsic presentation of the variety $W$; for example, when $W$ is (a subvariety in the sense of \S \ref{ss-QSD} of) a toric Fano variety $V$, the $I$-function has the form of a hypergeometric series. The (small) $I$-function itself is not equal to the (small) $J$-function, but is related by a change of variables known as the mirror map. More specifically there is a map $\tau(q)$ on the Novikov ring such that
	\begin{equation*}
		I_W(q,z) = J_W(\tau(q),z) \cdot \exp(f(q)/z) \,
	\end{equation*}
	where $f(q)$ is a certain scalar correction. The relation here allows one to read off the mirror map, presuming that both sides are known sufficiently. When $W$ is CY, the mirror map is often highly transcendental.
	\par Iritani, Mann, \& Mignon  \cite[Def. 5.11, Lem. 5.4]{MR3519131} showed that for the $\mathrm{SQDM}$ in (\ref{eq-SQDM_twisted}), both an euler-twisted $I$-function and local $I$-function can be defined. They are given as follows. Take $W=V$ as an almost minimal Fano, $\tau \in H^2(V,\mathbb{C})$. Then they defined
	\begin{equation}
		\label{eq-I-functions}
		\begin{aligned}
			I_\alpha^{\mathrm{eu}}\left(\tau, z\right) &:=e^{\tau / z} \sum_{d \in \operatorname{Eff}(V)} N_{\alpha, d}(z) \exp(\langle \tau,d \rangle)\prod_{k=1}^{\langle c_1(V),d \rangle}(c_1(V)+k z) \, , \\
			I_\alpha^{\mathrm{loc}}\left(\tau, z\right) &:=e^{\tau / z} \sum_{d \in \operatorname{Eff}(V)} N_{\alpha, d}(z) \exp(\langle \tau,d \rangle)\prod_{k=0}^{\langle c_1(V),d \rangle-1}(-c_1(V)-k z) \, , \\
			N_{\alpha, d}(z) &:= \begin{cases}\sum_{\beta=0}^s\left\langle T_\alpha, \frac{T_\beta}{z-\psi}\right\rangle_{0,2, d} T^\beta & \text { if } d \neq 0 ; \\ T_\alpha & \text { if } d=0 .\end{cases}
		\end{aligned}
	\end{equation}
	where we use the notation of the 2-point invariants as in Appendix \ref{ss-GW_invariants}. We shall be interested in the specific $I$-functions in (\ref{eq-I-functions}) corresponding to $\alpha=0$. Then one should interpret the expressions here as series expansions
	\begin{equation}
		\label{eq-I-expansions}
		\begin{aligned}
			& I_0^{\mathrm{eu}}\left(\tau, z\right)=I_0^{\mathrm{eu}}\left(\tau\right) \cdot 1+I_1^{\mathrm{eu}}\left(\tau\right) \cdot z^{-1}+I_2^{\mathrm{eu}}\left(\tau\right) \cdot z^{-2}+O\left(z^{-3}\right), \\
			& I_0^{\mathrm{loc}}\left(\tau, z\right)=1+I^{\mathrm{loc}}_1\left(\tau\right) \cdot z^{-1}+I_2^{\mathrm{loc}}\left(\tau\right) \cdot z^{-2}+O\left(z^{-3}\right),
		\end{aligned}
	\end{equation}
	where $I_0^{\mathrm{eu,loc}}(\tau)$ is a scalar valued function, $I_1^{\mathrm{eu,loc}}(\tau)$ is $H^2(V,\mathbb{C})$ valued, etc. From this presentation, the euler-twisted and local mirror maps are defined as 
	\begin{equation}
		\label{eq-I-mirror_maps}
		\tau_{\mathrm{eu}} := \frac{I^{\mathrm{eu}}_1(\tau)}{I^{\mathrm{eu}}_0(\tau)} \, , \quad \tau_{\mathrm{loc}} := \frac{I^{\mathrm{loc}}_1(\tau)}{1}  \equiv I^{\mathrm{loc}}_1(\tau) \, .
	\end{equation}
	From the context of Corollary \ref{cor-IMM_2}, the assumption of $W=V$ is an almost minimal Fano variety yields the following result.
	\begin{theorem}
		\label{thm-minimal_I-functions}
		For $V$ an almost minimal Fano variety and $X_\tau^\vee \in |-K_V|$ an anticanonical CY variety, the small $I$-function of $X^\vee_{\tau}$ is the euler-twisted $I$-function $I_\alpha^\mathrm{eu}(\tau,z)$ of $V$. For the local CY variety $Y^\vee_\tau$, the local $I$-function is given by $I_\alpha^\mathrm{loc}(\tau,z)$.
	\end{theorem}
	\begin{proof}
		The almost minimality assumption on $V$ guarantees that the $\mathrm{(S)QDM}$ morphisms in Corollary (\ref{cor-IMM_2}) are surjective. By \cite[Lem. 5.14]{MR3519131}, the expressions in (\ref{eq-I-functions}) agree with the appropriate fundamental solutions of the $\mathrm{SQDM}$ in (\ref{eq-SQDM_twisted}) up to composition with the mirror maps $\tau_{\mathrm{eu}}, \tau_\mathrm{loc}$, in (\ref{eq-I-mirror_maps}) respectively.
	\end{proof}
	\par The following more general results are elementary but will be useful later to point to. We do not claim any originality of their statements or proofs. 
	\begin{lemma}
		\label{lem-3-point_I-function}
		Let $W$ be a smooth, almost minimal variety of dimension $n \geq 2$ and Picard rank-1 with $-K_W \geq 0$. Suppose that $\mathrm{SQDM}(W)$ admits an $I$-function $I_W(\tau,z)$, written as 
		\begin{equation}
			\label{eq-I-function_expansion}
			I_W(q,z) = I_0(q) \cdot 1 + I_1(q) \cdot z^{-1} + I_2(q) \cdot z^{-2} + O(z^{-3}) \, .
		\end{equation}
		Here $I_0(q)$ is a scalar valued function, $I_1(q)$ is $H^2(W,\mathbb{C})$ valued, etc., and $q$ is the Novikov variable. Define $\tau=\tau(q)$ by
		\begin{equation}
			\label{eq-tau_I-function}
			\frac{I_1(q)}{I_0(q)}  =: \tau(q)
		\end{equation}
		and $\kappa = \langle \mathrm{P}^n,[W] \rangle$. Then the 3-point invariants $\langle \langle \mathrm{P},\mathrm{P},\mathrm{P}^{n-2}\rangle \rangle$ are determined by the second order derivative
		\begin{equation}
			\label{eq-sec_derivative_I-function}
			\frac{d^2}{d\tau^2}\left(\frac{I_2(\tau)}{I_0(\tau)}\right) \,  = \langle \langle \mathrm{P},\mathrm{P},\mathrm{P}^{n-2}\rangle \rangle = \kappa +  \sum_{d=1}^{\infty} d^2 \langle \mathrm{P}^{n-2} \rangle_{0,1,d} q^d \, . 
		\end{equation}
		Moreover, the small quantum cohomology ring $QH_{\mathrm{sm}}(W)$ decomposes as
		\begin{equation}
			\label{eq-QH_decomp}
			QH_{\mathrm{sm}}(W) = \mathcal{A}^\bullet(W) \oplus H^n(W,\mathbb{C})_{\mathrm{prim}} \, ,
		\end{equation}
		where $\mathcal{A}^{\bullet}(W) = \mathrm{span}_{\mathbb{C}}\left\{\mathrm{P}^0=1,\mathrm{P},\dots,\mathrm{P}^n\right\} \subseteq H^\bullet(W,\mathbb{C})$ is generated by the algebraic classes, $H^n(W,\mathbb{C})_{\mathrm{prim}}$ is the image of the primitive cohomology group in small quantum cohomology, and $\oplus$ is orthogonal with respect to the Poincar\'e pairing. The 3-point invariants in (\ref{eq-sec_derivative_I-function}) contain no contributions from primitive classes.
	\end{lemma}
	\begin{proof}
		Equation (\ref{eq-sec_derivative_I-function}) is standard and follows from the divisor equation (in the flat coordinate) applied twice, and depends only on the existence of $I_W(\tau,z)$ and an expansion of the form (\ref{eq-I-function_expansion}) with no restriction on the cohomology of $W$. (The equation (\ref{eq-sec_derivative_I-function}) should be thought of formally if there are no further assumptions on convergence in small $QH_{\mathrm{sm}}(W)$ dependent on $\tau$.) In fact, (\ref{eq-sec_derivative_I-function}) is simply the coefficient of $\mathrm{P}^2$ in the quantum product $\mathrm{P} \star_\tau \mathrm{P}$. 
		\par That $QH_{\mathrm{sm}}(W)$ has the decomposition (\ref{eq-QH_decomp}) follows from the assumption that  $H^\bullet(W,\mathbb{C})$ is almost minimal. It follows that $\mathcal{A}^{\bullet}(W)$ is stable under $\mathrm{P} \star_\tau$. For $\eta \in H^n(W,\mathbb{C})_{\mathrm{prim}}$, it can be shown that $\mathrm{P} \star_\tau \eta$ has no component in $\mathcal{A}^{\bullet}(W)$, so the structure constants of the algebraic subring $\mathcal{A}^{\bullet}(W) \subseteq QH_{\mathrm{sm}}(W)$ do not involve contributions from primitive classes. Thus the claim about the 3-point invariants in (\ref{eq-sec_derivative_I-function}) follows.  
	\end{proof}
	\begin{definition}
		\label{def-3-point_function}
		For an $n$-dimensional almost minimal variety $W$ as in Lemma \ref{lem-3-point_I-function}, the 3-point function in (\ref{eq-sec_derivative_I-function}) will be denoted by the symbol $\mathcal{Y}^{(3)}_{\mathrm{P}^{n-2}}(W)$. As a convention, we define $\mathcal{Y}^{(3)}_{\mathrm{P}^{-1}}(W) = 1$ in the case $\dim(W) = 1$.
	\end{definition}
	The actual geometric data that is extracted from (\ref{eq-sec_derivative_I-function}) depends heavily on whether or not $K_W = 0$ or $-K_W > 0$, as well as $\dim(W)$. Let us write $X^\vee$ to be an almost minimal CY variety ($K_{X^\vee} = 0$) and $V$ to be an almost minimal Fano variety ($-K_V > 0$). For notational purposes, let us write particular structure constants relevant to the 3-point invariants (\ref{eq-sec_derivative_I-function}) of $\mathrm{QH}^{\mathrm{sm}}(X^\vee)$ as
	\begin{equation*}
		C_k(q) := \langle \langle \mathrm{P}, \mathrm{P}^k, \mathrm{P}^{n-2} \rangle \rangle \, ,
	\end{equation*} 
	so that $C_1(q)$ is what appears in (\ref{eq-sec_derivative_I-function}).
	\begin{theorem}
		\label{thm-CY_I-function}
		Suppose that $X^\vee$ is an almost minimal CY variety of dimension $n \geq 2$ that satisfies the conditions of Lemma \ref{lem-3-point_I-function}.  
		\begin{enumerate}
			\item[a)] If $\dim(X^\vee) = 2$, then $\mathcal{Y}^{(3)}_{\mathrm{P}^{0}}(X^\vee) = 1$. 
			\item[b)] If $\dim(X^\vee) = n = 3, 4$, then $\mathcal{Y}^{(3)}_{\mathrm{P}^{n-2}}(X^\vee)$ generates $\mathcal{A}^\bullet(X^\vee) \subseteq QH_{\mathrm{sm}}(X^\vee)$. 
			\item[c)] If $\dim(X^\vee) \geq 5$, then $\mathcal{Y}^{(3)}_{\mathrm{P}^{n-2}}(X^\vee)$ determines only the 3-point invariants in (\ref{eq-sec_derivative_I-function}).
		\end{enumerate}
	\end{theorem}
	\begin{proof}
		Notice that since $K_{X^\vee} = 0$, $q$ has degree zero. Thus $\star_\tau$ is degree preserving and quantum multiplication collects all structure constants into the same cohomological slot.
		\begin{enumerate}
			\item[a)] It is well-known that all genus zero invariants vanish for $X^\vee$ a K3 surface; this is reviewed in \S \ref{ss-genus_zero_invariants}. Thus we have
			\begin{equation}
				\label{eq-K3_cup_product}
				\mathrm{P} \star_\tau \mathrm{P} = \mathrm{P} \cup \mathrm{P} = \kappa[\mathrm{pt}] \, ,
			\end{equation}
			whence the result follows since the prepotential is at most quadratic in $\tau$. 
			\item[b)] If $n = 3$, then $\mathrm{P}^{n-2} = \mathrm{P}$, and thus every $\langle \mathrm{P}, \mathrm{P}, \mathrm{P} \rangle_d$ is determined by the quantum product $\mathrm{P} \star_\tau \mathrm{P}$ since $\mathcal{A}^\bullet(X^\vee) \cong \mathbb{C}[[q]][\mathrm{P}]/\langle \mathrm{P}^{\star_\tau 4}\rangle$. Similarly if $n = 4$, $\mathrm{P}^{n-2} = \mathrm{P}^2$, and we have by symmetry that $\langle \mathrm{P}, \mathrm{P}^2, \mathrm{P} \rangle = \langle \mathrm{P},\mathrm{P}, \mathrm{P}^2\rangle$. Thus $\mathrm{P} \star_\tau \mathrm{P}^2 = C_1(q) \mathrm{P}^3$. By Lemma \ref{lem-3-point_I-function}, the remainder of the claim follows. 
			\item[c)] In the case $n \geq 5$, the quantum product contains genuine 2-point invariants that cannot be recovered from 1-point invariants from the divisor equation, so the higher order $C_k(q)$ with $k \geq 2$ are genuinely independent from $C_1(q)$. 
		\end{enumerate}
	\end{proof}
	\begin{theorem}
		\label{thm-Fano_I-function}
		Suppose that $V$ is an almost minimal Fano variety of dimension $n \geq 2$. Then $\mathcal{Y}^{(3)}_{\mathrm{P}^{n-2}}(V) = 1$.
	\end{theorem}
	\begin{proof}
		Since $V$ is Fano, $q$ has positive degree, and thus quantum product in $QH_{\mathrm{sm}}(V)$ agrees with the classical product except in top quantum degree $(=n+1)$. The second derivative (\ref{eq-sec_derivative_I-function}) thus only collects classical products; by almost minimality we have $\mathrm{P} \star_\tau \mathrm{P}$ is given by (\ref{eq-K3_cup_product}). 
	\end{proof}
	\par Summarizing, we have the following situation.
	\begin{corollary}
		\label{cor-3-point_agreement}
		For $V$ an almost minimal Fano variety of dimension $n \geq 3$ the second order derivatives as in (\ref{eq-sec_derivative_I-function}) obtained from $I^\mathrm{eu}_0$ with respect to the mirror maps $\tau_{\mathrm{eu}}$, respectively, determine the 3-point invariants encoded in $\mathcal{Y}_{\mathrm{P}^{n-2}}(X^\vee)$. In the case $\dim(V) \geq 2$, the analogous situation with $I^\mathrm{loc}_0$ and $\tau_{\mathrm{loc}}$ yields the 3-point invariants encoded in $\mathcal{Y}_{\mathrm{P}^{n-2}}(X^\vee)$.
	\end{corollary}
	%
	%\begin{remark}
	%	\label{rmk-local_I-function}
	%	The results of Lemma \ref{lem-3-point_I-function} also extend to the case of quasiprojective varieties; in particular, Theorem \ref{thm-CY_I-function} extends to the analogous situation with local CY varieties. Moreover, the $I$-functions relevant to the local and ambient $\mathrm{SQDM}$ can always be defined in the context of Iritani, Mann, \& Mignon. They are described in terms of the 2-point invariants $N_{\alpha,d}(z)$ in (\ref{eq-N_invariants}) and modified hypergeometric expressions \cite[Def. 5.11, Lem. 5.4]{MR3519131} for $\tau \in H^2(V,\mathbb{C})$
	%	\begin{equation}
	%		\label{eq-local_I-function}
	%		I_\alpha^{\mathrm{loc}}\left(\tau, %z\right):=e^{\tau / z} \sum_{d \in \operatorname{Eff}(V)} N_{\alpha, d}(z) e^{\tau(d)} \prod_{k=0}^{c_1(V)(d)-1}(-c_1(V)-k z) \, .
	%	\end{equation}
	%	It will follow from our Assumption \ref{ass-mirror} that such data can always be obtained through mirror symmetry.
	%\end{remark}
	%
	\subsection{Genus zero invariants of elliptic curves and K3 surfaces}
	\label{ss-genus_zero_invariants}
	\par In this section let $X^\vee$ be either an elliptic curve or a K3 surface. It is well-known that the genus zero Gromov-Witten theory of elliptic curves and K3 surfaces are essentially vacuous, since in the former case the moduli space of stable maps with source $\mathbb{P}^1$ and target an elliptic curve $X^\vee$ is empty as there are no nonconstant maps $\mathbb{P}^1 \to X^\vee$; in the latter case the moduli stack of stable maps fails to have a virtual fundamental class of positive dimension, for any effective curve class $\beta \in H_2(X^\vee,\mathbb{Z})$. This is a manifestation of the fact that a completely generic complex structure on a K3 surface $X^\vee$ is not algebraic, thus $X^\vee$ contains no curves. This implies in both cases that the quantum cohomology ring agrees with the classical cohomology ring \cite[Ex. 8.4]{MR1266766}. More generally, Kiem \& Li \cite{MR3073883} established important vanishing results for (quasi)projective varieties that carry holomorphic symplectic forms by cosection localization. 
	\par However, the higher genus theory for elliptic curves is nontrivial, and for example the Gromov-Witten potential in each genus $g \geq 2$ can be expressed as a quasi-modular form of weight $6g-6$ in terms of the classical Eisenstein series $G_2,G_4,G_6$ \cite{MR1363055}. For K3 surfaces, the so-called reduced descendant Gromov-Witten theory allows for explicit expression of the analogue of the Gromov-Witten potential for all $g \geq 0$ as a quasi-modular form in terms of $G_2$ and the reciprocal of the classical modular discriminant $\Delta$. This was first conjectured in the physics literature by Yau \& Zaslow \cite{MR1398633}, then subsequently proved in separate works by Beauville \cite{MR1682284}, Fantechi, G\"ottsche, \& van Straten \cite{MR1658220}, and Bryan \& Leung \cite{MR1750955}. A wide generalization of these results was conjectured by G\"ottsche for curves on any surface \cite{MR1645204} and subsequently proven by Kool, Shende, and Thomas \cite{MR2776848}. There is continued interest in such structures, such as the recent work of Oberdieck \cite{MR4953124}. In the context of this article, we do not use any of the higher genus or refined structure referenced here, in fact it will be of significance for us that the classical genus zero invariants in either case are trivial. This will be manifest through our expression of the nonequivariant limit of quantum Serre duality in terms of 3-point euler-twisted invariants and relating them to the ordinary 3-point invariants of anticanonical CY varieties in \S \ref{s-projective}, and is intimately connected to the growth rate of local 3-point invariants.
	\section{Mirror symmetry in almost minimal Fano varieties}
	\label{s-mirror}
	%
	%To do (Michael): Description of local mirror symmetry for local threefolds \& fourfolds, Picard-Fuchs equations, the twist construction, relevant results
	%
	\par In this section we recall pertinent aspects related to mirror symmetry of anticanonical Calabi-Yau varieties in a smooth, almost minimal Fano $n$-fold $V$. We will set $H = -K_V$ to be the anticanonical class and degree will be measured with respect to $H$.  We will be concerned only with the implications for the small quantum $D$-module $\mathrm{SQDM}(V)$. 
	
	%Quantum Serre duality as in \S \ref{ss-QSD} in particular provides a relation of the euler-twisted quantum $D$-module as embedding into the pullback of the local quantum $D$-module by the anticanonical shift of the 2-point Euler twisted invariants; the assumption that $V$ is of almost minimal guarantees that the ambient part of the small quantum $D$-module is the entire thing.
	%
	\subsection{Fano mirrors}
	\label{ss-fano_ms}
	\par The mirror of a Fano variety $V$ is a Landau-Ginzburg (LG) model. It is (locally) given as the data of a Laurent polynomial $\mathcal{W}$ in $n$-variables known as a superpotential, and is thought of as a holomorphic map $\mathcal{W} : (\mathbb{C}^\times)^n \to \mathbb{C}$. As in Assumption \ref{ass-nonempty}, we let $X^\vee \in |-K_V|$ be a very general smooth anticanonical CY variety. The quantum Lefschetz theorem of Coates \& Givental \cite{MR2276766} then tells us that the GW invariants of $X^\vee$ are determined by the invariants of the ambient $V$ twisted by $-K_V$ in the sense described in \S \ref{ss-QSD}. Let $\tau \in H^{2}(X^\vee,\mathbb{C})/H^2(X^\vee,\mathbb{Z})$ record the corresponding choice of (complexified) symplectic structure on $X^\vee$. It represents the complexification of a positive class in $H^{1,1}(X^\vee)$, and thus $\tau$ lives in a space that is naturally identified with $\mathbb{H}/\mathbb{Z}$, where $\mathbb{H} \subset \mathbb{C}$ is the upper half plane. Mirror symmetry predicts that the minimal, smooth compactification of the generic fibre $X_t := \overline{\mathcal{W}^{-1}(t)}$ is mirror to $X^\vee_\tau$ in a suitable way that is recalled in \S \ref{ss-closed_ms}. This is known as the LG-CY correspondence.
	\par There are precisely 17 deformation classes of smooth rank-1 Fano threefolds as classified by Iskovskikh \cite{MR0463151,MR0503430}, which were studied in the context of mirror symmetry by Golyshev \cite{MR2306141}. It is easy to see that every rank-1 Fano threefold is either minimal or almost minimal. There, \cite[\S 2.3]{MR2306141} he introduced the concept of a ``counting differential operator of type D3", a fourth order linear operator that is obtained from the $\mathrm{SQDM}(V)$, where the `D' stands for `determinantal'. It depends on a matrix $A_V = [a_{ij}]$ that records the (suitably normalized) 2-point invariants of $V$ as follows.
	\begin{definition}
		\label{def-A_V}
		[Golyshev \cite[\S 2.3]{MR2306141}, Golyshev \& Stienstra \cite[\S 1.1]{MR2346574}]
		Let $V$ be a smooth, almost minimal Fano $n$-fold and $H = -K_V$. Define a matrix $A_V = [a_{ij}]$ of size $(n+1) \times (n+1)$ with $0 \leq i,j \leq n$ by 
		\begin{equation}
			\label{eq-normalized_2-point}
			a_{ij} = \frac{j-i+1}{\mathrm{deg}(V)} \cdot \langle H^j,H^{n-i}\rangle_{0,2,j-i+1}  \quad \mathrm{if} \; \; j - i > 1 \, ,
		\end{equation}
		which records up to normalization the number of maps $\mathbb{P}^1 \to V$ of degree $j-i+1$ such that $0$ lies on a generic representative of $H^j$ and $\infty$ lies on a generic representative of $H^{n-i}$. Here $\mathrm{deg}$ is measured with respect to $H$. Furthermore, define $a_{ij} = 0$ if $i-j > 1$, $a_{ij} = 1$ if $i-j=1$, and $a_{ij} = a_{n-j,n-i}$.
	\end{definition}
	It follows from basic properties of quantum multiplication with $H$ that $a_{ij} = a_{n-j,n-i}$. The significance of the matrix $A_V$ is given by the following result, also due to Golyshev \cite[Prop. 2.2, 2.4]{MR2306141}. Assume for simplicity that $V$ has index $\mathrm{ind}(V)=1$. Let $\mathcal{A}^{\bullet}(V) = \mathrm{span}_{\mathbb{C}}\{\mathrm{P}^0=1,\mathrm{P},\dots,\mathrm{P}^n\} \subseteq H^\bullet(V,\mathbb{C})$ be the subspace generated by the algebraic classes, and $\mathbf{P} = (\mathrm{P}^0,\mathrm{P},\dots,\mathrm{P}^{n})$. Let $t$ be the co-ordinate dual to $\mathrm{P} \in H^2(V,\mathbb{C})$, and let $\theta = t \frac{d}{dt}$ be the usual logarithmic derivative operator. Define an $(n + 1) \times (n+1)$ matrix $\omega = [\omega_{ij}]$ from the matrix $A_V$ by 
	\begin{equation}
		\label{eq-A_side_connection}
		\omega_{ij} = 
		\begin{cases}
		a_{ij}t^{j-i+1}	\, , \quad j - i > 1 \, ,\\
		a_{ij} \, , \quad j - i \leq 1 \, .
		\end{cases}
	\end{equation}
	\begin{proposition}
		\label{prop-A_V_SQDM}
		Consider $V$, $A_V$ as in Definition \ref{def-A_V}. Then the small quantum $D$-module $\mathrm{SQDM}(V)$ of $V$ is determined entirely by the matrix $A_V$ by the relation 
		\begin{equation}
			\label{eq-A_side_SQDM}
			\theta \mathbf{P} = \mathbf{P}\omega \, .
		\end{equation}
	\end{proposition}
	\begin{proof}
		The `dimension axiom' of GW invariants implies that for $V$ almost minimal, the quantum connection $\nabla$ on $H^\bullet(V)$ restricts to the rank-$(n+1)$ sub-bundle generated by $\mathcal{A}^{\bullet}(V)$, since the image of $\mathcal{A}^{\bullet}(V)$ in $QH(V)$ is stable under quantum multiplication by $\mathrm{P}$. This precisely determines $\mathrm{SQDM}(V)$, which by definition is $\mathrm{QDM}(V)$ restricted to the divisorial directions in $H^\bullet(V,\mathbb{C})$. The 2-point invariants given in (\ref{eq-normalized_2-point}) (up to normalization) determine $\mathrm{SQDM}(V)$ by (\ref{eq-A_side_SQDM}) from the definition of $\nabla$ restricted to $\mathcal{A}^{\bullet}(V)$ and the divisor equation.
	\end{proof}
	\par The quantum connection given as Equation (\ref{eq-A_side_SQDM})  is known as the \emph{first} structure connection. In fact, $\mathrm{P}^0=1$ determines a cyclic vector for the connection; rewriting (\ref{eq-A_side_SQDM}) as a linear, homogeneous ODE operator $\check{\mathcal{D}}_{n+1}$ of order $n+1$, it follows that the operator $\check{\mathcal{D}}$ annihilates $\mathrm{P}^0$, and moreover, it can be shown that $\check{\mathcal{D}}_{n+1}$ has precisely two singular points at $t=0,\infty$, with the singular point at $t=0$ regular and the one at $t=\infty$ irregular as described in \S \ref{ss-SQDM}. In fact, the roots of the indicial equation at $t=0$ are all zero, so $t=0$ is a so-called point of maximally unipotent monodromy (MUM). It follows that at $t=0$, there is a canonical single valued holomorphic solution $\Psi_0(t)$ of the homogeneous equation  $\check{\mathcal{D}}_{n+1} =0$ that can be obtained by the Frobenius method. As written, $\mathrm{SQDM}(V)$, or equivalently $\check{\mathcal{D}}_{n+1}$, determines a $D$-module over $\mathbb{A}^1 = \mathrm{Spec} \, \mathbb{C}[t]$, whence writing 
	\begin{equation}
		\label{eq-SQDM_operator}
		\check{\mathcal{D}}_{n+1} = \sum_{j=0}^{n+1} a_j(t) \partial^j \, , \quad \partial = \frac{d}{dt} \, ,
	\end{equation}
	the coefficient functions $a_j(t) \in \mathbb{C}[t]$ are polynomials of suitable order and zeroes to have the required singularity properties.
	\par Because of the presence of an irregular singular point, the operator $\check{\mathcal{D}}_{n+1}$ cannot be of geometric origin, in the sense of being a Picard-Fuchs operator of some pencil of (CY) varieties. Picard-Fuchs operators are well known to Fuchsian \cite{deligne_equations_1970,MR229641,MR233825,MR282990} as they underlie the associated Gauss-Manin connection, and hence possessing only regular singular points. One interpretation of the mirror symmetry conjecture \cite[1.9, 1.10]{MR2306141} is that a suitable regularization $\nabla^{\mathrm{reg}}$ of $\mathrm{SQDM}(V)$ is of geometric origin. 
	In the case of $V$ a rank-1 Fano threefold, Golyshev showed \cite[1.8]{MR2306141} that such a regularization of $\mathrm{SQDM}(V)$ can be obtained by convolution of the small quantum connection $\nabla$ with the so-called canonical exponent object on the torus $T_{\mathrm{NS}(V)^{\vee}}$ obtained from the dual N\'eron-Severi lattice of $V$. Then, the regularized object $\mathrm{SQDM}(V)^{\mathrm{reg}}$ turns out to be a single, reducible $D$-module that is obtained from the Fourier-Laplace transform of the irregular operator $\check{\mathcal{D}}_4$; the desired Picard-Fuchs operator could be then realized as an irreducible subquotient thereof. The Fourier-Laplace transform is recalled in Remark \ref{rmk-fourier_laplace}.
	\par More precisely, the exponent object on $\mathbb{G}_m = \mathrm{Spec} \, \mathbb{C}[t,t^{-1}]$ is defined to be the $D$-module $E = \mathcal{D}_{\mathbb{G}_m}/\langle \theta - t \rangle$ on $\mathbb{G}_m$ determined by the $\mathbb{G}_m$-invariant logarithmic derivative operator $\theta$.\footnote{The exponent object $E$ is obtained from pullback of $\mathcal{D}_{\mathbb{A}^1}/\langle \partial -1 \rangle$ under the canonical open immersion $\mathbb{G}_m \hookrightarrow \mathbb{A}^1$.} There is a canonical morphism of tori $\iota : \mathbb{G}_m \to T_{\mathrm{NS}(V)^\vee}$ that comes from twice dualizing the inclusion $\mathrm{span}_{\mathbb{Z}} \{K_V\} \hookrightarrow \mathrm{NS}(V)$, and one defines the canonical exponent object on $T_{\mathrm{NS}(V)^\vee}$ as $E^\vee := \iota_*E$. Then the regularization $\mathrm{SQDM}(V)^{\mathrm{reg}}$ is defined in terms of convolution by
	\begin{equation}
		\label{eq-regularized_SQDM}
		\mathrm{SQDM}(V)^{\mathrm{reg}} := \mathrm{SQDM}(V) * E^\vee \, ,
	\end{equation}
	where $*$ is convolution of $D$-modules on $T_{\mathrm{NS}(V)^\vee}$ \cite[1.2]{MR2306141}. Crucially, this operation can be realized up to composition with the inversion $t \mapsto t^{-1}$ in terms of the Fourier-Laplace transform of $\mathrm{SQDM}(V)$  \cite[\S 5.2.3, 5.2.3.1]{MR1081536}, \cite[2.8]{MR2306141}. It was shown \cite[\S 3, 4]{MR2306141} that for each pair $(N,d)$ of level and Fano index in (\ref{eq-degree_level_invariants}), $\mathrm{SQDM}(V)^{\mathrm{reg}}$ is very closely related to the Gauss-Manin connection underlying a polarized variation of Hodge structure (VHS) of K3 surfaces that are polarized by the rank-19 lattice 
	\begin{equation}
		\label{eq-M_n-lattice}
		\mathrm{M}_N := \mathrm{H} \oplus \mathrm{E}_8(-1)^{\oplus 2}  \oplus \langle -2N \rangle \, .
	\end{equation}
	The Gauss-Manin connection is isomorphic to the desired irreducible subquotient of $\mathrm{SQDM}(V)^{\mathrm{reg}}$.
	\par Such lattice polarized K3 surfaces were first studied by Dolgachev \cite[\S 7]{MR1420220}, where they were shown to be mirror in a suitable sense to rank-1 K3 surfaces polarized by the lattice $\langle 2N \rangle$. More details will be recalled in the sequel \S \ref{ss-closed_ms}, and in (\ref{eq-K3_mirror}). Roughly speaking, these K3 surfaces are members of the anticanonical linear system of a Fano threefold of level $N$. The coarse moduli space of $\mathrm{M}_N$-polarized K3 surfaces is known to be isomorphic to the modular curve $X_0(N)^+$ \cite[Thm. 7.1]{MR1420220}.\footnote{The modular curve is defined as $X_0(N)^+ =\overline{\mathbb{H}/\Gamma_0(N)^+}$ , where $\Gamma_0(N)^+ \subset \mathrm{SL}(2,\mathbb{R})$ is obtained from the principal congruence subgroup $\Gamma_0(N)$ by adjoining the Fricke involution $\tau \mapsto -1/N\tau$.} For precisely the values of $N$ in (\ref{eq-degree_level_invariants}), $X_0(N)^+ \cong \mathbb{P}^1$. The associated VHS over $X_0(N)^+$ was studied at the level of Picard-Fuchs operators by Doran \cite[\S 5.2]{MR1738862}, and it was subsequently shown that up to pullback determined by the branched covering map $t \mapsto t^\delta$, the associated third order Picard-Fuchs operator $\mathcal{L}_3$ of each family with numerical invariants $(N,\delta)$ \cite[\S 4]{MR2306141}, which in this case carries the same data as the Gauss-Manin connection of the family. More details will be discussed in shortly and in \S \ref{ss-closed_ms}. This essentially proves the mirror symmetry conjecture for rank-1 Fano threefolds in the sense described in the preceding paragraph. 
	\begin{remark}
	\label{rmk-fourier_laplace}
		At the level of the operator $\check{\mathcal{D}}_{n+1}$ underlying $\mathrm{SQDM}(V)$, this can be computed concretely as follows. Writing $\check{\mathcal{D}}_{n+1}$ as in (\ref{eq-SQDM_operator}), we have the Fourier-Laplace transform $\mathcal{D}_{n+1} := \mathrm{FL}(\check{\mathcal{D}}_{n+1})$, where
	\begin{equation}
		\label{eq-fourier_laplace}
		\mathrm{FL}(\check{\mathcal{D}}_{n+1}) = \sum_{j=0}^{n+1} a_j(\partial)(-t)^j \, .
	\end{equation}
	Suppose that $\Psi_0(t) = \sum_{n=0}^\infty c_n \frac{t^n}{n!}$ is the fundamental holomorphic solution of $\check{\mathcal{D}}_{n+1}=0$ at the MUM point $t=0$. It can be shown that $t=0$ is also a MUM point, and that 
	\begin{equation}
		\label{eq-holomorphic_fundamental_solution}
		\omega_0(t) :=  \sum_{n=0}^{\infty} c_n t^n
	\end{equation}
	is the fundamental holomorphic solution of $\mathcal{D}_{n+1} = 0$ at that point. It follows that $\Psi_0(t) = \exp(t) \star \omega_0(t)$ is the exponential generating function of the regularized object $\omega_0(t)$. Here $\star$ is the Hadamard product\footnote{The Hadamard product of two series $f(t) = \sum_{n=0}^{\infty}a_nt^n$ $g(t) = \sum_{n=0}^{\infty}b_nt^n$ is defined as $(f \star g)(t) = \sum_{n=0}^{\infty}a_nb_nt^n$.} of series, in line with the definition of regularization given in (\ref{eq-regularized_SQDM}). Moreover, the irregular singular point at $t=\infty$ of $\check{\mathcal{D}}_{n+1}$ becomes a regular singular point for $\mathcal{D}_{n+1}$.  It is noted by Golyshev \& Zagier \cite[pp. 26]{MR3462676} that a more appropriate name might be `Borel' transform. In this vein, for explicit examples such as those discussed in \S \ref{s-projective}, this can be computed, for example using the commands \textbf{borel}, \textbf{invborel} in the `gfun' package of \texttt{Maple}.
	\end{remark}
	\par Golyshev \& Stienstra \cite{MR2346574} then studied general properties of the operator $\check{\mathcal{D}}_{n+1}$ obtained from matrices as in (\ref{eq-A_side_connection}, \ref{eq-A_side_SQDM}), and showed that suitable subquotients of their Fourier-Laplace transform (\ref{eq-fourier_laplace}) are candidates for studying mirror symmetry of almost minimal Fano varieties. They proved a key property of the regularized operator $\mathcal{D}_{n+1}$ that is crucial for us in the context of this article.
	\begin{theorem}
		\label{thm-operator_factorization}
		[Golyshev \cite[2.10]{MR2306141}, Golyshev \& Stienstra, \cite[Thm. 5.8]{MR2346574}] The regularized operator $\mathcal{D}_{n+1}$ factorizes as 
		\begin{equation}
			\label{eq-quantum_operator}
			\mathcal{D}_{n+1} = \theta \circ (t\mathcal{L}_n) \, ,
		\end{equation}
		where $\mathcal{L}_n$ is an irreducible Fuchsian operator of order $n$ that is maximally unipotent at $t=0$. Moreover, the operator $\mathcal{L}_n$ is (anti) self-adjoint, $\mathcal{L}_n^\vee = (-1)^n\mathcal{L}_n$.
	\end{theorem}
	\par Thus the operator $\mathcal{L}_n$, the irreducible subquotient of the regularized object $\mathcal{D}_{n+1}$, is a candidate operator to be of geometric origin. It is deemed a determinantal $\mathrm{D}n$ operator, and explicit formulae for example for the determinantal $\mathrm{D}3$ operators $\mathcal{L}_3$ may be found in \cite[2.11]{MR2306141}. Although one may naively believe that all of the relevant geometric content is contained in the irreducible operator $\mathcal{L}_n$, the purpose of this article is to show that the full regularized operator $\mathcal{D}_{n+1}$ naturally contains important geometric data about the associated Fano variety, namely, the genus zero Euler-twisted data discussed in \S \ref{ss-twisted_invariants}. This is discussed further in Theorem \ref{thm-3-point_euler-twisted}b) and Remark \ref{rmk-gamma_conjecture}. 
	\begin{remark}
		\label{rmk-second_structure}
		The regularized operator $\mathcal{D}_{n+1}$ is essentially Dubrovin's \emph{second} structure connection \cite[Def. 3.1]{MR1397274}, \cite[Ch. II, \S 1.1]{MR1702284}.
	\end{remark}
	We now turn towards mirror symmetry for the anticanonical CY varieties in $V$.
	\subsection{Closed mirror symmetry}
	\label{ss-closed_ms}
	\par For $V$ an almost minimal Fano, from Assumption \ref{ass-nonempty} we let $X^\vee \in |-K_V|$ be an anticanonical CY hypersurface in $V$. Then $X^\vee$ has a 1-dimensional complexified symplectic moduli space this is locally isomorphic to $\mathbb{H} / \mathbb{Z}$, where $\mathbb{H}$ is the upper half plane \cite{yang2019ap}. It is essentially the complexification of the positive cone in $H^{1,1}(X^\vee)$. Denote $X^\vee = X^\vee_\tau$ to be the variety $X^\vee$ equipped with such a complexified symplectic class from $\tau \in \mathbb{H}/\mathbb{Z}$. The limit $\tau \to i\infty$ is known as the large volume limit of $X^\vee$; as such it is convenient to analytically embed the moduli space into a copy of the projective line, which we denote by $\mathbb{P}^1_\tau$. We lose nothing by taking $\mathbb{P}^1_\tau$ to be the compactification of the full moduli space. Then we get a family $\pi^\vee : \mathcal{X}^\vee \to \mathbb{P}^1_\tau$, whose fibres are $(\pi^\vee)^{-1}(\tau) = X^\vee_{\tau}$ parameterized by variation in $\tau$. The singular locus of the fibration includes the large volume limit $\tau = i \infty$.
	\par If $V$ is a toric variety, the well-known mirror construction of Batyrev \& Borisov \cite{MR1416334} constructs a mirror variety $X=X_t$ that satisfies the mirrored relation of Hodge numbers $h^{p,q}(X^\vee_{\tau})=h^{n-1-p,q}(X_t)$. From the relatively minimal assumption, it follows that the middle Hodge numbers satisfy $h^{i,j}(X_t) = 1$ for all $i+j=n-1$, and thus has middle Betti number $b_{n-1}(X_t) = n$. Crucially, the mirror construction includes the data an analytic isomorphism $\tau$ that carries a neighborhood of the large volume limit in $\mathbb{P}^1_\tau$ to a suitable open subset of the complex structure moduli space of the mirror variety; in analogy to above we take $\mathbb{P}^1_t$ to be the compactification of this complex structure moduli space. The mirror map is further known to agree with the mirror maps described in \S \ref{ss-I_functions}. In fact, given the discussion in \S \ref{ss-fano_ms}, it is reasonable to take $\mathbb{P}^1_t$ as the target of the superpotential $\mathcal{W}$ associated to $V$. 
	\par The mirror map allows for the quantum corrected symplectic geometry of $X^\vee_{\tau}$ to be expressed in terms of the classical complex geometry of $X_t$, at least in a neighborhood of the large volume limit. More precisely, the map may be taken to be of the form $\tau : U \to \mathbb{P}^1_\tau$ defined in an analytic neighborhood $U \subseteq \mathbb{P}^1_t - \{0,\infty\}$ of the limiting value $t \to 0$ that maps $0 \mapsto i\infty$. The limiting value $t \to 0$ is called the large complex structure limit (LCSL), and corresponds to the variety $X_t$ acquiring maximal degeneration. Givental's toric mirror theorem \cite{MR1653024} then identifies the components of the hypergeometric $I$-function of $X^\vee$ as in Lemma \ref{lem-3-point_I-function} with the ratios of the classical period integrals of the mirror variety $X$ up to multiplication by an oscillatory integral. The period integrals are given as functions of the complex structure moduli by
	\begin{equation}
		\label{eq-periods}
		\omega_i(t) := \oint_{\Sigma_i} \Omega_t \, ,
	\end{equation}
	where $H^{n-1,0}(X_t) = \mathrm{span}_{\mathbb{C}}\{\Omega_t\}$ and $H_{n-1}(X_t,\mathbb{Z}) = \mathrm{span}_{\mathbb{Z}}\{\Sigma_0,\dots,\Sigma_{n-1}\}$\footnote{For simplicity assume that $H_{n-1}(X_t,\mathbb{Z})$ is torsion free. Some care must be taken here in the case that $n-1 \in 2\mathbb{Z}$ is even. For example, if $n-1=2$, we integrate over an integral basis of the transcendental lattice $\mathrm{T}(X_t) := \mathrm{Pic}(X_t)^\perp$ of rank $22-\rho$. Here $\perp$ is taken with respect to the K3 lattice (\ref{eq-K3_mirror}).}. The periods are well known to be generically multivalued functions that are annihilated by the Picard-Fuchs operator $\mathcal{L}_n$ associated to $X$, an $n$th order Fuchsian operator defined on $\mathbb{P}^1_t - \{\mathrm{singular \; points}\}$ that underlies the full Gauss-Manin connection for the family $X_t$. In the context of \S \ref{ss-fano_ms}, $\mathcal{L}_n$ is precisely the irreducible subquotient of the regularized operator $\mathcal{D}_{n+1}$. 
	\par The fundamental idea is that, when the mirror maps $\tau_{\mathrm{eu}},\tau_{\mathrm{loc}}$ in (\ref{eq-I-mirror_maps}) have geometric origin as a corresponding ratio of period integrals, whence the oscillatory factors cancel appropriately. The geometric form of these mirror maps is given in Definitions \ref{def-closed_mirror_map} and \ref{def-local_mirror_map}.  When there is a known mirror construction, such as when $V$ is a toric Fano variety as above, Iritani \cite{MR2553377,MR3112512} further clarified that  there is an integral structure in quantum cohomology that recovers the underlying integral variation of Hodge structure (VHS) for the mirror CY on the B-side. In other cases, such as when $V$ is a smooth Fano threefold or del Pezzo surface, the mirror construction for anticanonical K3 surfaces or elliptic curves is known to agree with predictions from the LG-model perspective \cite{MR4913468}. 
	\par By utilizing the Quantum Lefschetz theorem of Coates \& Givental \cite{MR2276766}, this perspective has found in particular powerful applications in the classification of 3-dimensional smooth Fano varieties. In \cite{MR3470714}, it was shown that the quantum period of $V$, a certain generating function of its 1-point genus zero GW invariants, provided a strong enough invariant so as to distinguish between deformation classes of Fano threefolds as classified by Mori \& Mukai \cite{MR641971,MR1102273,MR1969009}. LG models for all deformation classes were provided in \cite{MR4913468}, where it was in addition shown that the expected Dolgachev-Pinkham-Nikulin mirror symmetric relation for a smooth, very general K3 fibre held. That is,
	\begin{equation}
		\label{eq-K3_mirror}
		\mathrm{Pic}(X^\vee_\tau)^\perp \cong \mathrm{Pic}(X_t) \oplus \mathrm{H} \, ,
	\end{equation}
	where $\mathrm{H}$ is the canonical unimodular hyperbolic lattice of rank-2 and $\perp$ is taken with respect to the ambient K3 lattice $\Lambda_{\mathrm{K3}} = \mathrm{H}^{\oplus 3} \oplus \mathrm{E}_8(-1)^{\oplus 2}$. Moreover, the authors provided a LCSL for the K3 pencil. Thus we make the following assumption. 
	\begin{assumption}
		\label{ass-mirror}
		The anticanonical CY variety $X_\tau^\vee \in |-K_V|$ has a smooth mirror variety $X_t$ that satisfies the properties above. If $\pi : \mathcal{W} \to \mathbb{C} \subset \mathbb{P}^1_t$ is a superpotential for $V$, then $X_t$ can be obtained from the minimal compactification of the fibres of $\mathcal{W}$, $X_t = \overline{\pi^{-1}(t)}$, and the collection $\mathcal{X} \to \mathbb{P}^1$ of all such mirror varieties forms a holomorphic fibration of polarized CY varieties on the complement of discrete singular locus $\{0,\dots,\infty\}$ on $\mathbb{P}^1$. We take $t \to 0$ to be the large complex structure limit of the mirror family $\mathcal{X}$. There is a correspondence of integral structures on the A-side and B-side so that the mirror maps (\ref{eq-I-mirror_maps}) are determined precisely by appropriate ratios of the period integrals (\ref{eq-periods}) of the mirror variety.
	\end{assumption}
	\begin{remark}
		\label{rmk-empty}
		In the case that $|-K_V| = \emptyset$ but there is still a mirror construction for $V$, one may in general lose the LG-CY correspondence. For example, the $\mathbb{Q}$-del Pezzo surface described in Remark \ref{rmk-Q_Fano}, Corti \& Gugiatti \cite{MR4337923} showed that the generic fibre of the associated LG model is a smooth curve of genus 2.
	\end{remark}
	\par Expressed in terms of the logarithmic derivative $\theta := t \frac{d}{dt}$, near $t=0$ the Picard-Fuchs operator $\mathcal{L}_n$ of the mirror family $X_t$ is expressed as
	\begin{equation}
		\label{eq-pf_op}
		\mathcal{L}_n = \theta^n + t P_1(\theta) + \cdots + t^k P_k(\theta) 
	\end{equation}
	for some $k \in \mathbb{N}$, where each $P_i(\theta)$ is a polynomial of degree $n$ in $\theta$. That $t=0$ is a large complex structure limit is encoded by the leading order term of $\mathcal{L}_n$ being $\theta^n$ in (\ref{eq-pf_op}); by the Frobenius method, the periods $\omega_0,\dots,\omega_{n-1}$ become a canonical basis of solutions of the form \begin{align}
		\label{eq-MUM_basis}
		\omega_0(t) &= h_0(t) \, , \\
		\omega_1(t) &= \frac{1}{2\pi i}\left(h_0(t) \log(t)+h_1(t)\right) \, , \nonumber \\
		\vdots & \nonumber \\
		\omega_{n-1}(t) &= \frac{1}{(2\pi i)^{n-1}}\left(\sum_{k=0}^{n-1} \binom{n-1}{k} h_k(t)\log^{n-1-k}(t) \right) \,, \nonumber 
	\end{align}
	where $h_0,\dots, h_{n-1}$ are holomorphic functions satisfying $h_0(0)=1, h_1(0)=h_2(0)= \dots = h_{n-1}(0) = 0$. Going back to the original construction of Candelas et al \cite{MR1101784}, the mirror map $\tau$ is defined concretely by
	\begin{equation}
		\label{eq-mirror_map}
		\begin{aligned}
			\tau :& \; U \to \mathbb{H}/\mathbb{Z} \, , \\
			\tau(t) =& \, \frac{\omega_1(t)}{\omega_0(t)} \, ,\\
		\end{aligned}
	\end{equation}
	which is effectively the truncated period morphism of the family $X_t$. The map $\tau$ can be extended to $t = 0$ by setting $\tau(0) = i \infty \in \overline{\mathbb{H} / \mathbb{Z}}$, where the latter symbol denotes the appropriate compactification. This extended map, denoted by the same character, is invertible in an analytic neighborhood $\overline{U} \subseteq \mathbb{P}^1 - \{0,\infty\}$ of $t=0$. 
	\begin{definition}
		\label{def-closed_mirror_map}
		In the extended setting, the mapping $\tau(t)$ or its inverse $t(\tau)$ are referred to as the \emph{closed} mirror map. It is convenient then to define $q = \exp \circ \, 2 \pi i\tau : U \to \mathbb{C}$, so that
		\begin{equation}
			\label{eq-q_coordinate}
			q(t) = \exp(2 \pi i\tau(t)) = \exp\left(2 \pi i \frac{\omega_1(t)}{\omega_0(t)} \right) = t \exp\left(\frac{h_1(t)}{h_0(t)}\right) 
		\end{equation} and $q(0)=0$. In cases relevant to mirror symmetry, we have $q(t) \in t \mathbb{Z}[[t]]$. We will refer to $q(t)$ as the \emph{exponentiated} mirror map. 
	\end{definition}
	\par Under Assumption \ref{ass-mirror}, $\tau = \tau_{\mathrm{eu}}$ as in (\ref{eq-I-mirror_maps}). Then the fundamental expression of mirror symmetry from this perspective is the encoding of the degree-$d$ instanton numbers $N_d$ (= genus zero GV invariants, for $\dim(X) = 3$ \cite{gv1998-1,gv1998-2}, analogously defined invariants for $\dim(X) = 4$ by Klemm \& Pandharipande \cite{MR2415462}) of $X^\vee$ through a generating function $\mathcal{Y}_{n-1}(t)$ known in physics as the B-side Yukawa coupling. It takes the form of a Lambert series expansion of the $q$-series for $\mathcal{Y}_{n-1}(t)$ near $t=0$ and can be computed from the Picard-Fuchs operator (\ref{eq-pf_op}) as follows. 
	\par By definition, the (algebraic) B-side Yukawa coupling $\mathsf{Y}_{n-1}(t)$ is given by 
	\begin{equation}
		\label{eq-algebraic_yukawa}
		\mathsf{Y}_{n-1}(t) := \int_{X_t} \Omega_t \wedge \frac{d^{n-1}}{dt^{n-1}} \Omega_t \, .
	\end{equation}
	Write the Picard-Fuchs operator (\ref{eq-pf_op}) as the monic ordinary differential operator
	\begin{equation}
		\label{eq-monic_diff_op}
		\mathcal{L}_n = \frac{d^n}{dt^n} + p_{n-1}(t)\frac{d^{n-1}}{dt^{n-1}} + \cdots + p_0(t) \, .
	\end{equation}
	It follows from a straightforward computation with $\mathcal{L}_n$ and Griffiths transversality\footnote{Which is vacuous in the case $n-1=1$ for elliptic curves and holds automatically in the case $n-1=2$ for K3 surfaces} that $\mathsf{Y}_{n-1}(t)$ satisfies a first order linear ODE of Wronskian type
	\begin{equation}
		\label{eq-yukawa_ODE}
		\frac{d\mathsf{Y}_{n-1}}{dt} = -\frac{2}{n} p_{n-1}(t)\mathsf{Y}_{n-1} \, .
	\end{equation}
	Hence by Abel's theorem, $\mathsf{Y}_{n-1}(t) \propto W(t)^{2/n}$, which is a priori an algebraic function in $t$ determined by the Wronskian $W(t)$. The algebraicity follows from the fact that $\mathcal{L}_n$ is a Fuchsian operator whose coefficients are rational in $t$.  In fact, $\mathsf{Y}_{n-1}(t)$ is a rational function in $t$ with a pole of order $n-1$ at $t=0$. This stems ultimately from the fact that the relevant Picard-Fuchs operator is of type $\mathrm{D}n$. 
	\begin{definition}
		\label{def-general_algebraic_yukawa}
		For an operator $\mathcal{L}_n$ of type $\mathrm{D}n$ whose monic form is given in (\ref{eq-monic_diff_op}), regardless of if it is of geometric origin or not, define the algebraic Yukawa coupling $\mathsf{Y}_{n-1}(t)$ as the solution of (\ref{eq-yukawa_ODE}), with the integration constant normalized to $1$.
	\end{definition}
	\begin{lemma}
		\label{lem-algebraic_yukawa}
		Suppose that $\mathcal{L}_n$ is an operator of type $\mathrm{D}n$, and that near $t=0$, the symbol\footnote{That is, the leading order coefficient of $\frac{d^n}{dt^n}$ for $\mathcal{L}_n$} of the operator is given as $t^{n}\mathcal{Q}(t)$, where $\mathcal{Q}(t)$ is relatively prime from $t$. Then the algebraic Yukawa coupling $\mathsf{Y}_{n-1}(t)$ of $\mathcal{L}_n$ is given as
		\begin{equation}
			\label{eq-Dn_algebraic_yukawa}
			\mathsf{Y}_{n-1}(t) = \frac{1}{t^{n-1}\mathcal{Q}(t)} \, .
		\end{equation}
	\end{lemma}
	\begin{proof}
		The symbol has the desired form since $\mathcal{L}_n$ is MUM at $t=0$. Writing $\mathcal{L}_n$ as a monic operator as in (\ref{eq-monic_diff_op}), one finds by carefully examining the $\mathrm{D}n$ condition that after partial fraction decomposition,
		\begin{equation}
			\label{eq-par_frac}
			-\frac{2}{n} p_{n-1}(t) = -\frac{(n-1)}{t} - \frac{\mathcal{Q}'(t)}{\mathcal{Q}(t)} \, .
		\end{equation}
		Integrating (\ref{eq-yukawa_ODE}) then gives (\ref{eq-Dn_algebraic_yukawa}).
	\end{proof}
	\par While the specific form of the algebraic Yukawa coupling $\mathsf{Y}_{n-1}(t)$ ultimately plays an important role for us, more generally the solution of (\ref{eq-yukawa_ODE}) is of fundamental importance for the geometric content carried by the operator $\mathcal{L}_n$ and can be seen as follows. That $\mathcal{L}_n$ underlies a polarized VHS of a pencil of CY varieties  with $h^{n-2,1}=1$\footnote{The appropriate analogue of this condition for lattice polarized K3 surfaces is that a very general member of the pencil be of Picard rank-19} forces the monodromy group of the operator to be contained in an appropriate integral symplectic or special orthogonal group of certain signature, depending on the parity of $n$; e.g.,  $\mathrm{SO}(2,1;\mathbb{Z})$ and $\mathrm{Sp}(4,\mathbb{Z})$ for $n=3,4$. This  necessitates the existence of a function $\mathsf{Y}_{0}(t)$ satisfying (\ref{eq-yukawa_ODE}), such that under  conjugation the operator $\mathsf{Y}_{0}^{-1}\mathcal{L}_n\mathsf{Y}_{0}$ is (anti)self-adjoint, again depending on parity \cite[Rmk. 3.6]{MR4069107}. For $\mathcal{L}_n$ a CY Picard-Fuchs operator of the type above, we necessarily have $\mathsf{Y}_{0}(t) = \mathsf{Y}_{n-1}(t)$ from (\ref{eq-algebraic_yukawa}) up to scale \cite[\S 4.1]{MR1101784}.
	\par When $n = 2$, we have $\mathsf{Y}_0(t) = W(t)$, and we arrive at the well known fact that every such second order operator is formally self-adjoint. For $n \geq 3$, the condition of $\mathcal{L}_n$ being gauge equivalent to one that is (anti)self-adjoint can be detected from the vanishing of the first Wilczynski invariant $\mathsf{w}_1(\mathcal{L}_n) = 0$, \cite{MR1197867,MR1779161,MR0131232}. For $n = 3,4$, we have
	\begin{equation*}
		\begin{aligned}
			\mathsf{w}_1(\mathcal{L}_3) &= -2 p_0(t)+\frac{2}{3} p_1(t) p_2(t)-\frac{4}{27} p_2(t)^3-\frac{1}{3}  p_2''(t)+ p_1'(t)-\frac{2}{3} p_2'(t) p_2(t), \\
			\mathsf{w}_1(\mathcal{L}_4) &= -p_1(t)+\frac{1}{2} p_2(t) p_3(t)-\frac{1}{8} p_3(t)^3-\frac{1}{2} p_3''(t)+ p_2'(t)-\frac{3}{4} p_3'(t) p_3(t) .
		\end{aligned}
	\end{equation*}
	Note that this condition alone implies that the monodromy group of $\mathcal{L}_n$ is contained only in the complexified version of the appropriate aforementioned group. The condition $\mathsf{w}_1(\mathcal{L}_3) = 0$ is well-known to be both necessary and sufficient for $\mathcal{L}_3 = \mathcal{L}_2^{\otimes 2}$ to be the symmetric square of a second order operator, whence the monodromy representation of $\mathcal{L}_3$ coincides with the symmetric square representation for $\mathcal{L}_2$ of $\mathrm{SL}(2,\mathbb{C}) \cong \mathrm{Sp}(2,\mathbb{C})$ in $\mathrm{SO}(2,1;\mathbb{C})$. This is in perfect agreement with well-known results  on the Picard-Fuchs operators of pencils of lattice polarized K3 surfaces Picard rank-19 decomposing as the symmetric square of a second order operator \cite[\S 5.2]{MR1738862}, \cite[\S 3.1]{MR1779161}, and corresponds to the Hodge theoretic fact that the period domain of rank-19 lattice polarized K3 surfaces is a Zariski open subset of a conic in $\mathbb{P}^2$. More generally, the period domain of lattice polarized K3 surfaces of Picard rank-$k$ is a quasiprojective subvariety of a quadric hypersurface in $\mathbb{P}^{21-k}$ \cite{MR1420220}, which has significant implications for the geometry of the associated moduli spaces and full Picard-Fuchs PDE systems \cite{MR4928995,MR1136204,MR960834,MR996019,schultz_geometry_2021}. We will comment on the geometric implications of this decomposition for the K3 surfaces themselves momentarily.
	\par In light of the discussion in \S \ref{ss-I_functions}, it is appropriate to recalculate the algebraic Yukawa coupling $\mathsf{Y}_{n-1}(t)$ after the projective gauge transformation $y(t) \mapsto y(t)/\omega_0(t)$, where $y(t)$ is the dependent variable in the Picard-Fuchs ODE. Pulled back by the local inverse $t(\tau)$ of the mirror map, one obtains the following (normalized) Yukawa coupling. 
	\begin{definition}
		\label{def-normalized_yukawa_instanton}
		The normalized Yukawa coupling of the  pencil $\mathcal{X} \to \mathbb{P}^1_t$ of projective CY $(n-1)$-folds is given in an analytic neighborhood of the LCSL $t=0$ by
		\begin{equation}
			\label{eq-normalized_Yukawa}
			\mathcal{Y}_{n-1}(\tau):=\frac{1}{\omega_0(\tau)^2}\left(\frac{1}{2 \pi i} \frac{d t}{d \tau}\right)^{n-1} \mathsf{Y}_{n-1}(\tau)  \, . 
		\end{equation}
	When $n-1 = \dim(X_t) \geq 3$, the degree $d$-instanton numbers $N_d$, $d \geq 1$ are defined by the Lambert series expansion 
		\begin{equation}
			\label{eq-instanton_numbers}
			\mathcal{Y}_{n-1}(\tau) = \sum_{m=0}^{\infty} c_m q^m =: \kappa +\sum_{d=1}^{\infty} d^j N_d \frac{q^d}{1-q^d} \, ,
		\end{equation}
	for an appropriate integer $j \geq 1$; it is taken as $j = 3,2,1$ for $\dim(X) = 3,4,5$. Here $q = \exp(2\pi i \tau)$ is the exponentiated mirror map.
	\end{definition}
	\par That the expression in (\ref{eq-normalized_Yukawa}) appears is dependent entirely on the CY condition holding from the vanishing of Wilczynski invariant $\mathsf{w}_1(\mathcal{L}_n)=0$ \cite[Prop. 1, Rmk. (3.8)]{MR2282972}. This is elucidated in Theorem \ref{thm-3-point_invariants} below, which also explains the following. In addition to the vanishing of $\mathsf{w}_1(\mathcal{L}_n)$, the CY condition assumes strong integrality properties hold for the quantities above, $t(\tau), \,  \mathcal{Y}(\tau) \in \mathbb{Z}[[q]]$. These conditions are reasonable, given that the instanton numbers conjecturally compute the degree-$d$, 3-point invariants $N_d \in \mathbb{Z}$ of $X^\vee_\tau$. The integer $j$ appearing in (\ref{eq-instanton_numbers}) comes from the multicovering formula from 3-point GW invariants
	\begin{equation}
		\label{eq-covering_formula}
		N_d=\frac{1}{d^j} \sum_{m \mid d} \mu\left(\frac{d}{m}\right) c_m \in \mathbb{Z}
	\end{equation}
	where $\mu$ is the M\"obius function. The threefold case is the environment for classical mirror symmetry, in which  (\ref{eq-normalized_Yukawa}) - the so-called B-side Yukawa coupling determined from $X_t$ - is supposed to be identified with an analogous quantity called the A-side Yukawa coupling that encodes the structure constants of the small quantum cohomology ring of $X^\vee_\tau$. Famously, the calculation of the instanton numbers as agreeing with the degree-$d$, 3-point invariants holds for the quintic threefold \cite{MR1204770,MR1101784,MR1191426}, as well as other CY-3 pencils with $h^{1,1} = 1$ for which established mirror constructions hold. In this sense, we have the following.
	\begin{theorem}
		\label{thm-3-point_invariants}
		If $\mathcal{X} \to \mathbb{P}^1_t$ is a pencil of CY varieties of dimension-$(n-1)$ mirror to the anticanonical linear system $\mathcal{X}^\vee \to \mathbb{P}^1_{\tau}$ of an almost minimal Fano variety $V$ of dimension $n \geq 3$, then under the mirror map $t(\tau)$ in (\ref{eq-mirror_map}), in a neighborhood of the LCSL $t=0$, the 3-point function $\mathcal{Y}^{(3)}_{\mathrm{P}^{n-3}}(X^\vee_{\tau})$ of Theorem \ref{thm-CY_I-function} is given by $\mathcal{Y}_{n-1}(\tau)$ in (\ref{eq-normalized_Yukawa}). Thus in this case, the degree-$d$ instanton numbers $N_d$ in (\ref{eq-instanton_numbers}) count genuine 3-point invariants of $X^\vee_\tau$, under the identification (\ref{eq-covering_formula}) for $m \geq 1$,
		\begin{equation}
			\label{eq-3-point_invariants}
			c_m = \kappa \langle \mathrm{P},\mathrm{P},\mathrm{P}^{n-3}\rangle_{0,3,m} \, .
		\end{equation}
	\end{theorem}
	\begin{proof}
		Under the hypotheses we have assumed in Assumption \ref{ass-mirror}, the $z^{-k}$ component of the $I$-function $I_{X^\vee}(\tau,z)$ in (\ref{eq-I-function_expansion}) are identified under composition with the mirror map as the period $\omega_k(t)$ of $X_t$ as in (\ref{eq-MUM_basis}), up to a scalar correction factor which cancels in the ratio defining the mirror map. Such periods are solutions of the Picard-Fuchs operator $\mathcal{L}_n$, which by Theorem \ref{thm-operator_factorization} is (anti) self-adjoint, and hence $\mathsf{w}_1(\mathcal{L}_n)=0$. An argument similar to the one presented by van Straten \& Zudilin in \cite[Prop. 1, Rmk. (3.8)]{MR2282972} for the case $n-1=3$ shows that 
		\begin{equation}
			\label{eq-period_ratio}
			\kappa\frac{d^2}{d\tau^2} \left(\frac{\omega_2(\tau)}{\omega_0(\tau)}\right) = \mathcal{Y}_{n-1}(\tau) \,,
		\end{equation}
		which by the identification with $I$-function components above, yields the desired result.
	\end{proof}
	\subsubsection{Projective CY threefolds}
	\label{sss-projective_CY-threefolds}
	\par The Yukawa coupling is best understood in the case of mirrors $(X^\vee_\tau,X_t)$ of projective CY threefolds \cite{MR1101784}. As above consider $X^\vee$ as an anticanonical CY in an almost minimal Fano fourfold $V$. As is well known, the A-side Yukawa coupling is determined in a neighborhood of the large radius limit $\tau = i\infty \in \mathbb{P}^1_{\tau}$ by an analytic function of the form (as usual $q = \exp(2\pi i \tau)$)
	\begin{equation*}
		\mathcal{F}=\frac{1}{6} \kappa \tau^3+\frac{1}{2} \kappa_2 \tau^2+\frac{1}{2} \kappa_1 \tau+\frac{1}{6} \kappa_0 + \frac{1}{(2\pi i)^3} \sum_{m=1}^{\infty} a_m q^m
	\end{equation*}
	which is known in physics as the holomorphic prepotential. The cubic terms are related to the topology of $X^\vee$ as 
	\begin{equation}
		\label{eq-A-side_topological_data}
		\kappa=\int_{X^{\vee}} \mathrm{P} \wedge \mathrm{P} \wedge \mathrm{P} \, , \quad \kappa_0=3 e\left(X^{\vee}\right) \frac{\zeta(3)}{(2 \pi i)^3}
	\end{equation}
	so that $\kappa$ is the triple intersection number of the positive, primitive generator $\mathrm{P}$ of $H^2(X^\vee,\mathbb{Z})$ (the hyperplane class), and $e(X^\vee)$ is the topological euler characteristic. The other numbers $\kappa_1,\kappa_2$ are more subtle. For the quintic threefold in $\mathbb{P}^4$, the former is the degree $\kappa=5$. The prepotential determines the A-side Yukawa coupling $\mathcal{Y}^{(3)}_\mathrm{P}(X^\vee)$ as
	\begin{equation}
		\label{eq-prepotential_yukawa}
		\frac{d^3\mathcal{F}}{d\tau^3} = \mathcal{Y}^{(3)}_\mathrm{P}(X^\vee) \, .
	\end{equation}
	\par Suppose that $\mathcal{L}_4$ is the Picard-Fuchs operator of the mirror $X_t$. Then near the LCSL, there is a canonical Frobenius basis of solutions $\omega_0(t),\dots,\omega_3(t)$. Mirror symmetry predicts that the holomorphic prepotential $\mathcal{F}$ is expressible in the form
	\begin{equation*}
		\mathcal{F}=-\frac{1}{12} \kappa\left(-\frac{\omega_3}{\omega_0}+\frac{\omega_1 \omega_2}{\omega_0^2}\right)-\frac{1}{2} \kappa_2 \frac{\omega_1^2}{\omega_0^2}-\frac{1}{2} \kappa_1 \frac{\omega_1}{\omega_0}-\frac{1}{6} \kappa_0 .
	\end{equation*}
	This is at least a bit reasonable, given the mirror map $\tau = \frac{\omega_1}{\omega_0}$. On the other hand, $\mathcal{F}$ determines the periods as
	\begin{equation*}
		(\omega_0,\omega_1,\omega_2,\omega_3) = \left(\omega_0,\omega_0\cdot \tau,\omega_0 \cdot \frac{d\mathcal{F}}{d\tau}, \omega_0 \cdot\left(\tau\frac{d\mathcal{F}}{d\tau} - 2\mathcal{F}\right)\right) \, ,
	\end{equation*}
	and by the CY condition $\mathsf{w}_1(\mathcal{L}_4)=0$, the Picard-Fuchs operator transforms under the mirror map as $\mathcal{L}_4 \mapsto \mathcal{L}_{4}^{(\tau)}$, where
	\begin{equation}
		\label{eq-CY3_PF_tau}
		\mathcal{L}_{4}^{(\tau)} = \frac{1}{\mathcal{R}(\tau)} \frac{d^2}{d\tau^2} \circ \left(\frac{1}{\mathcal{Y}_{3}(\tau)}\frac{d^{2}}{d\tau^2}\right) \frac{1}{\omega_0(\tau)} \, .
	\end{equation}
	Here $\mathcal{Y}_3(\tau)$ is the B-side Yukawa coupling (\ref{eq-normalized_Yukawa}), which by Theorem \ref{thm-3-point_invariants} is equal to $\mathcal{Y}^{(3)}_\mathrm{P}(X^\vee)$. The quantity $\mathcal{R}(\tau)$ is related to the symbol of $\mathcal{L}_4$. It will play a significant role in \S \ref{s-twisted_invariants} as a generating function for 3-point euler-twisted invariants.
	\subsection{Mirrors of elliptic curves \& K3 surfaces}
	\label{ss-EC_K3_mirror}
	\par As discussed in \S \ref{ss-genus_zero_invariants}, the ordinary genus zero Gromov-Witten theory of elliptic curves and K3 surfaces are respectively vacuous and trivial. This necessitates that the associated normalized Yukawa couplings $\mathcal{Y}_{1}, \mathcal{Y}_2$ must be trivial; as opposed to the well-known CY threefold case, the constancy of these quantities alone is not sufficient to conclude that the associated genus zero invariants vanish. Yet there is still a well-defined notion of a mirror variety for both elliptic curves \cite{MR1363055,MR3652463,MR3822918} and K3 surfaces \cite{MR1420220,MR4913468,MR1779161}, and we will see that the necessity of this statement that  $\mathcal{Y}_{1}, \mathcal{Y}_2$ are trivial  will play an important role in the context of this article. 
	\par For ease of notation in what follows, let us write in this section only monic second and third order operators $\mathcal{L}_2, \, \mathcal{L}_3$ as 
	\begin{equation*}
		\begin{aligned}
			\mathcal{L}_2 &= \frac{d^2}{dt^2} + a_1(t)\frac{d}{dt} + a_0(t) \, , \\
			\mathcal{L}_3 &= \frac{d^3}{dt^3} +  b_2(t)\frac{d^2}{dt^2} + b_1(t)\frac{d}{dt} + b_0(t) \, .
		\end{aligned}
	\end{equation*}
	The third order symmetric square $\mathcal{L}_3$  is determined by the second order operator via the following relations:
	\begin{equation}
		\label{eq-sym_square_coeffs}
		b_2(t) = 3a_1(t) \, , \quad b_1(t) = a_1'(t) + 2a_1(t)^2 + 4a_0(t) \, , \quad b_0(t) = 2a_0'(t) + 4a_1(t)a_0(t) \, .
	\end{equation}
	Furthermore, if $\eta_0(t), \eta_1(t)$ are a local basis of solutions for $\mathcal{L}_2$, then a local basis of solutions $\mathcal{L}_3$ is given by $\omega_0(t) = \eta_0(t)^2$, $\omega_1(t) = \eta_0(t)\eta_1(t)$, and $\omega_2(t) = \eta_1(t)^2$. Hence the solutions satisfy the quadratic relation $\omega_0\omega_2 - \omega_1^2 = 0$.
	\par These equations carry geometric information about the generic elements of the underlying K3 pencil, and is intimately related to the existence of rational double covers of Picard rank-19 Kummer surfaces \cite[\S 3.3.3]{MR4494119}, \cite[\S 6.2.3]{schultz_geometry_2021}. For example, the pencils of $\mathrm{M}_N$-polarized K3 surfaces, mirror to anticanonical K3 surfaces in rank-1 Fano threefolds,  admit a canonical Shioda-Inose structure \cite[\S 7]{MR1420220}, and a generic member of the pencil admits a rational double cover of the Kummer surface $\mathrm{Kum}(E \times E')$, where $E$ and $E'$ are $N$-isogenous. This played an important role in Golyshev's analysis \cite{MR2306141}, where many geometric quantities relevant to mirror symmetry are expressible in terms of elliptic modular forms determined by the underlying K3 moduli space($\cong X_0(N)^+$). It follows from Doran \cite[Thm. 5.13]{MR1738862} that the K3 Picard-Fuchs operator $\mathcal{L}_3$ is the symmetric square of the underlying second order Picard-Fuchs operator $\mathcal{L}_2$ of the associated elliptic pencil whose generic fibre determines the Kummer surface. The quadratic period relation $\omega_0\omega_2 - \omega_1^2 = 0$ is necessitated by the K3 period domain lying on a quadric surface in $\mathbb{P}^3$ \cite[\S 3]{MR1420220}. Moreover, an important consequence observed by Doran \cite[Thm. 6]{MR1779161} can be summarized in the following. 
	\begin{proposition}
		\label{prop-square_mirror_map}
		The K3 mirror map $\tau_2(t)$ is projectively equivalent to the underlying elliptic curve mirror map $\tau_1(t)$.
	\end{proposition}
	\begin{proof}
		Observe:
		\begin{equation}
			\label{eq-K3_mirror_map}
			\tau_2(t) = \frac{\omega_1(t)}{\omega_0(t)} = \frac{\eta_0(t)\eta_1(t)}{\eta_0(t)^2} = \frac{\eta_1(t)}{\eta_0(t)} = \tau_1(t) \, .
		\end{equation}
	\end{proof}
	\par We conclude with the following elementary result on the normalized Yukawa couplings $\mathcal{Y}_{1}, \, \mathcal{Y}_{2}$ of elliptic curves and K3 surfaces. 
	\begin{proposition}
		\label{prop-EC_K3_Yukawa}
		For both $n-1 = 1, 2$, the normalized Yukawa coupling is identically constant, $\mathcal{Y}_{n-1}(\tau) \equiv \mathrm{const}$. Thus for all $m \geq 1$, the terms $c_m = 0$ in the $q$-expansion (\ref{eq-instanton_numbers}) vanish.
	\end{proposition}
	\begin{proof}
		It follows from (\ref{eq-yukawa_ODE}) that when $n-1=1$, $\mathsf{Y}_1(t)=C \cdot W(t)$ is simply a constant multiple the Wronskian $W(t)$ of the Picard-Fuchs operator $\mathcal{L}_2$. We take $C=1$. Using the notation in the Proposition \ref{prop-square_mirror_map}, the corresponding normalizing coefficient in (\ref{eq-normalized_Yukawa}) yielding $\mathcal{Y}_1(\tau)$ is precisely the reciprocal $1/W(\tau)$, since from (\ref{eq-mirror_map}) we have
		\begin{equation}
			\label{eq-wronskian_eq}
			\frac{d\tau_1}{dt} = \frac{W(t)}{\eta_0(t)^2} \, .
		\end{equation}
		By the inverse function theorem, $\mathcal{Y}_1(\tau) \equiv \mathrm{const}$.
		\par In the case $n-1=2$, applying (\ref{eq-yukawa_ODE}) with the symmetric square relations (\ref{eq-sym_square_coeffs}) for $\mathcal{L}_3 = \mathcal{L}_2^{\otimes 2}$, we have (up to a constant multiple, which we again normalize to be 1) that $\mathsf{Y}_2(t) = W(t)^2$, where $W(t)$ is the Wronskian of the second order operator $\mathcal{L}_2$. It follows from (\ref{eq-K3_mirror_map}) that 
		\begin{equation}
			\label{eq-K3_mirror_map_derivative}
			\frac{d \tau_2}{dt} = \frac{d \tau_1}{dt} \, ,
		\end{equation}
		so evaluating (\ref{eq-normalized_Yukawa}) we find the normalized Yukawa coupling $\mathcal{Y}_2(\tau) \equiv \mathrm{const}$. 
	\end{proof}
	\begin{remark}
		\label{rmk-elliptic_curve_Yuk}
		The argument above is consistent with Theorem \ref{thm-CY_I-function}a) and \S \ref{ss-genus_zero_invariants}. It also generalizes the well-known calculations of Nagura \& Sugiyama \cite[\S 3.2, \S 6]{Nagura1995MirrorSymmetry}. Moreover, Proposition \ref{prop-EC_K3_Yukawa} shows that Theorem \ref{thm-3-point_invariants} can be formally extended to the elliptic curve case $n-1=1$ by identifying
		\begin{equation}
			\label{eq-elliptic_curve_yukawa}
			\mathcal{Y}^{(3)}_{\mathrm{P}^{-1}}(X^\vee_{\tau}) \equiv \mathcal{Y}_1(\tau) \, ,
		\end{equation} 
		which is consistent with Definition \ref{def-3-point_function}. This will be important in the sequel.
	\end{remark}
	\subsection{Local mirror symmetry}
	\label{ss-local_ms}
	\par Local mirror symmetry adapts the techniques of closed mirror symmetry as in \S \ref{ss-closed_ms} to study the situation of the geometry of Fano subvariety $V \subset Y^\vee$ in an ambient CY variety $Y^\vee$ \cite{MR1797015}. By adjunction, it follows that the normal bundle $N_{V/Y^\vee} \cong K_V$, and so local mirror symmetry is intrinsic to $V$, and is thus modeled by the open CY variety $Y^\vee = \mathrm{Tot}(K_V)$. Suppose that $\mathcal{W}$ is a superpotential for $V$. Then the mirror construction a la Hori \& Vafa \cite{hori2000mirrorsymmetry} produces a mirror $Y$ of $Y^\vee$: it is an affine conic bundle defined over the level sets of $\mathcal{W}$, and is hence a quasiprojective variety of dimension-$(n+1)$ equipped with a \emph{meromorphic} $(n+1,0)$ form with nonvanishing (constant) residue obtained from integration over the fibre directions. Accordingly, the base of the conic bundle is precisely the mirror family $\mathcal{X}$ of a generic, smooth anticanonical divisor $X^\vee_\tau \in |-K_V|$ as in \S \ref{ss-closed_ms}, and the reduction of the meromorphic $(n+1,0)$ form yields the period integrals of $X$. A celebrated paper of Abouzaid, Aroux, \& Katzarkov \cite{MR3502098} shows that the local mirrors in this sense satisfy a version of mirror symmetry known as SYZ mirror symmetry \cite{MR1429831}.
	\par The regularized local $\mathrm{SQDM}^{\mathrm{reg}}(K_V)$ \cite[\S 6]{MR3519131} is then represented by the operator 
	\begin{equation}
		\label{eq-local_operator}
		\mathcal{D}^{\mathrm{loc}}_{n+1} = \mathcal{L}_n \circ \, \theta \, .
	\end{equation}
	It is defined over the same complex structure moduli space as the compact part of the mirror, and it is easy to see that the $\mathcal{D}_{n+1}^{\mathrm{loc}}$ is also MUM at $t = 0 \in \overline{U}$. The local operator, sometimes referred to in the literature as the Lerche-Mayr Picard-Fuchs operator, annihilates the periods $\eta_0,\eta_1,\dots,\eta_{n}$ of the mirror variety $Y$, which can easily be seen to have the following form from the discussion above: 
	\begin{equation}
		\label{eq-local_periods}
		\begin{aligned}
			\eta_0 & = 1 \, , \\
			\eta_j &= \frac{1}{2\pi i} \int_{t_0}^t \frac{\omega_{j-1}(s)}{s} \, ds 
		\end{aligned}
	\end{equation}
	for $j = 1,\dots,n$, for a fixed choice of $t_0 \in \overline{U}-\{0\}$ such that (\ref{eq-local_periods}) forms a Frobenius basis of $\mathcal{D}_{n+1}^{\mathrm{loc}}$ and integration is carried out with respect to a given fixed branch of the logarithm in $\overline{U}-\{0\}$. The mirror map $\tau_{\mathrm{loc}}$ may be computed in exactly the same manner as before (\ref{eq-mirror_map}, \ref{eq-q_coordinate}).
	\begin{definition}
		\label{def-local_mirror_map}
		The mapping $\tau_{\mathrm{loc}}(t)$ given by 
		\begin{equation}
			\label{eq-local_mirror_map_0}
			\tau_{\mathrm{loc}}(t) := \frac{\eta_1(t)}{\eta_0(t)} = \frac{1}{2\pi i}\int_{t_0}^t \frac{\omega_0(s)}{s} \, ds \, 
		\end{equation}
		is invertible in a sufficiently small analytic neighborhood $\widetilde{U}$ of $t=0$ in $\overline{U}$, and is called together with its inverse $t(\tau_{\mathrm{loc}})$ the \emph{local} mirror map. By shrinking $\widetilde{U}$ if necessary, we may assume that $\widetilde{U} = \overline{U}$. Define $Q = \exp \circ \, 2 \pi i\tau_{\mathrm{loc}} : \overline{U} \to \mathbb{C}$, so that
		\begin{equation}
			\label{eq-Q_coordinate}
			Q(t) := \exp(2 \pi i\tau_{\mathrm{loc}}(t)) = \exp\left(\int_{t_0}^t \frac{\omega_0(s)}{s} \, ds\right) 
		\end{equation} and $Q(0)=0$. We will refer to $Q(t)$ as the \emph{exponentiated} local mirror map. 
	\end{definition}
	\begin{remark}
		\label{rmk-B-field}
		The choice of $t_0 \in \overline{U}-\{0\}$ is not immaterial, and related to normalization of the periods and determining the so called $B$-field of the local variety.
	\end{remark}
	\par Under Assumption \ref{ass-mirror}, $\tau_{\mathrm{loc}}$ agrees with the local mirror map in (\ref{eq-I-mirror_maps}). The local Yukawa coupling and genus zero invariants are defined in the exact same manner as in \S \ref{ss-closed_ms}.
	\begin{definition}
		\label{def-local_yukawa}
		The local algebraic Yukawa coupling $\mathsf{Y}_{n+1}^{\mathrm{loc}}(t)$ is defined by the solution to the analogous equation to (\ref{eq-yukawa_ODE}) obtained from $\mathcal{D}_{n+1}^{\mathrm{loc}}$ with indices appropriately shifted. The normalized Yukawa coupling $\mathcal{Y}_{n+1}^{\mathrm{loc}}(\tau_{\mathrm{loc}})$ is given by
		\begin{equation}
			\label{eq-normalized_local_yukawa}
			\begin{aligned}
				\mathcal{Y}_{n+1}^{\mathrm{loc}}(\tau_{\mathrm{loc}}) :=& \frac{1}{\eta_0(\tau_{\mathrm{loc}})^2}\left(\frac{1}{2 \pi i} \frac{d t}{d \tau_{\mathrm{loc}}}\right)^{n+1}\mathsf{Y}_{n+1}^{\mathrm{loc}}(\tau_{\mathrm{loc}}) \\  \equiv& \left(\frac{1}{2 \pi i} \frac{d t}{d \tau_{\mathrm{loc}}}\right)^{n+1} \mathsf{Y}_{n+1}^{\mathrm{loc}}(\tau_{\mathrm{loc}})  \, . 
			\end{aligned}
			\end{equation}
			When $n+1 = \dim(Y) \geq 3$, the degree-$d$ local instanton numbers $N_d$, $d \geq 1$ are defined by the Lambert series expansion 
			\begin{equation}
			\label{eq-local_instanton_numbers}
			\mathcal{Y}_{n+1}^{\mathrm{loc}}(\tau_{\mathrm{loc}}) = \sum_{m=0}^{\infty} b_m Q^m =: \kappa^{\mathrm{loc}}+\sum_{d=1}^{\infty} d^j N_d \frac{Q^d}{1-Q^d} \, ,
			\end{equation}
			for an appropriate integer $j \geq 1$; it is again taken as $j = 3,2,1$ for $\dim(Y) = 3,4,5$. Here $Q = \exp(2\pi i \tau_{\mathrm{loc}})$ is the exponentiated local mirror map.
	\end{definition}
	Notice that from the mirror map $\tau(t)$ in (\ref{eq-mirror_map}), the parameter $t$ can be eliminated, and then we have by elementary calculus that 
	\begin{equation}
		\label{eq-local_mirror_map_derivative}
		\frac{d\tau_{\mathrm{loc}}}{d\tau}  =  \frac{1}{2\pi i} \frac{\omega_0(\tau)}{t(\tau)} \frac{dt}{d\tau}  \, .
	\end{equation} This mirror map derivative will arise naturally in the sequel as a generating function of 3-point euler-twisted invariants of the ambient Fano variety. In the local threefold case, the $Q$-series of the reciprocal mirror map derivative is known to play a fundamental role in determining the local genus zero invariants.
	\subsubsection{Local CY threefolds}
	\label{sss-local_CY-3}
	\par For local threefolds of the form $Y^\vee = \mathrm{Tot}(K_S)$, for $S$ an almost minimal del Pezzo, the local Yukawa coupling is also particularly well understood, in analogy to the case of projective CY threefolds in \S \ref{sss-projective_CY-threefolds}. It is known that the local Yukawa coupling $\mathcal{Y}_{3}^{\mathrm{loc}}(\tau_{\mathrm{loc}})$ is again determined by a holomorphic prepotential $\mathcal{F}^{\mathrm{loc}}(\tau_{\mathrm{loc}})$. However in the local case, since the periods fail to be homogeneous, in physics one defines the prepotential implicitly through the relation
	\begin{equation}
		\label{eq-local_prepotential}
		\mathcal{Y}_{3}^{\mathrm{loc}}(\tau_{\mathrm{loc}}) = \frac{d^3\mathcal{F}^{\mathrm{loc}}}{d\tau_{\mathrm{loc}}^3} = \kappa^{\mathrm{loc}}\frac{d^2}{d\tau_{\mathrm{loc}}^2} \left(\frac{\eta_2}{\eta_0}\right) = \kappa^{\mathrm{loc}}\frac{d}{d\tau_{\mathrm{loc}}}\left(\frac{\omega_1}{\omega_0}\right) =  \kappa^{\mathrm{loc}}\frac{d\tau}{d\tau_{\mathrm{loc}}} \, .
	\end{equation}
	This relation has been mathematically verified in certain cases - more precisely, for $V$ an almost minimal toric Fano surface - to be expressed solely in terms of the mirror map derivative $\frac{d\tau}{d\tau_{\mathrm{loc}}}$. The following result for local threefolds $Y^\vee = \mathrm{Tot}(K_V)$ is described in \cite[\S 2.3]{ms2024}. It is a summary of a statement first conjectured by Hosono \cite{MR2282969}, and proved by Bloch, Kerr, \& Vanhove \cite[Thm. 6.1]{MR3780269} on the local GW potential for $Y^\vee_{\tau} = \mathrm{Tot}(K_V)$, and combined with results of Doran \& Kerr \cite[\S 5]{MR2851155}.
	\begin{theorem}
		\label{thm-bkv}
		[Bloch, Kerr, Vanhove; Doran \& Kerr] Let $V$ be an almost minimal toric Fano surface corresponding to an LG model $\mathcal{W}$, and $E_t$ be the corresponding mirror family of elliptic curves to anticanonical elliptic curves $\check{E}_{\tau} \in |-K_V|$. Then the functional inverse $\tau(\tau_{\mathrm{loc}})$ of mirror maps determines the local invariants of $K_V$ as 
		\begin{equation}
			\label{eq-inversion_eq2}
			 \kappa^{\mathrm{loc}}\frac{d\tau}{d\tau_{\mathrm{loc}}} = \mathcal{Y}_{3}^{\mathrm{loc}}(\tau_{\mathrm{loc}}) = \kappa^{\mathrm{loc}} + \sum_{d=1}^{\infty} d^3 n_d Q^d \, ,
		\end{equation}
		where $\kappa^{\mathrm{loc}} \in \mathbb{Q}$ and $n_d$ are the degree-$d$ genus zero invariants of $Y^\vee = \mathrm{Tot}(K_V)$.\footnote{The normalization presented in the $Q$-expansion of (\ref{eq-inversion_eq2}) differs by a factor of $-1$ from that in \cite{MR3780269}.}
	\end{theorem}
	The goal of this article is to show that the non-equivariant limit of quantum Serre duality a la Iritani, Mann, \& Mignon \cite{MR3519131} allows one to greatly generalize this story. 
	\section{Reciprocal mirror map derivatives as generating functions of euler-twisted \& local 3-point invariants}
	\label{s-twisted_invariants}
	%
	%\par To do (Both): carefully formulate statement about local invariants and their expression in terms of the twisted invariants. Likely will need separate results for threefolds and fourfolds. LATER: it seems that a uniform proof can more or less be given. 
	%
	%Know: for a toric fano surface $F$, that $Q\frac{d}{dQ}\log q$ is the generating function for local $GV_0$ invariants of $K_F$ (Bloch, Kerr, Vanhove).
	
	%Also, it appears that $q\frac{d}{dq}f(q)$ = logarithmic derivative of the twisted invariants = $\pm q \frac{d}{dq}\log(Q)$ is the functional inverse of this quantity. Notice that $q$ is the exponentiated mirror map coming from the divisor, and $Q$ is the exponentiated mirror map coming from the local variety.

	%
	
	\par The non-equivariant limit of QSD restricted to the SQDM of the anticanonical Calabi-Yau hypersurface $X^\vee \in |-K_V|$ in an $n$-dimensional almost minimal Fano as described in \S \ref{ss-QSD} has a very simple expression in terms of exponentiated mirror maps and 2-point twisted invariants of $V$. Let $f(q)$ be the generating function of the 2-point euler-twisted invariants of $V$ that is obtained from the non-equivariant limit, 
	\begin{equation}
		\label{eq-twisted_gf}
		f(q) = \log(q)+\sum_{d=1}^{\infty} q^{d}\left\langle \mathrm{P}^{n-1}, \widetilde{1}\right\rangle_{0,2, d}^{\left(\mathrm{e}, \,  K_V^{-1}\right)} \, ,
	\end{equation}
	where $H^2(V,\mathbb{Z}) = \mathrm{span}_{\mathbb{Z}}\{\mathrm{P}\}$. Here $\widetilde{1}$ is meant to be understood as insertion of the class $\mathbf{e}(K_V^{-1})^{-1} = (c_1(V))^{-1}$ of degree $-2$. 
	The following result from Iritani, Mann, \& Mignon is crucial for the remainder of our analysis:
	\begin{lemma}
		\label{lem-exp_mirror_maps} 
		The exponentiated closed and local mirror maps are intertwined by the 2-point euler-twisted generating function $f(q)$ as
		\begin{equation}
			\label{eq-QSD_mirror}
			Q(q) = (-1)^{\delta}\exp(f(q)) \, 
		\end{equation}
		where $\delta$ is the Fano index of $V$.
	\end{lemma}
	\begin{proof}
		This is simply the exponentiated form of \cite[Lem. 5.18]{MR3519131}, up to a parity shift in the Novikov variable. 
	\end{proof}
	This very simple equation (\ref{eq-QSD_mirror}) is the fundamental expression of QSD that is relevant to the remainder of the article. It captures the essence of the QSD commutative diagram in Corollary (\ref{cor-IMM_2}) in the context of the mirror maps given in \S \ref{ss-I_functions}, \ref{ss-closed_ms}, \ref{ss-local_ms}. 
	\subsection{The mirror map derivative \texorpdfstring{$\frac{d\tau_{\mathrm{loc}}}{d\tau}$ as generating function of 3-point euler-twisted invariants}{d_\tau_{\mathrm{loc}}/d\tau}}
	\label{ss-mirror_map_derivative}
	From Lemma \ref{lem-exp_mirror_maps}, we obtain from the following result that expresses the regularized SQDM of the almost minimal Fano variety $V$ in terms of 3-point euler-twisted invariants.
	\begin{theorem}
		\label{thm-3-point_euler-twisted}
		Let $X^\vee_\tau \in |-K_V|$ be a smooth, anticanonical CY variety in a smooth, almost minimal Fano variety $V$ of dimension $n \geq 2$. Let $\tau,\tau_{\mathrm{loc}}$ be the closed and local mirror maps as in \S \ref{ss-closed_ms}, \ref{ss-local_ms}, convergent in an analytic neighborhood of the large volume limit for $X^\vee_\tau$.
		\begin{enumerate}
			\item[a)] A generating function of 3-point euler-twisted invariants is given in terms of the derivative of the local mirror map as
			\begin{equation}
				\label{eq-3-point_euler-twisted}
				\frac{d \tau_{\mathrm{loc}}}{d\tau} = q \frac{d}{dq}f(q) = 1+\sum_{d=1}^{\infty} q^{d}\left\langle \mathrm{P}, \mathrm{P}^{n-1}, \widetilde{1}\right\rangle_{0,3, d}^{\left(\mathbf{e}, \, K_V^{-1}\right)} \, .
			\end{equation}
			Consequently, we have 
			\begin{equation}
				\label{eq-period_3-point_euler-twisted}
				q \frac{d}{dq}f(q) = \frac{1}{2\pi i} \frac{\omega_0(\tau)}{t(\tau)} \frac{dt}{d\tau}  \, ,
			\end{equation}
			where $\omega_0$, $t$ are the holomorphic period and complex structure modulus in the corresponding analytic neighborhood of the large complex structure limit of the mirror variety $X$.
			\item[b)] The regularized $\mathrm{SQDM}^{\mathrm{reg}}(V)$ (\ref{eq-regularized_SQDM}) is represented in the flat co-ordinate $\tau$ by the inhomogeneous Picard-Fuchs equation 
			\begin{equation}
				\label{eq-inhomogeneous_PF_tau}
				\mathcal{L}_{n}^{(\tau)} = \frac{d \tau_{\mathrm{loc}}}{d\tau} \, .
			\end{equation}
		\end{enumerate}
	\end{theorem}
	\begin{proof}
		\begin{enumerate}
			\begin{flushleft}
			\item[a)] Equation (\ref{eq-3-point_euler-twisted}) follows readily from Lemma \ref{lem-exp_mirror_maps}, by taking the logarithm of (\ref{eq-QSD_mirror}) and then differentiating with respect to $\tau$ using (\ref{eq-q_coordinate}) and then using the divisor equation, which is known to hold for twisted invariants, see Theorem \ref{thm-twisted_invariants}.
			\item[b)] As in Theorem \ref{thm-operator_factorization}, the regularized $\mathrm{SQDM}^{\mathrm{reg}}(V)$ is represented in the co-ordinate $t$ dual to the primitive generator $\mathrm{P}$ by the reducible, regular operator $\mathcal{D}_{n+1} = \theta \circ (t\mathcal{L}_n)$, and hence by the inhomogeneous Picard-Fuchs equation
			\begin{equation}
				\label{eq-inhomogeneous_PF_t}
				\mathcal{L}_n = \frac{1}{t} \, .
			\end{equation}
			The $\mathrm{SQDM}(V)$ encodes the 3-point invariants of $V$, and hence by writing (\ref{eq-inhomogeneous_PF_t}) in the flat co-ordinate $\tau$, we find that the operator becomes one of the form 
			\begin{equation}
				\label{eq-PF_tau}
				\mathcal{L}_{n}^{(\tau)} = \frac{1}{\mathcal{R}(\tau)} \mathcal{L}^{(\tau)}_{n-2} \circ \left(\frac{1}{\mathcal{Y}^{(3)}_{n-1}(\tau)}\frac{d^{2}}{d\tau^2}\right) \frac{1}{\omega_0(\tau)} \, ,
			\end{equation}
			where $\mathcal{Y}^{(3)}_{n-1}(\tau)$ is the 3-point generating function of the anticanonical divisor $X^\vee$ as in \S \ref{ss-I_functions}, $\mathcal{L}^{(\tau)}_{n-2}$ is a Fuchsian operator pulled back by the mirror map (defined only when $n \geq 3$, in analogy with (\ref{eq-CY3_PF_tau})) and $\mathcal{R}(\tau)$ is the reciprocal of the symbol of $\mathcal{L}_n$ transformed by the mirror map. Equating symbols, we find from Lemma \ref{lem-algebraic_yukawa}, Theorem \ref{thm-3-point_invariants}, and Equation (\ref{eq-normalized_Yukawa}) that
			\begin{equation}
				\label{eq-symbols}
				\mathcal{R}(\tau) = \frac{1}{t^{n} \mathcal{Q}(\tau) \mathcal{Y}_{n-1}(\tau) \omega_0(\tau)} \left(\frac{1}{2\pi i}\frac{dt}{d\tau}\right)^{n-1} \, .
			\end{equation}
			Again applying (\ref{eq-normalized_Yukawa}), we find that 
			\begin{equation}
				\label{eq-R(tau)}
				\mathcal{R}(\tau) = \frac{1}{2\pi i} \frac{\omega_0(\tau)}{t(\tau)} \frac{dt}{d\tau}  = 	\frac{d \tau_{\mathrm{loc}}}{d\tau} \, .
			\end{equation}
			\end{flushleft}
		\end{enumerate}
	\end{proof}
	\begin{remark}
		\label{rmk-gamma_conjecture}
		To the best of our knowledge, this is the first instance in which the association of the full regularized $\mathrm{SQDM}^{\mathrm{reg}}(V)$ has been explicitly linked with such euler-twisted invariants, though it implicitly appears by Iritani, Mann, \& Mignon \cite[\S 6.2, \S 6.4]{MR3519131}. For $n=3$, an inhomogeneous Picard-Fuchs equation equivalent to (\ref{eq-inhomogeneous_PF_t}) appears in the analysis of Golyshev \& Zagier in their proof of the gamma conjecture\footnote{This is called `gamma conjecture I' in \cite{MR3536989}.} for rank-1 Fano threefolds \cite[\S 2.1]{MR3462676}. In that work, the full MUM basis $\omega_0(t),\omega_1(t),\omega_{2}(t),\omega_3(t)$ for the equation (\ref{eq-inhomogeneous_PF_t}) at $t=0$ plays an important role in determining the so-called Frobenius limits of the rank-1 Fano threefold $V$. It is clear from Theorem \ref{thm-3-point_euler-twisted} that a similar analysis extends to $n \geq 4$-dimensional almost minimal Fano varieties, and that via a variation of parameters type expression for the solution $\omega_n(t)$ of (\ref{eq-inhomogeneous_PF_t}) that the 3-point euler-twisted invariants of $V$ play an important but subtle role in the story.
	\end{remark}
	\begin{remark}
		\label{rmk-modularity}
		The $\mathrm{D}2$ and $\mathrm{D}3$ operators that appear in connection with mirror symmetry of (almost minimal) del Pezzo surfaces and rank-1 Fano threefolds are by this point well-known to admit modular parameterizations, going back to Stienstra \cite{MR2282958} \& Zagier \cite{MR2500571} for the $\mathrm{D}2$ operators, and Golyshev \cite{MR2306141} for the $\mathrm{D}3$ operators. Roughly speaking, this means that the mirror map $t(\tau)$ forms a Hauptmodul for some appropriate group $\Gamma \subset \mathrm{SL}(2,\mathbb{R})$, and that the variation of Hodge structure (VHS) of the CY fibres of the LG model pulled back by the Hauptmodul is isogenous to a symmetric power of the standard VHS on the relative $H^1$ of the elliptic modular surface over the (compactification of the) modular curve $\Gamma \backslash \mathbb{H}$. See for example \cite[App. A]{ms2024} for a relevant summary of these ideas. In particular, the holomorphic period $\omega_0(\tau)$ is a modular form of weight-$(n-1)$ of Eisenstein type in each respective case $n=2,3$. It was shown by Doran \& Kerr \cite[Prop. 9.1]{MR2851155} that again when $n=2,3$, what we call $\mathcal{R}(\tau)$ in (\ref{eq-R(tau)}) is a $\mathbb{Q}$-Eisenstein series of weight-$(n+1)$, meaning $\mathcal{R}(\tau)$ is a modular form of weight-$(n+1)$ of Eisenstein type with constant terms that are rational at the cusp points. 
		\par Such modularity is deeply connected with the description of the complementary solution $\omega_n(t)$ in terms of a higher normal function \cite{MR0527828} and extension regulator classes \cite{MR2851155,MR3780269}, and thus an algebraic cycle on the CY fibre of the LG model. Eisenstein series appear since the relevant forms are the pullback of a motivic cohomology class on a Kuga variety to an LG model; such cycles are necessarily Eisenstein symbols \cite[\S 10]{MR2851155}. For LG models, a modular description is only possible for $n=2,3$ since the presence of both MUM and conifold points in the discriminant locus rules out $\mathrm{SL}(2)$ arising as a Mumford-Tate / geometric monodromy / differential Galois group more generally \cite[\S 10.6]{MR2851155}.
	\end{remark}
	\subsubsection{Modularity in low dimensions}
	\label{sss-modularity}
	\par Let us establish the following notation found in \cite[\S 4.1]{ms2024}.  All of the following data can be found in \cite{MR2282958} or \cite[App. B.3]{ms2024}. For the principal congruence subgroup $\Gamma_0(m) \subseteq \mathrm{SL}(2,\mathbb{Z})$, let $\mathcal{L}_{2,m}$ be the Picard-Fuchs operator annihilating the periods of the mirror family of elliptic curves for the anticanonical elliptic curves in the associated del Pezzo surface, and let $t(\tau)$ be the corresponding Hauptmodul of the modular curve $X_0(m) = \overline{\Gamma_0(m) \backslash \mathbb{H}}$. Then the associated weight-3 $\mathbb{Q}$-Eisenstein series $\mathcal{R}(\tau) \equiv \mathcal{R}_{3,m}(\tau)$ described above in Remark \ref{rmk-modularity} and determined by (\ref{eq-R(tau)}) has the following Lambert series expansion
	\begin{equation}
		\label{eq-weight-3_expansion}
		\mathcal{R}_{3,m}(\tau) = 1 + \sum_{k=1}^{\infty} b_k q^k = 1 + \sum_{k=1}^{\infty} k^2 \mathsf{N}_{m,k} \frac{q^k}{1-q^k} \, .
	\end{equation}
	The coefficients $b_k$ in the $q$-series on the left-hand side of (\ref{eq-weight-3_expansion}) can be seen to be integral, $b_k \in \mathbb{Z}$ for all $k \geq 1$, in the following way. The corresponding numbers $\mathsf{N}_{m,k} \in \mathbb{Z}$ are themselves integral, which follows by the $\mathbb{Q}$-Eisenstein condition. They are called (integral) Eisenstein weights, and appear in \cite[Tab. 1]{ms2024}. There are only finitely many such numbers, as they can be seen to exhibit the periodicity $\mathsf{N}_{m,k+m} = \mathsf{N}_{m,k}$ for all $k \in \mathbb{N}_{\geq 1}$.
	\par For concreteness, consider the degree-9 del Pezzo surface $S_9 = \mathbb{P}^2$. Then the mirror family of the anticanonical linear system $\mathcal{X}^\vee \to \mathbb{P}^1_\tau$ consisting of the usual cubic elliptic curves in $\mathbb{P}^2$ is the rational elliptic surface $\mathcal{X} \to \mathbb{P}^1_t \cong X_0(3)$ \cite{MR1732409}. The associated elliptic fibration is described by a Weierstrass model and associated Picard-Fuchs operator
	\begin{equation}
		\label{eq-RES_WM}
		X_t : \quad  y^2 = 4x^3 - g_2(t)x - g_3(t) \, , \quad
		\begin{aligned}
			g_2(t) &= 27 - 648 t \, , \\
			g_3(t) &= 5832 t^2 - 972 t + 27 \, , \\
			\mathcal{L}_{2,3} &= \theta^2 - 3t (3\theta + 1)(3\theta + 2) \, , 
		\end{aligned}
	\end{equation}
	with Hauptmodul $t(\tau)$, $q$-series of holomorphic period $\omega_0(t(\tau))$, and weight-3 $\mathbb{Q}$-Eisenstein series are given by
	\begin{equation}
		\label{eq-P2_modular}
		\begin{aligned}
			t(\tau) & = \frac{\eta(3 \tau)^{12}}{\eta(\tau)^{12}+27 \eta(3 \tau)^{12}} \\
			&= q-15 q^2+171 q^3-1679 q^4+15054 q^5-\cdots \\
			\omega_0(t(\tau)) &=1+6 \sum_{k=1}^\infty \chi_{-3}(k) \frac{q^k}{1-q^k} \\
			&= 1+6 q+6 q^3+6 q^4+12 q^7+\cdots \\
			\mathcal{R}_{3,3}(\tau) & = 1-9 \sum_{k=1}^{\infty}\chi_{-3}(k) \frac{ k^2q^k}{1-q^k} \, \\
			& = 1-9 q+27 q^2-9 q^3-117 q^4+216 q^5+\cdots
		\end{aligned}
	\end{equation}
	We use here the normalizations found in \cite[Tab. B.3]{ms2024}. Here $\eta(\tau) = q^{\frac{1}{24}} \prod_{n=1}^{\infty}(1-q^n)$ is the usual Dedekind $\eta$-function,  and $\chi_{-3}(k)$ is the quadratic Dirichlet character $\chi_{-3}(k) = 0,1,-1 \; \; \mathrm{if} \; \; k \equiv 0,1,2 \; \mathrm{mod} \; 3$. Accordingly, the numbers $\mathsf{N}_{3,k} = -9,9,0,\dots$ are the integral Eisenstein weights associated to $\mathbb{P}^2$. Moreover, the 3-point euler-twisted invariants of $\mathbb{P}^2$ are the coefficients in the $q$-expansion of $\mathcal{R}_{3,3}(\tau)$ in (\ref{eq-P2_modular}).
	\par Similarly, one may define the analogous data (third order Picard-Fuchs operator $\mathcal{L}_{3,N}$, Hauptmodul $t(\tau)$ of $X_0(N)^+$, weight-4 $\mathbb{Q}$-Eisenstein series $\mathcal{R}_{4,N}(\tau)$) associated to the mirror family $\mathcal{X} \to \mathbb{P}^1_t \cong X_0(N)^+$ for a rank-1 Fano threefold of type $(N,\delta)$. One has then
	\begin{equation}
		\label{eq-weight-4_expansion}
		\mathcal{R}_{4,N}(\tau) = 1 + \sum_{k=1}^{\infty} a_k q^k = 1 + \sum_{k=1}^{\infty} k^3 \mathsf{N}_{N,k} \frac{q^k}{1-q^k} 
	\end{equation}
	which gives the analogous Eisenstein weights $\mathsf{N}_{N,k}$. In this case, the weight-4 Eisenstein weights are generically rational, $\mathsf{N}_{N,k} \in \mathbb{Q}$, though again there are only finitely many such numbers and exhibit the periodicity $\mathsf{N}_{N,k+N} = \mathsf{N}_{N,k}$ for all $k \in \mathbb{N}_{\geq 1}$, both by $\mathbb{Q}$-Eisenstein. For each $N$, there is a single integer $j_N \in \mathbb{N}_{\geq 1}$ such that $j_N\mathsf{N}_{N,k} \in \mathbb{Z}$, and once again one finds that the coefficients $a_k$ in the left-hand side of (\ref{eq-weight-4_expansion}) are integral, $a_k \in \mathbb{Z}$ for all $k \geq 1$. 
	\par For the rank-1 Fano threefold $V = \mathbb{P}^3$ of type $(N,\delta) = (2,4)$ one has the following: 
	\begin{equation}
		\label{eq-P3_modular}
		\begin{aligned}
			\mathcal{L}_{3,2} &= \theta^3 - 8t (2\theta + 1)(4\theta + 1)(4\theta + 3) \, , \\
			t(\tau) &= \frac{\eta(\tau)^{24} \eta(2 \tau)^{24}}{\eta(\tau)^{48}+128 \eta(\tau)^{24} \eta(2 \tau)^{24}+4096 \eta(2 \tau)^{48}} \, , \\
			&= q-104 q^2+6444 q^3-311744 q^4+13018830 q^5-\cdots \\
			\omega_0(t(\tau)) &= 24 G_2(\tau)-2 \cdot 24 G_2(2 \tau) \\
			&= 1+24 q+24 q^2+96 q^3+24 q^4+\cdots \\
			\mathcal{R}_{4,2}(\tau) &=-80 G_4(\tau)+2^2 \cdot 80 G_4(2 \tau) \, , \\
			&= 1 - 80q - 400q^2 - 2240q^3 - 2960q^4 - 10080q^5 - \cdots
		\end{aligned}
	\end{equation}
	where $G_2(\tau), G_4(\tau)$ are the Hecke-normalized Eisenstein series of weight-2\footnote{$G_2(\tau)$ is only a quasi-modular form, but the combination appearing in (\ref{eq-P3_modular}) for $\omega_0(t(\tau))$ determines a modular form of weight-2 for $\Gamma_0(2)^+$.} and weight-4, respectively
	\begin{equation*}
	G_2(\tau) = -\frac{1}{24} + \sum_{k=1}^{\infty} k\frac{q^k}{1-q^k} \, ,	\quad \quad G_4(\tau) = \frac{1}{240} + \sum_{k=1}^{\infty} k^3\frac{q^k}{1-q^k} \, .
	\end{equation*}
	Accordingly, one has the integral Eisenstein weights $\mathsf{N}_{2,k}=-80,-40,\dots$; similarly to above, the 3-point euler-twisted invariants of $\mathbb{P}^3$ are the coefficients in the $q$-expansion of $\mathcal{R}_{4,2}(\tau)$ in (\ref{eq-P3_modular}).
	\par More generally, we obtain the following result for low dimensions. Set $n =2,3 = \dim(V)$.
	\begin{theorem}
		\label{thm-euler-twisted_modular}
		For $V$ an almost minimal del Pezzo surface or rank-1 Fano threefold, the generating function of the 3-point euler-twisted invariants (\ref{eq-3-point_euler-twisted}) obtained as the $q$-logarithmic derivative of (\ref{eq-twisted_gf}) is the $\mathbb{Q}$-Eisenstein series $\mathcal{R}(\tau)$ of weight-$(n+1)$ in each respective case $n=2,3$. In particular, the 3-point euler-twisted invariants
		\begin{equation*}
			\left\langle \mathrm{P}, \mathrm{P}^{n-1}, \widetilde{1}\right\rangle_{0,3, d}^{\left(\mathbf{e}, \, K_V^{-1}\right)}  \in \mathbb{Z}
		\end{equation*}
		are integral, and completely determined by the integral Eisenstein weights of \cite[\S 4.1]{ms2024}. Equivalently, the linear combinations of 3-point euler-twisted invariants given by the multicovering formula yield periodic integer sequences as described above.
	\end{theorem}
	Notice that the general integrality of such invariants for higher dimensional varieties, while reasonable to expect, is not obvious from the (meromorphic) factorization shown in (\ref{eq-period_3-point_euler-twisted}).
	\subsection{The reciprocal mirror map derivative \texorpdfstring{$\frac{d\tau}{d\tau_{\mathrm{loc}}}$ as generating function of 3-point local invariants}{d\tau/d_\tau_{\mathrm{loc}}}}
	\label{ss-reciprocal_mirror_map_derivative}
	\par Stienstra \cite[pp. 2]{MR2282958}, using the perspective of Mahler measures rather than local mirror maps, observed on an experimental basis that upon inverting the weight-3 $\mathbb{Q}$-Eisenstein series $\mathcal{R}_{3,m}(\tau)$ in (\ref{eq-weight-3_expansion}) associated to a del Pezzo surface $S$, a Lambert series $\mathsf{Q}$-expansion of 
	\begin{equation}
		\label{eq-reciprocal_Q-expansion}
		\frac{d\tau}{d\tau_{\mathrm{loc}}} = \mathsf{Q}\frac{d}{d\mathsf{Q}} \log \mathsf{q} =: 1 + \sum_{d=1}^{\infty} d^3 \widetilde{\mathsf{N}}_{m,d} \frac{\mathsf{Q}^d}{1-\mathsf{Q}^d}
	\end{equation}
	appeared to agree with the instanton expansion (\ref{eq-local_instanton_numbers}) computed by physicists \cite[Eq. (4.10)]{MR1468700} of the local del Pezzo, after normalization by an overall rational multiple $\mu \in \mathbb{Q}^*$.  Here we have $\mathsf{Q} = Q^\nu$ and $\mathsf{q} = q^\nu$, for some $\nu \in \mathbb{N}_{\geq 1}$. It implies that there is a functional relationship of the form
	\begin{equation}
		\label{eq-functional_rel}
		\beta^{-1}\mathsf{Q} = \mathsf{q} \prod_{d=1}^{\infty} (1- \mathsf{q}^d)^{-d \mathsf{N}_{m,d}} \, , \hspace{10mm} \beta \mathsf{q} = \mathsf{Q} \prod_{d=1}^{\infty} (1- \mathsf{Q}^d)^{-d^2 \widetilde{\mathsf{N}}_{m,d}} \, 
	\end{equation}
	for some undetermined $\beta \in \mathbb{Q}^*$. As before, these are reciprocal quantities. 
	\par Theorem \ref{thm-bkv} essentially proves Stienstra's observation. In \cite[\S 4]{ms2024}, it was observed by the second named author and A. Malmendier that since the weight-4 $\mathbb{Q}$-Eisenstein series $\mathcal{R}_{4,N}(\tau)$ associated to a rank-1 Fano threefold $V$ of type $(N,\delta)$ naturally admits a Lambert series of the form (\ref{eq-weight-4_expansion}), there is also naturally a functional relationship between the exponentiated closed and local mirror maps $q$ and $Q$ as in (\ref{eq-functional_rel}), by swapping the role of the numbers and weights. Here we have $\mathsf{Q} = Q^\delta$ and $\mathsf{q} = q^\delta$, again for some undetermined $\beta \in \mathbb{Q}^*$. Then we have in analogy to (\ref{eq-reciprocal_Q-expansion}):
	\begin{equation}
		\label{eq-local_fourfold_Q_expansion1}
		\frac{d\tau}{d\tau_{\mathrm{loc}}} = \mathsf{Q} \frac{d}{d\mathsf{Q}} \log \mathsf{q} =: 1 + \sum_{d=1}^{\infty} d^2 \widetilde{\mathsf{N}}_{N,d} \frac{\mathsf{Q}^d}{1-\mathsf{Q}^d}
	\end{equation}
	\begin{definition}
		\label{def-reciprocal_Q-Eisenstein}
		The a priori rational numbers $\widetilde{\mathsf{N}}_{m,d}, \widetilde{\mathsf{N}}_{N,d} \in \mathbb{Q}$ are called the degree-$d$ reciprocal $\mathbb{Q}$-Eisenstein numbers. The rational number $\beta \in \mathbb{Q}^*$ is called the $B$-field.
	\end{definition}
	\par Determining the $B$-field is important in physics, and can be quite subtle - for example, the difference between the local CY threefold $\mathrm{Tot}(K_{S_6})$ and the so-called conifold is respectively taking $\beta = \pm 1$, the remainder of the modular data is exactly the same. Again see \cite[Eq. (4.10)]{MR1468700}. In \cite[Def. 4.1]{ms2024}, the numbers in Definition \ref{def-reciprocal_Q-Eisenstein} were called local instanton numbers. We chose not to call them that, since we will prove in Theorem \ref{thm-local_3-point} below that they are local instanton numbers in the sense of our Definition \ref{def-local_yukawa}.
	\par This hinges on understanding the reciprocal mirror map derivative, in the sense of the inverse function theorem. This is perfectly natural, as we saw for local CY threefolds in \S \ref{sss-local_CY-3}. We show that it can be seen to generate local 3-point invariants, as in Theorem \ref{thm-bkv}.
	\begin{theorem}
		\label{thm-local_3-point}
		Let $V$ be an almost minimal Fano variety of dimension $n \geq 2$, and $Y^\vee_\tau = \mathrm{Tot}(K_V)$ be the associated local CY $(n+1)$-fold. Let $\tau,\tau_{\mathrm{loc}}$ be the closed and local mirror maps as in \S \ref{ss-closed_ms}, \ref{ss-local_ms}, convergent in an analytic neighborhood of the large volume limit $\tau = i \infty \in \mathbb{P}^1_\tau$ for $X^\vee_\tau \in |-K_V|$. Let $I^{\mathrm{loc}}(\tau_{\mathrm{loc}},z)$ be the local $I$-function for $Y^\vee_\tau$ as in \S \ref{ss-I_functions}. Then the local 3-point invariants of $Y^\vee_\tau$ are determined (up to normalization by a rational multiple $\mu \in \mathbb{Q}^*$) by the reciprocal mirror map derivative as
		\begin{equation}
			\label{eq-reciprocal_mirror_map}
			\mathcal{Y}^{(3)}_{\mathrm{P}^{n-1}}(Y^\vee_\tau) = \kappa^{\mathrm{loc}} \frac{d\tau}{d\tau_{\mathrm{loc}}} \, .
		\end{equation}
		In particular, the local 3-point invariants are obtained in this case, and the degree-$d$ local instanton numbers $N_d$ in (\ref{eq-local_instanton_numbers}) count genuine 3-point local invariants of $Y^\vee_\tau$, under the identification (\ref{eq-covering_formula}) for $m \geq 1$,
		\begin{equation}
			\label{eq-local_3-point_invariants}
			b_m = \kappa^{\mathrm{loc}}\langle \mathrm{P},\mathrm{P},\mathrm{P}^{n-1}\rangle_{0,3,m} \, .
		\end{equation}
	\end{theorem}
	\begin{proof}
		This is similar to the analysis in Theorem \ref{thm-3-point_invariants} but done in the local case. Then we have, after identifying the appropriate ratio of components of the local $I$-function \S \ref{ss-I_functions} with the local periods $\eta_k(\tau_{\mathrm{loc}})$ in (\ref{eq-local_periods}), we have from Lemma \ref{lem-3-point_I-function} that 
		\begin{equation}
			\label{eq-local_I-function_derivative}
			\mathcal{Y}^{(3)}_{\mathrm{P}^{n-1}}(Y^\vee_\tau) = \kappa^{\mathrm{loc}} \frac{d^2}{d\tau_{\mathrm{loc}}^2}\left(\frac{\eta_2(\tau_{\mathrm{loc}})}{\eta_0(\tau_{\mathrm{loc}})}\right) = \kappa^{\mathrm{loc}} \frac{d}{d\tau_{\mathrm{loc}}}\left(\frac{\omega_1(\tau_\mathrm{loc})}{\omega_0(\tau_{\mathrm{loc}})}\right) = \kappa^{\mathrm{loc}} \frac{d\tau}{d\tau_{\mathrm{loc}}} \, .
		\end{equation}
	\end{proof}
	This simple observation reproves, modulo details of normalization, Theorem \ref{thm-bkv}. It also allows us to prove a conjecture made by A. Malmendier and the second named author on genus zero invariants of local CY fourfolds. Since the proof is similar in spirit, we only prove the latter here.  We have the following.
	\begin{corollary}
		\label{cor-local_fourfold_instantons} \cite[Conjecture 4.3]{ms2024}
		Let $V$ be a smooth rank-1 Fano threefold of type $(N, \delta)$, as given in (\ref{eq-degree_level_invariants}). In the notation of Theorem \ref{thm-local_3-point},
		\begin{itemize}
			\item[a)] there are rational numbers $\beta,\mu \in \mathbb{Q}^*$ such that the normalized reciprocal $\mathbb{Q}$-Eisenstein numbers $\mu \widetilde{\mathsf{N}}_{N,d}$ expressed in $\mathsf{Q} = Q^{\delta}$ with choice of $B$-field determined by $\beta$ are the degree-$d$ genus zero invariants of the local Calabi-Yau fourfold $Y^\vee_\tau = \mathrm{Tot}(K_{V})$. 
			\item[b)] the $\mathsf{Q}$-expansions are expressible in terms of quasi-modular forms for some discrete subgroup of $\mathrm{SL}(2,\mathbb{R})$ commensurable with $\mathrm{SL}(2,\mathbb{Z})$.
		\end{itemize}
	\end{corollary}
	\begin{proof} This is mostly aligning the data from Theorem \ref{thm-3-point_invariants} with the details for mirror symmetry of rank-1 Fano threefolds discussed in \S \ref{ss-fano_ms}.
		\begin{itemize}
			\item[a)] Given $V$ of type $(N,\delta)$, the irreducible component of the regularized $\mathrm{SQDM}^{\mathrm{reg}}(V)$ is represented by the Picard-Fuchs operator $\mathcal{L}_{3,N}$ annihilating the periods of the $\mathrm{M}_N$-polarized pencil of Picard rank-19 K3 surfaces (\ref{eq-M_n-lattice}) pulled back by the morphism $t \mapsto t^\delta$ \cite[Thm. 5.1]{MR2306141}. These operators are listed in \cite[\S 5.8]{MR2306141}. It follows that the exponentiated closed mirror map also transforms as $q \mapsto q^\delta$. Then by Theorem \ref{thm-3-point_invariants}, the $\mathsf{Q} = Q^\delta$-Lambert series expansion (\ref{eq-local_fourfold_Q_expansion1}) of $\frac{d\tau}{d\tau_{\mathrm{loc}}}$ is, up to normalization by a rational multiple $\mu \in \mathbb{Q}^*$, the local 3-point generating function $\mathcal{Y}^{(3)}_{\mathrm{P}^{2}}(Y^\vee_\tau)$. The choice of $B$-field comes from the quantum Serre relation between the closed and local mirror maps (\ref{eq-QSD_mirror}).
			\item[b)] The reciprocal mirror map derivative
			\begin{equation*}
				\frac{d\tau_{\mathrm{loc}}}{d\tau} =  \mathcal{R}_{4,N}(\tau) 
			\end{equation*}
			is a $\mathbb{Q}$-Eisenstein series of weight-4 for the modular group $\Gamma_0(N)^+$, up to pullback by the morphism $q \mapsto q^\delta$. Hence the reciprocal quantity
			\begin{equation*}
				\frac{d\tau}{d\tau_{\mathrm{loc}}} = \frac{1}{\mathcal{R}_{4,N}(\tau_{\mathrm{loc}})}
			\end{equation*}
			is determined by quasi-modular forms for $\Gamma_0(N)^+$, in analogy to the case for local del Pezzo surfaces and the modular group $\Gamma_0(m)$ as discussed in \cite[\S 4.2]{MR3273318}.
		\end{itemize}
	\end{proof}
	\begin{remark}
		\label{rmk-twisted_invariants}
		Alternatively, one may see that the mirror map derivative $\frac{1}{2\pi i}\frac{d\tau}{d\tau_{\mathrm{loc}}}$ is a generating function of a priori rational numbers that are each highly nontrivial algebraic combinations of the 3-point euler-twisted invariants determined by the reciprocal $\frac{1}{2\pi i}\frac{d\tau_{\mathrm{loc}}}{d\tau}$ as in (\ref{eq-3-point_euler-twisted}), together with the 1-point invariants determined by the mirror map $\tau_{\mathrm{loc}}$. Such expressions may be viewed for leading order terms in Table \ref{tab-general_local_coeffs}. Using the machinery of GW-invariants, these numbers may be seen to themselves be 3-point invariants, and then traced to recover such quantities as genuine local 3-point invariants. In the sequel, we push in the other direction, and explicitly determine invariants of the anticanonical divisor in terms of twisted invariants.
	\end{remark}
	\section{Anticanonical Calabi-Yau varieties in \texorpdfstring{$\mathbb{P}^2 \, , \mathbb{P}^3,\mathbb{P}^4$}{P^2,P^3,P^4} and instanton counts}
	\label{s-projective}
	\par As a demonstration of the utility afforded from the perspective of QSD as established in \S \ref{s-twisted_invariants}, we explicitly calculate low degree euler-twisted invariants of (\ref{eq-3-point_euler-twisted}) for $\mathbb{P}^2,\mathbb{P}^3,$ and $\mathbb{P}^4$. This is an extension of calculations done by Iritani, Mann, \& Mignon \cite[\S 6.5]{MR3519131}. This allows us compute the Yukawa couplings for respective anticanonical CYs in terms of the euler-twisted invariants. We have the following generalities.
	\subsection{Yukawa couplings in terms of twisted invariants}
	\label{ss-twisted_yukawa}
	\par As usual let $V$ be an almost minimal Fano variety of dimension $n \geq 2$, etc. Theorem \ref{thm-3-point_euler-twisted} a) yields the simple observation that 
	\begin{equation}
		\label{eq-dt/dtau}
		\frac{dt}{d\tau} = 2\pi i \frac{t(\tau)}{\omega_0(\tau)} \frac{d\tau_{\mathrm{loc}}}{d\tau} \, .
	\end{equation}
	Hence the normalized Yukawa coupling $\mathcal{Y}_{n-1}(\tau)$ of the anticanonical CY (\ref{eq-normalized_Yukawa}) is given by 
	\begin{equation}
		\label{eq-twisted_yukawa}
		\mathcal{Y}_{n-1}(\tau) = \frac{1}{\omega_0(\tau)^{n+1}\mathcal{Q}(t(\tau))}\left(\frac{d\tau_{\mathrm{loc}}}{d\tau}\right)^{n-1} \, ,
	\end{equation}
	where $\mathcal{Q}(t)$ is obtained from the symbol of the Picard-Fuchs operator of the mirror CY as in Lemma \ref{lem-algebraic_yukawa}. Since (\ref{eq-3-point_euler-twisted}) is a generating function of 3-point euler-twisted invariants, we can explicitly determine the degree-$d$ instanton numbers (\ref{eq-instanton_numbers}) of the anticanonical divisor in terms of the twisted invariants and mirror data. Conversely, by reciprocity, using (\ref{eq-reciprocal_mirror_map}), we have that
	\begin{equation}
		\label{eq-yukawa_relations}
		\mathcal{Y}_{n-1}(\tau_{\mathrm{loc}}) = \frac{1}{\omega_0(\tau_{\mathrm{loc}})^{n+1}\mathcal{Q}(t(\tau_{\mathrm{loc}}))}\left(\frac{1}{2\pi i}\mathcal{Y}^{(3)}_{\mathrm{P}^{n-1}}(Y^\vee_\tau)\right)^{1-n} \, ,
	\end{equation}
	where $Y^\vee_{\tau} = \mathrm{Tot}(K_V)$ is the associated local CY.
	\par This explicit relationship between the Yukawa couplings is another expression of QSD. We demonstrate in the remainder of the article how known invariants for anticanonical CY varieties in $V = \mathbb{P}^n$, $n=2,3,4$ can be expressed in terms of euler-twisted invariants, and recovered from this perspective. In what follows, let us write 
	\begin{equation*}
		\begin{aligned}
			b_d &:= \left\langle \mathrm{P}, \mathrm{P}^{n-1}, \widetilde{1}\right\rangle_{0,3, d}^{\left(\mathbf{e}, \, K_V^{-1}\right)} \, ,\\
			\kappa^{-1}\mathcal{Y}_{n-1}(\tau) &= 1 + \sum_{m=1}^{\infty} a_m \mathsf{q}^m \, , \\
			(\kappa^{\mathrm{loc}})^{-1}\mathcal{Y}^{\mathrm{loc}}_{n+1}(\tau_{\mathrm{loc}}) &= 1 + \sum_{m=1}^{\infty} c_m \mathsf{Q}^m \, ,
		\end{aligned} 
	\end{equation*}
	and $\mathcal{X}^\vee \to \mathbb{P}^1_{\tau}$ be the pencil of anticanonical CY varieties over the complexified symplectic moduli space. Moreover, let us declare that local $\mathbb{P}^n$ is the local CY variety $Y^\vee = \mathrm{Tot}(K_{\mathbb{P}^n})$. We have the following, which holds more generally than just for $V = \mathbb{P}^n$:
	\begin{proposition}
		\label{prop-general_Q-expansion}
		For arbitrary choice of $B$-field $\beta \in \mathbb{Q}^*$, up to normalization by rational multiple $\mu \in \mathbb{Q}^*$, the local 3-point invariants $c_m$ in the local Yukawa coupling $\mathcal{Y}^{\mathrm{loc}}_{n+1}(\tau_{\mathrm{loc}})$ are expressible in terms of rational polynomial in the 3-point euler-twisted invariants $b_d$. Up to degree-6, they are given in Table \ref{tab-general_local_coeffs}.
	\begin{table}[ht]
		\centering
		$\begin{array}{|c|r|}
			\hline
			m & c_m \\
			\noalign{\hrule height 0.8pt}
			1 & -\frac{b_1}{\beta} \\
			2 & \frac{2 b_1^2-b_2}{\beta^2}\\
			3 & \frac{1}{2} \frac{-9 b_1^3+9 b_1 b_2-2 b_3}{\beta^3} \\
			4 & \frac{1}{3} \frac{32 b_1^4-48 b_1^2 b_2+16 b_1 b_3+6 b_2^2-3 b_4}{\beta^4} \\
			5 & \frac{1}{24} \frac{-625 b_1^5+1250 b_1^3 b_2-500 b_1^2 b_3+\left(-375 b_2^2+150 b_4\right) b_1+100 b_2 b_3-24 b_5}{\beta^5} \\
			6 & \frac{1}{10} \frac{648 b_1^6-1620 b_1^4 b_2+720 b_1^3 b_3+\left(810 b_2^2-270 b_4\right) b_1^2+\left(-360 b_2 b_3+72 b_5\right) b_1-45 b_2^3+45 b_2 b_4+20 b_3^2-10 b_6}{\beta^6} \\
			\hline
		\end{array}$
		\caption{Coefficients in $Q$-series expansion of $(\kappa^{\mathrm{loc}})^{-1}\mathcal{Y}^{\mathrm{loc}}_{n+1}(\tau_{\mathrm{loc}})$ in terms of the $b_d$.}
		\label{tab-general_local_coeffs}
	\end{table}
	\end{proposition}
	\subsubsection{\texorpdfstring{$\mathbb{P}^2$}{P^2}}
	\label{sss-P2}
	The 3-point euler-twisted invariants of $\mathbb{P}^2$ are determined by (\ref{eq-P2_modular}). We list out to degree $d=8$ in Table \ref{tab-P2}.\footnote{We truncated the tables at different degrees only to retain readability.}
	\begin{table}[ht]
		\centering
		$\begin{array}{|c|r|}
			\hline
			d & b_d \\
			\noalign{\hrule height 0.8pt}
			1 & -9 \\
			2 & 27\\
			3 & -9 \\
			4 & -117 \\
			5 & 216 \\
			6 & 27 \\
			7 & -450 \\
			8 & 459 \\
			\hline
		\end{array}$
		\caption{3-point $\left(\mathbf{e}, \, K_V^{-1}\right)$ euler-twisted invariants of $\mathbb{P}^2$. They are the coefficients of the $q$-series expansion of the weight-$3$ modular form $\mathcal{R}_{3,3}(\tau)$ of Eisenstein type given in (\ref{eq-P2_modular}).}
		\label{tab-P2}
	\end{table}
	\par We know that the normalized Yukawa coupling for the pencil $\mathcal{X}^\vee \to \mathbb{P}^1_{\tau}$ is identically constant as in Proposition \ref{prop-EC_K3_Yukawa}. Thus the expression (\ref{eq-twisted_yukawa}) must be constant. We will demonstrate for low degree terms that each combination of the twisted invariants sequentially vanish. One finds from the Picard-Fuchs operator $\mathcal{L}_{2,3}$ in (\ref{eq-RES_WM}) that the relevant symbol component is $\mathcal{Q}(t) = 1-27t$, and so expanding out (\ref{eq-twisted_yukawa}) term by term in the $b_d$ the following expressions in Table \ref{tab-P2_yukawa} below.
	\begin{table}[ht]
		\centering
		$\begin{array}{|c|r|}
			\hline
			m & a_m \\
			\noalign{\hrule height 0.8pt}
			1 & b_1+9 =0 \\
			2 & 9 b_1+b_2+54 =0\\
			3 & 54 b_1+9 b_2+b_3+252 =0 \\
			4 & 252 b_1+54 b_2+9 b_3+b_4+1008 =0 \\
			5 & 1008 b_1+252 b_2+54 b_3+9 b_4+b_5+3591 =0 \\
			6 & 3591 b_1+1008 b_2+252 b_3+54 b_4+9 b_5+b_6+11718 =0 \\
			7 & 11718 b_1+3591 b_2+1008 b_3+252 b_4+54 b_5+9 b_6+b_7+35604 =0 \\
			8 & 35604 b_1+11718 b_2+3591 b_3+1008 b_4+252 b_5+54 b_6+9 b_7+b_8+102006 =0 \\
			\hline
		\end{array}$
		\caption{Leading order coefficients in $q$-series expansion of $\mathcal{Y}_1(\tau)$ in terms of the $b_d$.}
		\label{tab-P2_yukawa}
	\end{table}
	\par It is elementary to check that each term in Table \ref{tab-P2_yukawa} vanishes sequentially upon substituting the values for $b_d$ found in Table \ref{tab-P2}, in perfect agreement with the discussions in \S \ref{ss-genus_zero_invariants}. Moreover, the genus zero invariants of local $\mathbb{P}^2$ can be computed as in Theorem \ref{thm-local_3-point}. It was shown in \cite[\S 4.2.1]{ms2024} that the genus zero GV invariants of local $\mathbb{P}^2$ can be computed from (\ref{eq-reciprocal_Q-expansion}), using the datum $\beta = -1, \mu = \frac{1}{3}$, and $\mathsf{Q} = Q^3$. These are given below in Table \ref{tab-local_P2} and agree with those computed by Coates \& Iritani \cite[Tab. 2]{MR4301054}. One finds that evaluating the expressions in Table \ref{tab-general_local_coeffs} with the $b_d$ in Table \ref{tab-P2}, the values $c_m$ determine the genus zero GV invariants in Table \ref{tab-local_P2} using the Lambert series for the local instanton expansion (\ref{eq-reciprocal_Q-expansion}). The second author believes that this can be used to understand from a more general perspective a result of Doran \& Kerr \cite[Cor. 5.3]{MR2851155} on the growth rate of local GW invariants. This will be the subject of future work.
	\begin{table}[ht]
		\centering
		$\begin{array}{|c|r|}
			\hline
			d & \frac{1}{3}\widetilde{\mathsf{N}}_{3,d} \\
			\noalign{\hrule height 0.8pt}
			1 & 3 \\
			2 & -6\\
			3 & 27 \\
			4 & -192 \\
			5 & 1695 \\
			6 & -17064 \\
			7 & 188454 \\
			8 & -2228160 \\
			\hline
		\end{array}$
		\caption{Genus zero GV invariants of local $\mathbb{P}^2$ computed the datum $\beta = -1, \mu = \frac{1}{3}$, and $\mathsf{Q} = Q^3$.}
		\label{tab-local_P2}
	\end{table}
	\subsubsection{\texorpdfstring{$\mathbb{P}^3$}{P^3}}
	\label{sss-P3}
	\par The analogous story for $\mathbb{P}^3$ and the anticanonical quartic K3 surfaces is similar in light of \S \ref{ss-EC_K3_mirror}, so we provide less exposition. The symbol component of the Picard-Fuchs operator $\mathcal{L}_{3,2}$ in \ref{eq-P3_modular} is $\mathcal{Q}(t) = 1-256t$. The genus zero invariants of local $\mathbb{P}^3$ were shown in \cite[\S 4.2.2]{ms2024} to be determined by the datum $\beta =1$, $\mu = -\frac{1}{4}$, and $\mathsf{Q} = Q^4$, in line with Corollary \ref{cor-local_fourfold_instantons} above. They agree with the invariants computed by Klemm \& Pandharipande \cite[Tab. 1]{MR2415462}. The 3-point euler-twisted invariants again conspire to sequentially force every term of order $m \geq 1$ in $\mathcal{Y}_2(\tau)$ to vanish identically, as required by Proposition \ref{prop-EC_K3_Yukawa}. The local invariants are recovered from the euler-twisted ones, as in \S \ref{sss-P2}. This is more nontrivial than even those results, as the expressions in Table \ref{tab-P3_yukawa} are quadratic in the $b_d$. This again should have explicit implications for the growth rates of local GW invariants.
	\begin{table}[ht]
		\centering
		$\begin{array}{|c|r|}
			\hline
			d & b_d \\
			\noalign{\hrule height 0.8pt}
			1 & -80 \\
			2 & -400\\
			3 & -2240 \\
			4 & -2960 \\
			5 & -10080 \\
			6 & -11200 \\
			7 & -27520 \\
			8 & -23440 \\
			\hline
		\end{array}$
		\caption{3-point $\left(\mathbf{e}, \, K_V^{-1}\right)$ euler-twisted invariants of $\mathbb{P}^3$. They are the coefficients of the $q$-series expansion of the weight-$4$ modular form $\mathcal{R}_{4,2}(\tau)$ of Eisenstein type given in (\ref{eq-P3_modular}).}
		\label{tab-P3}
	\end{table}
	\begin{table}[ht]
		\centering
		$\begin{array}{|c|r|}
			\hline
			m & a_m \\
			\noalign{\hrule height 0.8pt}
			1 & 2 b_1+160 =0\\
			2 & b_1^2+320 b_1+2 b_2+20000 =0\\
			3 & 160 b_1^2+2 b_1 b_2+40000 b_1+320 b_2+2 b_3+2244480 =0 \\
			4 & 20000 b_1^2+320 b_1 b_2+2 b_1 b_3+b_2^2+4488960 b_1+40000 b_2+320 b_3+2 b_4+237081120 =0\\
			\hline
		\end{array}$
		\caption{Leading order coefficients in $q$-series expansion of $\mathcal{Y}_2(\tau)$ in terms of the $b_d$.}
		\label{tab-P3_yukawa}
	\end{table}
	\begin{table}[ht]
		\centering
		$\begin{array}{|c|r|}
			\hline
			d & -\frac{1}{4}\widetilde{\mathsf{N}}_{4,d} \\
			\noalign{\hrule height 0.8pt}
			1 & -20 \\
			2 & -820\\
			3 & -68060 \\
			4 & -7486440 \\
			5 & -965038900 \\
			6 & -137569841980 \\
			7 & -21025364147340 \\
			8 & -3381701440136400 \\
			\hline
		\end{array}$
		\caption{Genus zero invariants of local $\mathbb{P}^3$ computed using the datum $\beta =1$, $\mu = -\frac{1}{4}$, and $\mathsf{Q} = Q^4$.}
		\label{tab-local_P3}
	\end{table}
	\subsubsection{\texorpdfstring{$\mathbb{P}^4$}{P^4}}
	\label{sss-P4}
	\par This is the realm of classical mirror symmetry a la Candelas et al \cite{MR1101784}, for $\mathbb{P}^4$ and the anticanonical quintic Calabi-Yau threefolds. Let $X^\vee_\tau \subset \mathbb{P}^4$ be a very general smooth quintic threefold. Then $X^\vee_\tau$ has nontrivial 3-point invariants, which are well-known and can be computed using (\ref{eq-instanton_numbers}) from the Picard-Fuchs data in (\ref{eq-quintic_PF}). We list the leading order expansion of $\mathcal{Y}_3(X^\vee_\tau)$ and instanton numbers in Table \ref{tab-quintic_instantons}. The same analysis works, and allows one to simultaneously express the genus zero GV invariants and local 3-point invariants in terms of the 3-point euler-twisted invariants. The Picard-Fuchs operator $\mathcal{L}_4$, mirror map $t(\tau)$ $q$-series, and holomorphic period $\omega_0(\tau)$ are given below: 
	\begin{equation}
		\label{eq-quintic_PF}
		\begin{aligned}
			\mathcal{L}_4 &= \theta^4 -5t(5\theta+1)(5\theta+2)(5\theta+3)(5\theta+4) \, , \\
			t(\tau) &= q - 770q^2 + 171525q^3 -81623000q^4 -35423171250q^5 \cdots \\
			\omega_0(\tau) &= 1 + 120q + 21000 q^2 +14115000 q^3 + 13414125000 q^4+ \cdots
		\end{aligned}
	\end{equation}
	with symbol component $\mathcal{Q}(t) = 1-3125t$. To the best of our knowledge, the 3-point local invariants of local $\mathbb{P}^4$ displayed in Table \ref{tab-local_P4} have not appeared in the literature before. In analogy with above, we take $\beta = -1$, $\mu = \frac{1}{5}$, and $\mathsf{Q} = Q^5$. The 3-point euler-twisted invariants were essentially computed in \cite[\S 6.5]{MR3519131}, in the form of the relevant 2-point euler-twisted invariants (\ref{eq-twisted_gf}). This again should have implications for the growth rates of local GW invariants.
	\begin{table}[ht]
		\centering
		$\begin{array}{|c|r|r|}
			\hline
			m,d & 5\cdot a_m & N_d \\
			\noalign{\hrule height 0.8pt}
			1 & N_1 = 2875 & 2875 \\
			2 & N_1 + 8N_2 = 4876875 & 609250 \\
			3 & N_1 + 27N_3 = 8564575000 & 317206375 \\
			4 & N_1 + 8N_2 + 64N_4 = 15517926796875 & 242467530000\\
			\hline
		\end{array}$
		\caption{Leading order coefficients in $q$-series expansion of $\mathcal{Y}_3(\tau)$ in terms of the the instanton numbers $N_d$. As in \ref{sss-projective_CY-threefolds}, we normalize with the triple intersection number $\kappa =5$.}
		\label{tab-quintic_instantons}
	\end{table}
	\begin{table}[ht]
		\centering
		$\begin{array}{|c|r|}
			\hline
			d & b_d \\
			\noalign{\hrule height 0.8pt}
			1 & -650\\
			2 & -321250\\
			3 & -337216250 \\
			4 & -435997681250 \\
			5 & -626257527494400 \\
			6 & -958598821553843650 \\
			7 & -1531227197616745455000 \\
			8 & -2521899259208568536561250 \\
			\hline
		\end{array}$
		\caption{3-point $\left(\mathbf{e}, \, K_V^{-1}\right)$ euler-twisted invariants of $\mathbb{P}^4$. These invariants also appear to be integral $b_d \in \mathbb{Z}$ for all $d \geq 1$, but we do not have proof.}
		\label{tab-P4}
	\end{table}
	\begin{table}[ht]
		\centering
		$\begin{array}{|c|r|}
			\hline
			m & 5\cdot a_m \\
			\noalign{\hrule height 0.8pt}
			1 & 5(3 b_1+2525) = 2875 \\
			2 & 5(3 b_1^2+7575 b_1+3 b_2+5595375) = 4876875\\
			3 & 5(b_1^3+7575 b_1^2+6 b_1 b_2+16786125 b_1+7575 b_2+3 b_3+11890326250) = 8564575000\\
			4 & 5(2525 b_1^3+3 b_1^2 b_2+16786125 b_1^2+15150 b_1 b_2+6 b_1 b_3+3 b_2^2+35670978750 b_1 \\ 
			&+16786125 b_2+7575 b_3+3 b_4 + 24764887590625) = 15517926796875\\
			\hline
		\end{array}$
		\caption{Leading order coefficients in $q$-series expansion of $\mathcal{Y}_3(\tau)$ in terms of the $b_d$. As in \S \ref{sss-projective_CY-threefolds}, we normalize with the triple intersection number $\kappa =5$. By equating terms in Table \ref{tab-quintic_instantons}, the $\left(\mathbf{e}, \, K_V^{-1}\right)$ euler-twisted invariants are explicitly related to the genus zero GV invariants of the quintic.}
		\label{tab-P4_yukawa}
	\end{table}
	\begin{table}[ht]
		\centering
		$\begin{array}{|c|r|}
			\hline
			d & \frac{1}{5}\widetilde{\mathsf{N}}_d \\
			\noalign{\hrule height 0.8pt}
			1 & -130 \\
			2 & 116690\\
			3 & -167512290 \\
			4 & 294356632000 \\
			5 & -579271009849750 \\
			6 & 1226054132473883250 \\
			7 & -2730358342178513883410 \\
			8 & 6313088053137552766432000 \\
			\hline
		\end{array}$
		\caption{Genus zero invariants of local $\mathbb{P}^4$ using (\ref{eq-local_instanton_numbers}) computed using the datum $\beta = -1$, $\mu = \frac{1}{5}$, and $\mathsf{Q} = Q^5$. These invariants again appear to be integral, but we do not have proof.}
		\label{tab-local_P4}
	\end{table}
	\newpage 
	\bibliographystyle{abbrv}
	\bibliography{QSD_mirror}
	\printaddresses 
	\appendix
	\section{Twisted invariants \& quantum Serre duality}
	\label{app-quantum_serre}
	\subsection{Gromov-Witten Invariants}
	\label{ss-GW_invariants}
	Let $V$ be a smooth projective complex variety of dimension $r$,
	and let $\operatorname{Eff}(V)\subset H_2(V;\mathbb Z)$ denote
	the effective curve classes, including $0$.
	For $g,n\in\mathbb Z_{\geq0}$ and $d\in\operatorname{Eff}(V)$,
	let $\overline{M}_{g,n,d}(V)$ denote the moduli stack
	of stable maps
	$$
	f:(C,p_1,\ldots,p_n)\longrightarrow V,
	\qquad f_*[C]=d,
	$$
	where $C$ is a connected projective nodal curve of arithmetic
	genus $g$, the marked points are distinct and smooth, and
	stability means that the automorphism group of the marked map
	is finite.
	
	The virtual fundamental class
	$[\overline{M}_{g,n,d}(V)]^{\mathrm{vir}}$
	has complex dimension 
	\begin{equation*}
	D := (1-g)(r-3)+\int_d c_1(TV)+n \, .
	\end{equation*}
	For each marking, define the $i$th evaluation map $\operatorname{ev}_i$ and $i$th cotangent class $\psi_i$ by
	$$
	\operatorname{ev}_i([C,p_1,\ldots,p_n,f])=f(p_i),
	\qquad
	\psi_i=c_1(\mathcal L_i),
	$$
	where $\mathcal L_i$ is the line bundle with fibre $T^*_{p_i}C$.
	
	For $a_i\in H^*(V;\mathbb C)$ and $k_i\in\mathbb Z_{\geq0}$,
	the descendant Gromov--Witten invariant is
	$$
	\left\langle
	a_1\psi^{k_1},\ldots,a_n\psi^{k_n}
	\right\rangle_{g,n,d}
	:=
	\int_{[\overline{M}_{g,n,d}(V)]^{\mathrm{vir}}}
	\prod_{i=1}^{n}
	\left(\operatorname{ev}_i^*a_i\cup\psi_i^{k_i}\right).
	$$
	These are called primary invariants when every $k_i=0$.
	For homogeneous insertions, the invariant vanishes unless
	$$
	\sum_{i=1}^{n}\bigl(\deg a_i+2k_i\bigr)=2D,
	$$
	where $\deg$ denotes cohomological degree.
	
	\paragraph{Pairing and generating series.}
	For the quantum product, work with the even cohomology $H:=H^{\mathrm{ev}}(V;\mathbb C),$ equipped with the Poincar\'e pairing$(a,b):=\int_V a\cup b.$
	
	Choose a homogeneous basis $T_0=1,\ldots,T_N$ and its
	Poincar\'e dual basis $T^0,\ldots,T^N$, so that
	$(T_\alpha,T^\beta)=\delta_\alpha^\beta$. Write  $$ \tau=\sum_{\alpha=0}^N t^\alpha T_\alpha,\qquad
	\partial_\alpha=\frac{\partial}{\partial t^\alpha}.$$
	
	The Novikov ring $\mathbb C[[Q]]$ consists of formal sums
	$\sum_{d\in\operatorname{Eff}(V)}c_dQ^d$, with
	$Q^dQ^{d'}=Q^{d+d'}$, subject to the following finiteness condition:
	for a fixed ample divisor class $A$ and each real number $R$,
	only finitely many nonzero terms have $\int_d A\leq R$.
	We abbreviate $\mathbb C[[Q,\tau]]:=
	\mathbb C[[Q]][[t^0,\ldots,t^N]].$
	
	For primary or descendant insertions, define
	$$
	\begin{aligned}
		&\left\langle\!\left\langle
		a_1\psi^{k_1},\ldots,a_n\psi^{k_n}
		\right\rangle\!\right\rangle_\tau
		\\
		&\quad:=
		\sum_{d\in\operatorname{Eff}(V)}\sum_{m\geq0}
		\frac{Q^d}{m!}
		\left\langle
		a_1\psi^{k_1},\ldots,a_n\psi^{k_n},
		\underbrace{\tau,\ldots,\tau}_{m}
		\right\rangle_{0,n+m,d}.
	\end{aligned}
	$$
	Terms with $d=0$ and $n+m<3$ are omitted.
	The notation extends linearly to formal series of insertions.
	
	The quantum product $\bullet_\tau$ is defined by
	\begin{equation}
		\label{eq-quantum_product}
		(a\bullet_\tau b,c)
		=
		\left\langle\!\left\langle a,b,c
		\right\rangle\!\right\rangle_\tau,
	\end{equation}
	or equivalently,
	$$
	a\bullet_\tau b
	=
	\sum_{\alpha=0}^N
	\left\langle\!\left\langle a,b,T_\alpha
	\right\rangle\!\right\rangle_\tau T^\alpha.
	$$
	It extends $\mathbb C[[Q,\tau]]$ bilinearly to
	$H\otimes_{\mathbb C}\mathbb C[[Q,\tau]]$.
	The product is associative and commutative, with unit $1$, and
	$$
	(a\bullet_\tau b,c)=(a,b\bullet_\tau c),
	\qquad
	\left.a\bullet_\tau b\right|_{Q=0}=a\cup b.
	$$
	
	For a formal parameter $z$, the quantum $D$-module is
	$$
	\operatorname{QDM}(V):=(\mathcal F,\nabla,S),
	\qquad
	\mathcal F
	:=
	H\otimes_{\mathbb C}\mathbb C[z][[Q,\tau]],
	$$
	where the quantum connection and pairing are
	$$
	\nabla_\alpha
	:=
	\partial_\alpha+\frac{1}{z}(T_\alpha\bullet_\tau),
	\qquad
	S(u,v):=\bigl(u(-z),v(z)\bigr).
	$$
	Here the connection differentiates in $\tau$, with
	$\nabla_\alpha\mathcal F\subset z^{-1}\mathcal F$.
	The pairing is nondegenerate and satisfies
	$$
	S(fu,hv)=f(-z)h(z)S(u,v)
	$$
	for scalar series $f,h$.
	The connection is flat and preserves the pairing:
	$$
	[\nabla_\alpha,\nabla_\beta]=0,
	\qquad
	\partial_\alpha S(u,v)
	=
	S(\nabla_\alpha u,v)+S(u,\nabla_\alpha v).
	$$
	
	Define
	\begin{equation}
		\label{eq-fundamental_solution}
		L(\tau,z)a
		:=
		a-\sum_{\alpha=0}^N
		\left\langle\!\left\langle
		\frac{a}{z+\psi},T_\alpha
		\right\rangle\!\right\rangle_\tau T^\alpha,
		\qquad
		\frac{1}{z+\psi}
		=
		\sum_{k\geq0}(-1)^k\psi^kz^{-k-1}.
	\end{equation}
	Then $L(\tau,z)=\mathrm{id}+O(z^{-1})$, and for constant
	$a,b\in H$,
	$$
	\nabla_\alpha\bigl(L(\tau,z)a\bigr)=0,
	\qquad
	\bigl(L(\tau,-z)a,L(\tau,z)b\bigr)=(a,b).
	$$
	Thus $L$ is a fundamental solution in the $\tau$ directions.
	The $J$-function is
	$$
	J(\tau,z):=zL(\tau,z)^{-1}1
	=z1+\tau+O(z^{-1}),
	\qquad
	\partial_\alpha J(\tau,z)=L(\tau,z)^{-1}T_\alpha.
	$$
	
	Define endomorphisms $\delta,\rho$ of $H$ by
	$$
	\delta(a):=\frac{\deg a}{2}a
	\quad\text{for homogeneous }a,
	\qquad
	\rho(a):=c_1(TV)\cup a.
	$$
	Under the identification $T_\tau H=H$, the Euler vector field is
	$$
	\mathcal E(\tau)
	:=
	c_1(TV)+
	\sum_{\alpha=0}^N
	\left(1-\frac{\deg T_\alpha}{2}\right)t^\alpha T_\alpha.
	$$
	The quantum connection extends to a flat connection in the $z$-direction by
	$$
	\nabla_{z\partial_z}
	:=
	z\partial_z-\frac{1}{z}(\mathcal E(\tau)\bullet_\tau)+\delta,
	\qquad
	[\nabla_\alpha,\nabla_{z\partial_z}]=0.
	$$
	A fundamental solution for the extended connection is
	$$
	L(\tau,z)z^{-\delta}z^\rho,
	$$
	where $z^B:=\exp((\log z)B)$ for an endomorphism $B$,
	after adjoining $\log z$ and inverting $z$.
	
	With this grading convention, the pairing satisfies
	$$
	S(\nabla_{z\partial_z}u,v)+S(u,\nabla_{z\partial_z}v)
	=
	z\partial_zS(u,v)+rS(u,v).
	$$
	
	\par The small quantum $D$-module $\mathrm{SQDM}(V)$ is the datum above, restricted to the divisorial directions $\tau \in H^2(V,\mathbb{C})$.
	
	\subsection{Twisted GW invariants}
	\label{ss-twisted_invariants}
Twisted Gromov-Witten invariants, introduced in general by Coates and Givental in \cite{MR2276766}, are modifications of Gromov-Witten invariants made by integrating against a modified virtual fundamental cycle. 

For this work, we limit ourselves to considering twistings of type I in Coates-Givental's taxonomy. Let $\overline{M}_{g,n,d}(V)$ denote Kontsevich's moduli space of degree-$d$ stable maps from $n$-pointed genus-$g$ curves to $V$. These spaces have $n$ evaluation maps $\rm{ev}_i$ to $V$, representing evaluation of the map at the $i$th marked point. The universal family over $\overline{M}_{g,n,d}(V)$ is $\overline{M}_{g,n+1,d}(V)$, with the map $\pi$ defined by forgetting the last marked point. 

As such, the universal stable map to $V$ is described by the following diagram:
\[
\begin{tikzcd}
	\overline{M}_{g,n+1,d}(V)
	\arrow[r, "\mathrm{ev}_{n+1}"]
	\arrow[d, "\pi"']
	& V \\
	\overline{M}_{g,n,d}(V)
	&
\end{tikzcd}
\]

Let $\mathbf{c}$ be an invertible multiplicative characteristic class, and let $V$ be a vector bundle (or $K$-theory class) on the target space $V$. $(\mathbf{c},V)$-twisted invariants are defined by integrating with respect to the class:

$$[\overline{M}_{g,n,d}(V)]^{vir}\cap \mathbf{c}(\pi_*\rm{ev}_{n+1}^*V)$$
In the same way that an input in the $i$th seat of a Gromov-Witten invariant corresponds to imposing constraints on the stable map at the $i$th marked point, a type-I twisting can be thought of as imposing a global constraint on the stable map itself. 

\begin{theorem}
	\label{thm-twisted_invariants}
	Type I-twisted Gromov-Witten invariants satisfy the Gromov-Witten axioms, provided the Poincare pairing on $V$ is changed to the twisted pairing, defined by $(a,b)^{tw}:=\int_{V} ab \mathbf{c}(V)$. In particular, the genus-0 twisted invariants have the structure of a Frobenius manifold. 
\end{theorem}

As a consequence, twisted invariants have a quantum product, and a resulting quantum $D$-module: write $Q$ for the Novikov variable and $\mathbf s$ for the
parameters of $\mathbf c$. For primary or descendant insertions, set
$$
\left\langle\!\left\langle
a_1,\ldots,a_\ell
\right\rangle\!\right\rangle_\tau
:=
\sum_{d\in\operatorname{Eff}(V)}\sum_{m\geq0}
\frac{Q^d}{m!}
\left\langle
a_1,\ldots,a_\ell,
\underbrace{\tau,\ldots,\tau}_{m}
\right\rangle^{(\mathbf c,E)}_{0,\ell+m,d},
$$
where $\operatorname{Eff}(V)$ denotes the effective curve classes
and unstable terms are omitted.

The twisted pairing $g$ and quantum product $\bullet_\tau$
are defined by
$$
g(a,b):=\int_V a\cup b\cup\mathbf c(E),
\qquad
g(a\bullet_\tau b,c)
:=
\left\langle\!\left\langle a,b,c
\right\rangle\!\right\rangle_\tau.
$$
The product is associative and commutative, with unit $1$, and
$$
g(a\bullet_\tau b,c)=g(a,b\bullet_\tau c).
$$

The twisted quantum $D$-module is
$$
\operatorname{QDM}_{(\mathbf c,E)}(V)
:=(\mathcal F,\nabla,S),
\qquad
\mathcal F:=H\otimes_{\mathbb C}
\mathbb C[z][[Q,\mathbf s,\tau]],
$$
where
$$
\nabla_\alpha
:=\partial_\alpha+\frac{1}{z}(T_\alpha\bullet_\tau),
\qquad
S(u,v):=g\bigl(u(-z),v(z)\bigr).
$$
All series are formal. The connection differentiates only in
$\tau$, with
$\nabla_\alpha\mathcal F\subset z^{-1}\mathcal F$.
The pairing is nondegenerate and, for scalar series $f,h$, satisfies
$$
S(fu,hv)=f(-z)h(z)S(u,v).
$$
The connection and pairing satisfy
$$
[\nabla_\alpha,\nabla_\beta]=0,
\qquad
\partial_\alpha S(u,v)
=
S(\nabla_\alpha u,v)+S(u,\nabla_\alpha v).
$$

Define
$$
L(\tau,z)a
:=
a-\sum_{\alpha=0}^N
\left\langle\!\left\langle
\frac{a}{z+\psi},T_\alpha
\right\rangle\!\right\rangle_\tau
\bigl(\mathbf c(E)^{-1}\cup T^\alpha\bigr),
\qquad
\frac{1}{z+\psi}
=
\sum_{k\geq0}(-1)^k\psi^k z^{-k-1}.
$$
Then $L(\tau,z)=\mathrm{id}+O(z^{-1})$, and for constant
$a,b\in H$,
$$
\nabla_\alpha\bigl(L(\tau,z)a\bigr)=0,
\qquad
g\bigl(L(\tau,-z)a,L(\tau,z)b\bigr)=g(a,b).
$$

\begin{remark}
	For nontrivial twists, there is no natural extension of the connection into the $z$-direction. 
\end{remark}

\subsubsection{The twisted $J$-function}
Define the twisted $J$-function by:
$$
J_{(\mathbf c,E)}(\tau,z)
:=
z1+\tau+
\sum_{\alpha=0}^N
\left\langle\!\left\langle
\frac{T_\alpha}{z-\psi}
\right\rangle\!\right\rangle_\tau^{(\mathbf c,E)}
\bigl(\mathbf c(E)^{-1}\cup T^\alpha\bigr),
$$
where the fraction denotes the descendant expansion
$$
\frac{1}{z-\psi}
=
\sum_{k\geq0}\psi^kz^{-k-1}.
$$
Writing $L_{(\mathbf c,E)}$ for the twisted fundamental solution,
we have
$$
J_{(\mathbf c,E)}(\tau,z)
=
zL_{(\mathbf c,E)}(\tau,z)^{-1}1,
\qquad
\partial_\alpha J_{(\mathbf c,E)}(\tau,z)
=
L_{(\mathbf c,E)}(\tau,z)^{-1}T_\alpha.
$$
In particular,
$J_{(\mathbf c,E)}(\tau,z)=z1+\tau+O(z^{-1})$.

\subsection{The Lagrangian Cone \& Quantum Serre duality}
\label{ss-lagrangian_cone}
Set
$$
\mathcal H
:=
H\otimes_{\mathbb C}\mathbb C[z,z^{-1}][[Q,\mathbf s]],
\qquad
\Omega_{(\mathbf c,E)}(f,h)
:=
\operatorname{Res}_{z=0}
\bigl(f(-z),h(z)\bigr)_{(\mathbf c,E)}\,dz.
$$
The standard polarization is the decomposition
$$
\mathcal H=\mathcal H_+\oplus\mathcal H_-,
\qquad
\begin{aligned}
	\mathcal H_+
	&:=H\otimes_{\mathbb C}\mathbb C[z][[Q,\mathbf s]],\\
	\mathcal H_-
	&:=H\otimes_{\mathbb C}
	z^{-1}\mathbb C[z^{-1}][[Q,\mathbf s]].
\end{aligned}
$$
Both summands are Lagrangian for $\Omega_{(\mathbf c,E)}$.
This decomposition identifies $\mathcal H$ with
$T^*\mathcal H_+$: an element $p\in\mathcal H_-$ represents
the covector
$v\mapsto\Omega_{(\mathbf c,E)}(p,v)$ on $\mathcal H_+$.

Introduce a formal element
$$
\mathbf t(z)=\sum_{k\geq0}t_kz^k,
\qquad
t_k=\sum_{\alpha=0}^N t_k^\alpha T_\alpha.
$$
The change of variables
$$
\mathbf q(z):=\mathbf t(z)-z1
$$
is called the dilaton shift.
Define the twisted genus zero descendant potential, regarded
as a function of $\mathbf q$ near $-z1$, by
$$
\mathcal F_V^{0,(\mathbf c,E)}(\mathbf q)
:=
\sum_{d\in\operatorname{Eff}(V)}\sum_{n\geq0}
\frac{Q^d}{n!}
\left\langle
\mathbf t(\psi_1),\ldots,\mathbf t(\psi_n)
\right\rangle_{0,n,d}^{(\mathbf c,E)},
\qquad
\mathbf t(\psi_i):=\sum_{k\geq0}t_k\psi_i^k.
$$
All expressions are interpreted formally in $Q$, $\mathbf s$,
and the coordinates $t_k^\alpha$.
In every sum, invariants with $d=0$ and fewer than three
marked points are omitted.

Under the identification $\mathcal H=T^*\mathcal H_+$
defined by $\Omega_{(\mathbf c,E)}$, set
$$
\mathcal L_V^{(\mathbf c,E)}
:=
\left\{
\mathbf q+d_{\mathbf q}\mathcal F_V^{0,(\mathbf c,E)}
\right\}
\subset\mathcal H,
$$
where $\mathbf q$ varies formally near $-z1$.
Thus $\mathcal L_V^{(\mathbf c,E)}$ is the graph of the
differential of the shifted twisted descendant potential.

Equivalently, its points have the form
$$
\begin{aligned}
	f_{\mathbf t}^{(\mathbf c,E)}(z)
	={}&-z1+\mathbf t(z)\\
	&+
	\sum_{\alpha=0}^N
	\sum_{d\in\operatorname{Eff}(V)}\sum_{n\geq0}
	\frac{Q^d}{n!}
	\left\langle
	\frac{T_\alpha}{-z-\psi},
	\underbrace{\mathbf t(\psi),\ldots,\mathbf t(\psi)}_{n}
	\right\rangle_{0,n+1,d}^{(\mathbf c,E)}
	\bigl(\mathbf c(E)^{-1}\cup T^\alpha\bigr).
\end{aligned}
$$
Here each occurrence of $\psi$ refers to the cotangent class
at the corresponding marking.

The submanifold $\mathcal L_V^{(\mathbf c,E)}$ is Lagrangian
for $\Omega_{(\mathbf c,E)}$ and is a cone with vertex at the origin.
For every tangent space $T$ to $\mathcal L_V^{(\mathbf c,E)}$,
$$
zT=T\cap\mathcal L_V^{(\mathbf c,E)},
$$
and the tangent space is constant along this intersection.
These statements are understood in the sense of formal germs.

Setting $\mathbf t(z)=\tau$ gives $J_{(\mathbf c,E)}(\tau,-z)$.
It is the unique point of $\mathcal L_V^{(\mathbf c,E)}$ whose
$\mathcal H_+$ component is $-z1+\tau$:
$$
\mathcal L_V^{(\mathbf c,E)}
\cap\bigl(-z1+\tau+\mathcal H_-\bigr)
=
\{J_{(\mathbf c,E)}(\tau,-z)\}.
$$
Its tangent space is
$$
\begin{aligned}
	T_{J_{(\mathbf c,E)}(\tau,-z)}
	\mathcal L_V^{(\mathbf c,E)}
	&=
	\bigoplus_{\alpha=0}^N
	\mathbb C[z][[Q,\mathbf s,\tau]]\,
	\partial_\alpha J_{(\mathbf c,E)}(\tau,-z)\\
	&=
	L_{(\mathbf c,E)}(\tau,-z)^{-1}\mathcal F_{(\mathbf c,E)},
\end{aligned}
$$
where
$$
\mathcal F_{(\mathbf c,E)}
=
H\otimes_{\mathbb C}\mathbb C[z][[Q,\mathbf s,\tau]]
$$
is the underlying module of $\operatorname{QDM}_{(\mathbf c,E)}(V)$.

The quantum Serre duality of Coates-Givental is the following statement:

$$c(E)\mathcal{L}_{(\mathbf c,E)}=\mathcal{L}_{(\mathbf c*,E^*)}.$$

The work of Iritani-Mann-Mignon \cite{MR3519131} translates the above statement into the following statement regarding quantum $D$-modules. 

This works in the following way: 

Define the map  
\begin{equation}
	\label{eq-equivariant_map}
	f_\circ(\tau):=\sum_{\alpha=0}^s\sum_{n,d} \frac{Q^d}{n!}\langle\langle T^\alpha, \mathbf{c}^*(E^*)\rangle \rangle_{\tau}^{(\mathbf{c},E)} T_\alpha
\end{equation}

Define the quantum Serre pairing $QS$ by $QS(a(z),b(z)):=\int_{V} a(z)b(-z)$. 

Theorem 2.11 of \cite{MR3519131} states the following:

\begin{theorem}\hfill
	\begin{itemize}
		\item[(1)] $QDM_{(\mathbf{c},E)}$ and $f_\circ^*QDM_{(\mathbf{c}^*,E^*)}$ are dual under $QS$
		\item[(2)] Multiplication by $\mathbf{c}(E)$ intertwines the quantum connections and quantum Serre pairings on the above two $D$-modules. 
		%\item There is a third point about fundamental solutions but unclear if we need it?
	\end{itemize}
\end{theorem}

\subsection{Specialization to Euler-twisted theory}
\label{ss-euler_twisted}

Let $\mathbb{C}^*$ act trivially on $X$, let $E$ now be a vector bundle on $X$ such that $\mathbb{C}^*$ acts on $E$ by scaling the fibers. Let $\mathbf{e}_\lambda(E)$ denote the $\mathbb{C}^*$-equivariant Euler class on $E$, where $\lambda$ is the equivariant parameter. We can alternatively define the class $\mathbf{e}_\lambda$ on line bundles by $\mathbf{e}_\lambda(L)=c_1(L)+\lambda$, and extending via the splitting principle. If the coefficient algebra is $\mathbb{C}((\lambda))$, this is an invertible characteristic class, so  it can be used to define twistings. Taking the non-equivariant limit $\lambda\to 0$ recovers the usual Euler class. 

Let $g=0$ and let $E$ be a convex vector bundle, in this instance, the $(\mathbf{e}_\lambda,E)$ twisted invariants of $X$ have the following interpretation. Let $Y\subset X$ be a complete intersection cut out by $\text{rk}(E)$ generic sections of $E$. $\pi_*E$ is a vector bundle, and it is in fact the normal bundle of $\overline{M}_{0,n,d|_Y}(Y)$ inside $\overline{M}_{0,n,d}(X)$,  

As such, in the non-equivariant limit, $(\mathbf{e}_\lambda,E)$-twisted Gromov-Witten invariants become Gromov-Witten invariants of $Y$, with classes and curve degrees inherited from the ambient space $X$. 

Let $\mathbf{e}^*_\lambda(E)$ be the class defined by $\mathbf{e}^*_\lambda(E^\vee)\mathbf{e}_\lambda(E)=1$. For $g=0$ and $E$ a concave bundle (i.e. $E^\vee$ is convex), the $(\mathbf{e}^*_\lambda,E^\vee)$-twisted invariants of $X$ have the following interpretation. $\overline{M}_{0,n,d}(X)$ embeds into $\overline{M}_{0,n,d}(E^\vee)$ as the $\mathbb{C}^*$-fixed locus. The moduli stack $\overline{M}_{0,n,d}(E^\vee)$ is not complete, however invariants of $E^\vee$ can be defined using the right hand-side of the Atiyah-Bott fixed point localization formula, which divides by the Euler class of the virtual normal bundle. So $(\mathbf{e}^{-1}_\lambda,E^\vee)$-twisted invariants can be thought of as the equivariant GW invariants of the noncompact space $E^\vee$. The twisting class defined by $(\mathbf{e}^{-1}_\lambda,E^\vee)$ differs from the previous class by a sign, which can be corrected for.

Iritani-Mann-Mignon show that in this setting, both non-equivariant limits can be taken, and are well-defined (Prop 3.6 and 3.11 of \cite{MR3519131}).

\begin{theorem}\hfill
	\begin{itemize}
		\item[(1)] For $E$ a convex bundle, $QDM_{(\mathbf{e}_\lambda,E)}$ admits a non-equivariant limit $QDM_{(\mathbf{e},E)}$. Restriction of classes to the complete intersection $Y$ defined by sections of $E$ defines a morphism of quantum $D$-modules whose target is the (untwisted) quantum $D$-module of $Y$, but defined with the Novikov variables from the ambient space. 
		
		\item[(2)] In the same situation, the quantum $D$-module of the total space of $E^\vee$ admits a  non-equivariant limit, denoted $QDM(E^\vee)$. Let $h$ be the map $\tau\mapsto \tau +\pi i c_1(E)$. Then the non-equivariant limit of $\nabla^{(\mathbf{e}_\lambda^*,E^\vee)}$ coincides with $h^*\nabla^{E^\vee}$. 
		
		\item[(3)] The map $f_\circ$ has a non-equivariant limit, denoted $f$. It is expressed in terms of 2-point euler-twisted invariants, and in the case of $E=K_V^{-1}$, is given in (\ref{eq-twisted_gf}). 
	\end{itemize}

\end{theorem}

As a result, we can state a non-equivariant version of quantum Serre duality for quantum $D$-modules (3.14 of \cite{MR3519131}):

\begin{theorem}\hfill
	\begin{itemize}
		\item[(1)] $QS$ is flat in $\tau$ and of weight $-\operatorname{dim} X$. 
		\item[(2)] Multiplication by $\mathbf{e}(E)$ is map between $QDM_{\mathbf{e},E}$ and $h\circ f^*QDM(E^\vee)$ which respects the connection and has weight $\text{rk}(E)$.
		%\item There is more stuff about fundamental solutions that we may or may not need.
	\end{itemize}
\end{theorem}

	%\appendixaddress{}{}

\end{document}